\input amstex
\mag=\magstep1
\documentstyle{amsppt}
\NoBlackBoxes

\def\fbox#1#2{\vbox to 0pt{\hrule\hbox{\hsize=#1\vrule\kern2pt
  \vbox{\kern2pt#2\kern2pt}\kern2pt\vrule}\hrule}}

\def \-#1{{#1}^{-1}}
\def \supp{\text{Supp} \ }

\def \t#1{\tilde{#1}}
\def \+n#1{{#1}^{+\nu}}
\def \o#1{\overline{#1}}
\def \ex#1#2#3#4#5#6
{0\to \Cal O_{#1}(#2) \to \Cal O_{#3}(#4) \to \Cal O_{#5}(#6) \to 0}
\def \tea#1{\text{#1} \ }
\def \teb#1{\ \text{#1} \ }
 

\topmatter
\title
On classification of $\Bbb Q$-Fano $3$-folds with Gorenstein index $2$
and Fano index $\frac 12$
\endtitle
\rightheadtext {$\Bbb Q$-Fano $3$-folds}
\leftheadtext {Hiromichi Takagi}
\author Hiromichi Takagi$^*$ 
\endauthor
\rightheadtext {$\Bbb Q$-Fano $3$-folds}
\leftheadtext {Hiromichi Takagi}
\address 
Department of Mathematical Science,
University of Tokyo,
Komaba,
Meguro-ku,
Tokyo 153-0041,
Japan\endaddress
\email
htakagi\@ms.u-tokyo.ac.jp
\endemail
\thanks 
\endthanks

\abstract
We generalize the theory developed by K. Takeuchi in [T1]
and restrict the birational type of
a $\Bbb Q$-factorial $\Bbb Q$-Fano $3$-fold $X$ with the following properties:
\roster 
\item the Picard number of $X$ is $1$;
\item the Gorenstein index of $X$ is $2$;
\item the Fano index of $X$ is $\frac 12$;
\item $h^0 (-K_X) \geq 4$;
\item there exists an index $2$ point $P$ such that
$$(X,P)\simeq (\{xy +f(z^2,u)=0\} / \Bbb Z _2 (1,1,1,0), o)$$
with $\text{ord}f(Z,U)=1$.

\endroster  

This gives an effective bound of the value of $(-K_X)^3$ and $h^0(-K_X)$
for a $\Bbb Q$-factorial $\Bbb Q$-Fano $3$-fold $X$ with 
(1)$\sim$(4) by a deformation theoretic result of T. Minagawa in [Mi2].

Furthermore based on the main result,  
we prove that 
if $X$ is a $\Bbb Q$-factorial 
$\Bbb Q$-Fano $3$-fold with (1)$\sim$(4) and with only
$\frac 12(1,1,1)$-singularities as its non Gorenstein points,
then
\roster
\item
$|-K_X|$ has a member with only canonical singularities;
\item
for any $\frac 12 (1,1,1)$-singularity,
there is a smooth rational curve $l$ through it such that $-K_X.l=\frac 12$;
\item
by a blow up at a $\frac 12(1,1,1)$-singularity,
a flopping contraction and a smoothing of Gorenstein singularities,
$X$ can be transformed to a $\Bbb Q$-Fano $3$-fold with
(1)$\sim$(4) and 
with only $\frac 12(1,1,1)$-singularities as its singularities;
\item 
$X$ can be embedded into
a weighted projective space $\Bbb P(1^h, 2^N)$, where $h:= h^0 (-K_X)$
and $N$ is the number of $\frac 12(1,1,1)$-singularities on $X$.
\endroster
\endabstract
\endtopmatter
\document

\head 0. Introduction  
\endhead

We start by the definition of $\Bbb Q$-Fano variety.
\definition{ Definition 0.0 ($\Bbb Q$-Fano variety)}
Let $X$ be a normal projective variety.
We say that $X$ is a $\Bbb Q$-Fano variety (resp. weak $\Bbb Q$-Fano variety)
if $X$ has only terminal singularities and $-K_X$ is ample (resp. nef and big).

Let $I(X):=\min \{I| IK_X \ \text{is a Cartier divisor} \}$
and we call $I(X)$ the Gorenstein index of $X$.

Write $I(X)(-K_X)\equiv r(X)H(X)$, where $H(X)$ is a primitive Cartier divisor
and $r(X) \in \Bbb N$.
(Note that $H(X)$ is unique since $\text{Pic} X$ is torsion free.)
Then we call $\frac {r(X)} {I(X)}$ the Fano index of $X$ and denote it by 
$F(X)$.
\enddefinition

\definition{Remark}
\roster
\item
We can allow that a $\Bbb Q$-Fano variety (resp. a weak $\Bbb Q$-Fano variety)
has worse singularities than terminal.
When we have to treat such a variety in this paper,
we indicate singularities which we allow, e.g.,
'a $\Bbb Q$-Fano $3$-fold with only canonical singularities';
\item
if $X$ is Gorenstein in Definition 0.0, we say that $X$ is a Fano variety
(resp. a weak Fano variety).
\endroster
\enddefinition

As an output of the minimal model program, a $\Bbb Q$-factorial 
$\Bbb Q$-Fano variety with Picard number $1$ is an important class.
We are interested in the classification of 
$\Bbb Q$-factorial $\Bbb Q$-Fano $3$-folds with Picard number $1$.
Here we mention the known result about the classification of 
$\Bbb Q$-Fano $3$-folds.
G. Fano started the classification of 
smooth Fano $3$-folds 
and it was completed by V. A. Iskovskih, V. V. Shokurov, T. Fujita,
S. Mori and S. Mukai.
S. Mukai considered Gorenstein Fano $3$-folds with canonical singularities
and classified them under mild assumptions.
In non Gorenstein case,
if Fano index is greater than $1$, 
then the classification was obtained by T.Sano [San1] and independently by
F. Campana and H. Flenner [CF] and  
if Fano index is $1$ and 
only cyclic quotient terminal singularity is allowed,
the classification was obtained by T.Sano [San2]
(recently  
T. Minagawa [Mi1] proved that any non Gorenstein $\Bbb Q$-Fano $3$-fold
with Fano index $1$ can be deformed to one with only cyclic quotient terminal
singularities).
But if Fano index is less than $1$,
the only systematic result is the classification of 
$\Bbb Q$-Fano $3$-folds which are weighted complete intersections of 
codimension $1$ or $2$ by A. R. Fletcher [Fl].  

In [T1], Kiyohiko Takeuchi developed a theory to give a simple way
of restricting birational type of a Fano $3$-fold $X$ 
with $\rho(X)=1$ and $F(X)=1$ and derived the bound of the genus of $X$
and the existence of lines.
In this paper, we formulate a generalization of Takeuchi's theory
for a $\Bbb Q$-factorial 
$\Bbb Q$-Fano $3$-fold $X$ with $\rho(X) =1$.
We expect that it is useful for the classification
of $\Bbb Q$-factorial $\Bbb Q$-Fano $3$-folds with $\rho (X)=1$ and
$F(X)<1$.
As a test case we show that it works well under the additional assumptions
that $I(X)=2$, $F(X)=\frac 12$, $(-K_X)^3 \geq 1$ and $h^0 (-K_X)\geq 1$.

Here we explain a generalization of Takeuchi's theory.
Let $X$ be a $\Bbb Q$-factorial $\Bbb Q$-Fano $3$-fold with $\rho(X)=1$.
First we seek a birational morphism $f:Y\to X$ with the following properties:
\roster
\item
$Y$ is a weak $\Bbb Q$-Fano $3$-fold;
\item
$f$ is an extremal divisorial contraction such that
$f$-exceptional divisor is a prime $\Bbb Q$-Cartier divisor. 
\endroster

Fix a $f:Y\to X$ as above and let $E$ be the exceptional divisor of $f$. 
Then we obtain the following diagram:
$$\matrix
\ & Y & \dashrightarrow & {Y'} \\
\ & {f\swarrow} & \ &  {\searrow f'} \\
X & \ & \ & \ & X' & ,
\endmatrix $$
where $Y \dashrightarrow Y'$ is an isomorphism or a composition
of possibly one flop and flips and
$f'$ is a non small contraction.
If $Y\dashrightarrow Y'$ is not an isomorphism,
let $$Y_0:=Y \overset g_0 \to \dashrightarrow Y_1 \dots
\overset g_{k-1} \to \dashrightarrow Y_k:=Y'$$ be the decomposition
of $Y \dashrightarrow Y'$ into flops and flips.
We will see that if there is a flop while $Y\dashrightarrow Y'$,
then it is $g_0$. Let $E_i$ be the strict transform of $E$ on $Y_i$
and $\t{E}$ the strict transform of $E$ on $Y'$.
Note the following:
\roster
\item
The values of $(-K_Y)^2 E$, $(-K_Y) E^2$ and $E^3$ are given. 
We know the value of $(-K_{Y_i})^3$, $(-K_{Y_i})^2 E_i$, 
$(-K_{Y_i}) {E_i}^2$ and ${E_i}^3$
are decreased by $f$, (possibly one) flop and flips and we can
express how they are decreased with some unknown quantities
associated to flop and flips and so do we the value
$(-K_{Y'})^3$, $(-K_{Y'})^2 \t{E}$, $(-K_{Y'}) \t{E}^2$ and $\t{E}^3$
with such quantities and $(-K_X)^3$.
\item
On the other hand
the value or the relation of the value (expressed by $z$ and $u$ below) 
of $(-K_{Y'})^3$, $(-K_{Y'})^2 \t{E}$, $(-K_{Y'}) \t{E}^2$ and $\t{E}^3$
are restricted by the properties of $f'$.
\endroster
By these (1) and (2), we obtain equations of Diophantine type. 
By solving these equalities, we can derive various properties of $X$
(see $\S 6$, $\S 7$ and $\S 8$).

In this paper we solve the equations in the following case:

\proclaim{Main Theorem (see Theorem 5.0 and Section 6)}
Let $X$ be a $\Bbb Q$-factorial
$\Bbb Q$-Fano $3$-fold with the following properties:
\roster 
\item $\rho(X)=1$;
\item $I(X)=2$;
\item $F(X)=\frac 12$;
\item $h^0 (-K_X) \geq 4$;
\item there exists an index $2$ point $P$ such that
$$(X,P)\simeq (\{xy +f(z^2,u)=0\} / \Bbb Z _2 (1,1,1,0), o)$$
with $\text{ord}f(Z,U)=1$.
\endroster
Then $X$ is isomorphic to one in the following table:
\endproclaim

\definition{Notation for the tables}
Let $f:Y\to X$ 
be the weighted blow up at $P$ with weight
$\frac 12 (1,1,1,2)$.
We will also use the notation as above explanation.

$$h:=h^0(-K_X).$$
$N:=\text{aw}(X)$ (see Definition 1.1 for the definition of $\text{aw}(X)$).

$e$ is defined as follows:

If there is a flop while $Y\dashrightarrow Y'$, then $e:= E^3 - {E_1}^3$,
where $E_1$ is the strict transform of $E$ on $Y_1$
(we will know that $e>0$).
If there is no flop while $Y\dashrightarrow Y'$, then $e:=0$.

$n:=\sum \text{aw}(Y_i,P_{ij})$, where the summation is taken over
the index $2$ points on flipping curves.

$z$ and $u$ is defined as follows:
If $f'$ is birational, then let $E'$ be the exceptional divisor of $f'$
and set $E' \equiv z(-K_{Y'}) - u\t{E}$ or
if $f'$ is not birational, then let $L$ be the pull back of an ample generator
of $\text{Pic} X'$
and set $L \equiv z(-K_{Y'}) - u\t{E}$,
where $\t{E}$ is the strict transform of $E$ on $Y'$.

In case $f'$ is of type $E_1$, then
let $C:=f' (E')$. 
\comment and $\o{C} \to C$ the normalization. 
$m:=\frac{8(1-g(\o{C}))-(-K_E)^2}2$ which is found to be a non negative
integer (see Proposition 2.2).
\endcomment
$$l_C := (-K_{X'}.C).$$

In case $f'$ is of type $C$,
let $\Delta$ be the discriminant divisor of $f'$.

In case $f'$ is of type $D$,
let $F$ be a general fiber of $f'$.

$Q_3$ means the smooth $3$-dimensional quadric.

$B_i$ ($1\leq i \leq 5$) 
means the $\Bbb Q$-factorial Gorenstein terminal Fano $3$-fold $X$
with $\rho(X)=1$, $F(X)=2$ and $(-K_X)^3 =8i$.

$V_{2i}$ ($1\leq i \leq 11$ and $i\not = 10$) 
means the $\Bbb Q$-factorial Gorenstein terminal Fano $3$-fold $X$
with $\rho(X)=1$, $F(X)=1$ and $(-K_X)^3 =2i$. 

Type [i] means the $\Bbb Q$-Fano $3$-fold of type [i]
which was classified by T.Sano in [San2].

The mark $\bigcirc$ means that there is an example.
\enddefinition

$$\vbox{
\offinterlineskip
\halign{\strut#&&\vrule#&\quad\hfil#\hfil\quad\cr
\multispan7 \hfil $h=4$ \hfil \cr
\noalign{\hrule}
&&exists ?&&$(-K_X)^3$&&$N$&&$e$&&$n$&&$z$&&$l_C$&&$f', X'$&\cr
\noalign{\hrule}
&&$\bigcirc$&&$\frac52$&&$1$&&$15$&&$0$&&$1$&&/&&$E_5,
(-K_{X'})^3=\frac52,I(X')=2$&\cr
\noalign{\hrule}
&&?&&$\frac 52$&&$1$&&$15$&&$0$&&$1$&&/&&crep. div., $(-K_{X'})^3=2,
I(X')=1$&\cr
\noalign{\hrule}
&&$\bigcirc$&&$3$&&$2$&&$12$&&$1$&&$1$&&/&&$E_9$, $V_4$&\cr
\noalign{\hrule}
&& &&$\frac72$&&$3$&&$10$&&$0$&&$1$&&$1$&&$E_1$, $V_6$&\cr
\noalign{\hrule}
&& &&$4$&&$4$&&$8$&&$0$&&$1$&&$2$&&$E_1$, $V_8$&\cr
\noalign{\hrule}
&& &&$4$&&$4$&&$9$&&$3$&&$1$&&/&&$E_2$, $V_{10}$&\cr
\noalign{\hrule}
&& &&$\frac 92$&&$5$&&$6$&&$0$&&$1$&&$3$&&$E_1$, $V_{10}$&\cr
\noalign{\hrule}
&& &&$\frac 92$&&$5$&&$12$&&$3$&&$1$&&/&&$E_6$, $V_{16}$&\cr
\noalign{\hrule}
&& &&$\frac 92$&&$5$&&$9$&&$0$&&$2$&&/&&$D,\deg  F =6$&\cr
\noalign{\hrule}
&& &&$5$&&$6$&&$4$&&$0$&&$1$&&$4$&&$E_1$, $V_{12}$&\cr
\noalign{\hrule}
&& &&$5$&&$6$&&$8$&&$1$&&$2$&&/&&$D,\deg F=8$&\cr
\noalign{\hrule}
}}$$

$$z=u.$$
\comment
$$g=0 \teb{in case} f' \teb{is of type} E_1.$$
\endcomment

$$\vbox{
\offinterlineskip
\halign{\strut#&&\vrule#&\quad\hfil#\hfil\quad\cr
\multispan7 \hfil $h=5$ \hfil \cr
\noalign{\hrule}
&&exists ?&&$(-K_X)^3$&&$N$&&$e$&&$n$&&$z$&&$\deg \Delta$&& 
$\deg \ F$&& $f'$ &\cr
\noalign{\hrule}
&&$\bigcirc$&&$\frac92$&&$1$&&$9$&&$0$&&$1$&&/&&$3$&&$D$&\cr
\noalign{\hrule}
&& &&$5$&&$2$&&$8$&&$1$&&$1$&&/&&$4$&&$D$&\cr
\noalign{\hrule}
&& &&$\frac{11}2$&&$3$&&$7$&&$2$&&$1$&&/&&$5$&&$D$&\cr
\noalign{\hrule}
&& &&$\frac{11}2$&&$3$&&$8$&&$0$&&$2$&&$8$&&/&&$C, \Bbb F_{2,0}$&\cr
\noalign{\hrule}
&& &&$6$&&$4$&&$7$&&$1$&&$2$&&$6$&&/&&$C, \Bbb F_{2,0}$&\cr
\noalign{\hrule}
&& &&$6$&&$4$&&$6$&&$3$&&$1$&&/&&$6$&&$D$&\cr
\noalign{\hrule}
&& &&$\frac{13}2$&&$5$&&$6$&&$2$&&$2$&&$4$&&/&&$C, \Bbb F_{2,0}$&\cr
\noalign{\hrule}
&& &&$7$&&$6$&&$5$&&$3$&&$2$&&$2$&&/&&$C, \Bbb F_{2,0}$&\cr
\noalign{\hrule}
&& &&$\frac{15}2$&&$7$&&$4$&&$4$&&$2$&&$0$&&/&&$C, \Bbb F_{2,0}$&\cr
\noalign{\hrule}
}}$$

$$z=u.$$

$$\vbox{
\offinterlineskip
\halign{\strut#&&\vrule#&\quad\hfil#\hfil\quad\cr
\multispan7 \hfil $h=6$ \hfil \cr
\noalign{\hrule}
&&exists ?&&$(-K_X)^3$&&$N$&&$e$&&$n$&&$z$&&$\deg \Delta$&&$l_C$&&
$f', X'$&\cr
\noalign{\hrule}
&&$\bigcirc$&&$\frac{13}2$&&$1$&&$7$&&$0$&&$1$&&$7$&&/&&$C, \Bbb P^2$&\cr
\noalign{\hrule}
&&$\bigcirc$&&$7$&&$2$&&$7$&&$0$&&$4$&&/&&$35$&&$E_1$,$[5]$&\cr
\noalign{\hrule}
&& &&$7$&&$2$&&$6$&&$1$&&$1$&&$6$&&/&&$C,\Bbb P^2$&\cr
\noalign{\hrule}
&&$\bigcirc$&&$\frac{15}2$&&$3$&&$7$&&$0$&&$2$&&/&&$9$&&$E_1,[2]$&\cr
\noalign{\hrule}
&&$\bigcirc$&&$\frac{15}2$&&$3$&&$6$&&$1$&&$4$&&/&&$30$&&$E_1$,$[5]$&\cr
\noalign{\hrule}
&& &&$\frac{15}2$&&$3$&&$5$&&$2$&&$1$&&$5$&&/&&$C,\Bbb P^2$&\cr
\noalign{\hrule}
&& &&$8$&&$4$&&$4$&&$3$&&$1$&&$4$&&/&&$C,\Bbb P^2$&\cr
\noalign{\hrule}
&& &&$\frac{17}2$&&$5$&&$3$&&$4$&&$1$&&$3$&&/&&$C,\Bbb P^2$&\cr
\noalign{\hrule}
&& &&$9$&&$6$&&$2$&&$5$&&$1$&&$2$&&/&&$C,\Bbb P^2$&\cr
\noalign{\hrule}
&& &&$\frac{19}2$&&$7$&&$1$&&$6$&&$1$&&$1$&&/&&$C,\Bbb P^2$&\cr
\noalign{\hrule}
}}$$
\comment
$$m=0 \teb{and} u=z+1 \teb{in case} f' \teb{is of type} E_1.$$
\endcomment
$$u=z \teb{in case} f' \teb{is of type} C.$$
 
$$\vbox{
\offinterlineskip
\halign{\strut#&&\vrule#&\quad\hfil#\hfil\quad\cr
\multispan7 \hfil $h=7$ \hfil \cr
\noalign{\hrule}
&&exists ?&&$(-K_X)^3$&&$N$&&$e$&&$n$&&$z$&&$l_C$&&$f', X'$&\cr
\noalign{\hrule}
&&$\bigcirc$&&$\frac{17}2$&&$1$&&$6$&&$0$&&$3$&&$36$&&$E_1,\Bbb P^3$&\cr
\noalign{\hrule}
&&$\bigcirc$&&$9$ 
&&$2$&&$6$&&$0$&&$2$&&$18$&&$E_1$,$[3]$&\cr
\noalign{\hrule}
&&$\bigcirc$&&$9$&&$2$&&$5$&&$1$&&$3$&&$32$&&$E_1,\Bbb P^3$&\cr
\noalign{\hrule}
&&$\bigcirc$&&$\frac{19}2$&&$3$&&$5$&&$1$&&$2$&&$15$&&$E_1$,$[3]$&\cr
\noalign{\hrule}
&&$\bigcirc$&&$\frac{19}2$&&$3$&&$4$&&$2$&&$3$&&$28$&&$E_1,\Bbb P^3$&\cr
\noalign{\hrule}
}}$$

$$u=z+1.$$
\comment $m=0$ \endcomment

$$\vbox{
\offinterlineskip
\halign{\strut#&&\vrule#&\quad\hfil#\hfil\quad\cr
\multispan7 \hfil $h=8$ \hfil \cr
\noalign{\hrule}
&&exists ?&&$(-K_X)^3$&&$N$&&$e$&&$n$&&$z$&&$l_C$&&$f, X'$&\cr
\noalign{\hrule}
&&$\bigcirc$&&$\frac{21}2$&&$1$&&$6$&&$0$&&$1$&&$6$&&$E_1, B_3$&\cr
\noalign{\hrule}
&&$\bigcirc$&&$\frac{21}2$&&$1$&&$5$&&$0$&&$2$&&$27$&&$E_1, Q_3$&\cr
\noalign{\hrule}
&&$\bigcirc$&&$11$&&$2$&&$4$&&$1$&&$2$&&$24$&&$E_1, Q_3$&\cr
\noalign{\hrule}
&& &&$\frac{23}2$&&$3$&&$3$&&$2$&&$2$&&$21$&&$E_1, Q_3$&\cr
\noalign{\hrule}
}}$$

$$u=z+1.$$
\comment $m=0$ \endcomment

$$\vbox{
\offinterlineskip
\halign{\strut#&&\vrule#&\quad\hfil#\hfil\quad\cr
\multispan7 \hfil $h=9$ \hfil \cr
\noalign{\hrule}
&&exists ?&&$(-K_X)^3$&&$N$&&$e$&&$n$&&$z$&&$u$&&$l_C$&&$f', X'$&\cr
\noalign{\hrule}
&&$\bigcirc$
&&$\frac{25}2$&&$1$&&$5$&&$0$&&$1$&&$2$&&$10$&&$E_1,B_4$&\cr
\noalign{\hrule}
}}$$

$$\vbox{
\offinterlineskip
\halign{\strut#&&\vrule#&\quad\hfil#\hfil\quad\cr
\multispan7 \hfil $h=10$ \hfil \cr
\noalign{\hrule}
&&exists ?&&$(-K_X)^3$&&$N$&&$e$&&$n$&&$\deg \Delta$&&$l_C$&&$f', X'$&\cr
\noalign{\hrule}
&&$\bigcirc$&&$\frac{29}2$&&$1$&&$4$&&$0$&&/&&$14$&&$E_1,B_5$&\cr
\noalign{\hrule}
&& &&$\frac{29}2$&&$1$&&$6$&&$0$&&$0$&&/&&$C, \Bbb P^2$&\cr
\noalign{\hrule}
&&$\bigcirc$&&$15$&&$2$&&$3$&&$1$&&/&&$12$&&$E_1, B_5$&\cr
\noalign{\hrule}
}}$$

$$z=1\teb{and} u=2.$$
\comment
$$m=0 \teb{in case} f' \teb{is of type} E_1.$$
\endcomment
\qed

Based on this result, we can derive the following properties for 
$X$ as in the main theorem if $X$ has 
only $\frac 12(1,1,1)$-singularities as its non Gorenstein points:
 
\proclaim{Theorem A (See Corollary 7.2)}
$|-K_X|$ has a member with only canonical singularities.
\endproclaim

So the general elephant conjecture by Miles Reid 
is affirmative for $X$.

\proclaim{Theorem B (See Corollary 8.1)}
Let $f:Y \to X$ be as in the main theorem,
i.e., $f$ is the blow up at a $\frac 12(1,1,1)$-singularity
and
$g: Y\to Z$ the anti-canonical model
(by the main theorem, $g$ is found to be not an isomorphism).
Then if $N>1$ (resp. $N=1$), 
$Z$ can be deformed to a $\Bbb Q$-Fano $3$-fold $Z'$ 
with $\rho (Z')=1$ and $F(Z')=\frac 12$
which has only $N-1$ $\frac 12 (1,1,1)$-singularities as its singularities
and $h^0 (-K_{Z'}) = h$ (resp. a smooth Fano $3$-fold $Z'$ 
with $\rho (Z')=1$, $F(Z')=1$ and $h^0 (-K_{Z'}) = h$.)
\endproclaim
 
This is an analogue to the Reid's fantasy about Calabi-Yau $3$-folds [RM3].

\proclaim{Theorem C (See Corollary 8.3)}
$X$ can be embedded into
a weighted projective space $\Bbb P(1^h, 2^N)$, where $h:= h^0 (-K_X)$
and $N$ is the number of $\frac 12(1,1,1)$-singularities on $X$.
\endproclaim

We hope that this fact can be used for the classification of Mukai's
type (see [Mu1], [Mu2] and [Mu3]).

\definition{Acknowledgment}
I express my hearty thanks to Professor Kiyohiko Takeuchi for
giving me several examples of $\Bbb Q$-Fano $3$-folds.
I also thank for his theory developed in [T1].
Without these, this paper could not be born.
I am grateful to Professor Shigefumi Mori for giving me many useful
comments.
I am thankful to Professor Keiji Oguiso whose lecture at University of Tokyo
is helpful for me. 
I also thank Professor Yujiro Kawamata for giving me useful comments.
I also thank Mr. Tatsuhiro Minagawa for stimulus discussions on this subject.
I express my indebtedness to Professor Keiichi Watanabe
for giving me useful comments about graded rings of $\Bbb Q$-Fano $3$-folds. 
I am thankful to Professor Tie Luo for kindly sending me his paper [Lu]
before publication.
\enddefinition

\definition{Notation and Conventions}
\roster
\item In this paper, we will work over $\Bbb C$, the complex number field;
\item we denote the linear equivalence by $\sim$
and the numerical equivalence by $\equiv$.
The equality $=$ in an adjunction formula means the $\Bbb Q$-linear
equivalence;
\item we denote the Hirzebruch surface of degree $n$ by $\Bbb F_n$
and the surface which is obtained by the contraction of the negative section
of $\Bbb F_n$ by $\Bbb F_{n,0}$.
\endroster
\enddefinition

\head 1. Preliminaries \endhead

\proclaim{Theorem 1.0 (Vanishing theorem)}
Let $f: X\to Y$ be a projective morphism from a normal variety
$X$ with only Kawamata log terminal singularities.
Let $D$ be a $\Bbb Q$-Cartier integral Weil divisor
such that $D-K_X$ is $f$-nef and $f$-big.
Then $R^i f_* \Cal O_X (D) =0$ for all $i>0$.

We will quote this theorem as KKV vanishing theorem.
\endproclaim

\demo{Proof} 
See [Kod1], [KY1] and [V].
\qed
\enddemo

\definition{Definition 1.1}
Let $(X,P)$ be a germ of $3$-dimensional terminal singularity of index $>1$.
By the classification of such a singularity [Mo2], we can easily see
that a general deformation of $(X,P)$ has only cyclic quotient singularities. 
We call the number of these cyclic quotient singularities the axial weight
of $(X,P)$ and denote it by aw$(X,P)$.
Let $X$ be a $3$-fold with only terminal singularities.
We define aw$(X) :=\sum$ aw$(X,P)$, where 
the summation takes place over points of index $>1$.
\enddefinition

\proclaim{Theorem 1.2 (Special case of the singular Riemann-Roch Theorem)}
Let $X$ be a $3$-fold with at worst index $2$ terminal singularities
and $D$ an integral Weil divisor on $X$.
Then the following formula holds:
$$\chi(\Cal O_X(D)) = \chi(\Cal O_X) + \frac 1{12} D(D-K_X)(2D-K_X) +
\frac 1{12} D.c_2(X) +\sum c_Q(D),$$
where the summation takes place over non Gorenstein points where $D$ is not
Cartier and  $\sum c_Q(D)=-\frac n8$ for some non negative integer $n$.
(See [RM2, Theorem 10.2] for the definition of $c_Q(D)$.) 
\endproclaim
\demo{Proof}
See [RM2, Theorem 10.2].
\qed
\enddemo

\proclaim{Theorem 1.3}
Let $X$ be a projective $3$-fold 
with at worst index $2$ terminal singularities.
Then $-K_X.c_2(X)= 24 - \frac {3N}2$,
where $N:=$aw$(X)$.
Furthermore assume that $X$ is a $\Bbb Q$-factorial 
$\Bbb Q$-Fano $3$-fold with $\rho(X)=1$. 
Then $-K_X.c_2(X) >0$.
In particular $N\leq 15$.
\endproclaim
\demo{Proof}
See [KY2, Lemma 2.2 and Lemma 2.3] and [KY3, Proposition 1].
\qed
\enddemo

\proclaim{Corollary 1.4}
Let $X$ be a weak $\Bbb Q$-Fano $3$-fold with $I(X)=2$.
Then $h^0(-K_X) = 3 + \frac 12 (-K_X)^3 - \frac N4$,
where $N:=$aw$(X)$.
\endproclaim
\demo{Proof}
This directly follows from Theorem 1.0, Theorem 1.2 and Theorem 1.3.
\qed
\enddemo

\proclaim{Proposition 1.5}
Let $f:X \to (Y,Q)$ be a flopping contraction 
from a $3$-fold $X$ with only terminal singularities to a germ 
$(Y,Q)$ and $f^+ : X^+ \to Y$ the flop of $f$ constructed as
in [Kol1, Theorem 2.4].
Let $D$ be a Cartier divisor on $X$
and $D^+$ the strict transform of $D$ on $X^+$.
Then $D^+$ is a Cartier divisor.
\endproclaim

\demo{Proof}
By passing to the analytic category and 
taking algebraization [Ar, Theorem 3.8],
we may assume that $C:= \text{excep}f$ is irreducible.
Furthermore since we can deform $X$ to a $3$-fold with only cyclic quotient
terminal singularities [Mo3, (1b.8.2) Corollary] and 
such a deformation lifts to one of $f:X\to Y$ [KoMo, (11.4) Proposition],
we may assume that $X$ has only cyclic quotient
terminal singularities.
\comment
Let $D'$ be an effective divisor on $X$ such that $D'.C=1$ and $D'$
and $C$ intersect at a smooth point of $X$.
By $\text{Pic} X \simeq \text{Pic} C$, $D\sim (D.C) D'$.
Since linear equivalence is preserved by the flop, it suffices to prove
the assertion for $D'$. So we may assume that $D$ is effective and $D.C=1$. 
\endcomment
Let $H'$ be a general hyperplane section through $Q$ and $H:=f^* H'$. 
Then it is well known 
that 

(1.5.1) $H'$ and $H$ have only canonical singularities
and $H$ is dominated by the minimal resolution of $H'$.

We show that there are at most $2$ singularities on $C$.
Assume the contrary. Then there are $3$ singularities on $C$
and they coincide singularities of $H$ on $C$ by (1.5.1).
Let $p:\t{Y} \to Y$ be the canonical cover,
$\t{X} :=X\times _Y \t{Y}$, $\t{C}$ (resp. $\t{H'}$, $\t{H}$) 
the pull back of $C$ (resp. $H'$, $H$) on $\t{X}$ 
and $\t{f} : \t{X}\to \t{Y}$ the induced morphism.
Then $\t{X}$ is smooth
and $\t{f}$ is also a flopping contraction.
We will prove that $\t{C}$ is irreducible. 
If $Q$ is not of the exceptional type ([RM2, Theorem (6.1) (2)) and
$\t{C}$ is reducible, then there are components which intersect
at $3$ points, a contradiction to $R^1 \t{f} _* \Cal O _{\t{X}} =0$.
Hence $\t{C}$ is irreducible in this case.
If $Q$ is of the exceptional type and $\t{C}$ is reducible,
then $\t{C}$ has two component $\t{C_1}$, $\t{C_2}$
and they intersect at one point transversely.
But $\t{C_i} \to C$ is a double cover between $\Bbb P^1$'s and 
is branched at three points, a contradiction to Hurwitz's formula.
Hence in any case, $\t{C}$ is irreducible. 
By [RM2, (4.10)], $\t{H}$ must be smooth.
Hence $\t{H'}$ has only ODP whence $H'$ has a canonical singularity
of type A. But then $H$ has at most $2$ singularities, a contradiction.
So we have the assertion.  

Furthermore $H$ has exactly two singularities. For otherwise
$\text{aw}(Y,Q)=1$ since $\text{aw}(Y,Q)=\text{aw}(X)$.
Hence $Q$ is a cyclic quotient singularity but then there is no flopping
contraction to $Q$, a contradiction.

We can prove as above 
that $\t{C}$ is irreducible if $Q$ is not of the exceptional type
or $\t{C}$ has at most 2 components if $Q$ is of the exceptional type.

Assume that $Q$ is not of the exceptional type. Let $r$ be the index of $Q$.
Let $P$ be a non Gorenstein point on $C$ and $\t{P}$ the inverse image on 
$\t{X}$. Then $P$ is also of index $r$ and by [RM2, (4.10)], we have
locally analytically
$$(\t{P} \in \t{C} \subset \t{X}) \simeq 
(o \in \{x=y=0 \} \subset \Bbb C^3),$$
where $x,y,z$ are coordinates of $\Bbb C^3$
which are semi-invariants of $\Bbb Z _r$-action. 
Let $\t{E}$ be a Cartier divisor which is localized to
$z=0$ and $E$ the image of $\t{E}$ on $X$. Then we have
$E.C=\frac 1r$.
Since $rE$ is a Cartier divisor and $\text{Pic} X \simeq \text{Pic} C$,
we have $D\sim r(D.C)E$.
Then we have $D^+ \sim r(D.C)E^+$, 
where $E^+$ is the strict transform of $E$ on $X^+$
because linear equivalence is preserve by a flop.
Since the analytic types of $X$ and $X^+$ are the same by [Kol1, Theorem 2.4],
$r(D.C)E^+$ is Cartier and so is $D^+$.

Assume that $Q$ is of the exceptional type.
Then $X$ has one index $2$ point and one index $4$ point.
Using [Kol2, Proposition 2.2.6],
we can determine $Q$ as follows:
$$(Y,Q) \simeq \{x^2 -y^2 + (z^2 -u^{4k+4})(z -au^{4l+2}) \}
/ \Bbb Z _4 (1,3,2,1), o),$$
where $k$ and $l$ are non negative integers and $a \in \Bbb C$.
Let $\t{f_1}:\t{X_1} \to \t{Y}$ be the blow up 
along $\{ x+y=z+u^{2k+2}=0\}$
and $\t{f_2}:\t{X_2} \to \t{X_1}$ 
the blow up along the strict transform of $\{ x-y =z-u^{2k+2}=0 \}$.
Let $\t{F_1}$ (resp. $\t{F_2}$) be the strict transform 
of $\{ x+y=z+u^{2k+2}=0\}$ (resp. $\{ x-y =z-u^{2k+2}=0 \}$) on $\t{Y_2}$.
Then $-\t{F_1}-\t{F_2}$ is $\t{f_1}\circ \t{f_2}$-ample
and  $(\t{f_1}\circ \t{f_2})_* (\t{F_1}+\t{F_2})$ is $\Bbb Z _4$-invariant.
Hence $\Bbb Z _4$ acts on $\t{X_2}$ regularly
and $\t{f_1}\circ \t{f_2}$ is equivariant, i.e.,
we can identify $\t{f} :\t{X} \to \t{Y}$ and 
$\t{f_1}\circ \t{f_2} :\t{X_2} \to \t{Y}$.
Let $F':=p_* \{ x+y=z+u^{2k+2}=0\} _{\text{red}}$ and $F$ its strict transform
on $X$. We can see that $F$ is a Cartier divisor and $F.C =-1$.
Note that the involution $(x,y,z,u)\to (-x,y,z,u)$ induce an involution $\iota$
on $Y$. Then we have 
$F'+ \iota _* F' \sim 0$ whence $F^+$ is a Cartier divisor.
Since $D\sim -(D.C)F$ by $\text{Pic} X \simeq \text{Pic} C$,
$D^+$ is also a Cartier divisor.

\comment
By [Ko1, Theorem 2.4], $D^+.C^+ =-D.C \in \Bbb Z$, where $C^+ :=\text{f^+} C$.
Hence by $\text{Pic} X^+ \simeq \text{Pic} C^+$, there is a Cartier divisor $E$
such that $E\equiv D^+$.
We may assume that $E$ is effective.
Take a compactification ${X^+}'$ of $X^+$. 
Since $E$ and $D^+$ are effective and
$rE \sim rD^+$, where $r$ is the index of $Q$,
the numerical equivalence holds also on ${X^+}'$.
\endcomment

\qed
\enddemo

\head 2. extremal contractions 
from $3$-folds with only index $2$ terminal singularities
\endhead

\definition{Definition 2.0 (Extremal contraction)}
Let $X$ be an analytic $3$-fold with only terminal singularities
and $f: X\to (Y, Q)$ a projective morphism onto a germ of a normal variety 
with only connected fibers.
Let $\text{excep} f$ be the locus where $f$ is not isomorphic. 
Assume that $-K_X$ is $f$-ample.
\roster
\item 
If $\dim Y = 3$ and $\dim \ \text{excep} f = 1$,
then we say that $f$ is an extremal contraction of flipping type
(or in short a flipping contraction).
\item 
Only in this case, we assume that $-K_X$ is $f$-numerically trivial
instead that $-K_X$ is $f$-ample.
If $\dim Y = 3$ and $\dim \ \text{excep}  f = 1$,
then we say that $f$ is a flopping contraction.
\item 
Assume that $\dim Y=3$, $\text{excep}  f$ is purely $2$-dimensional
and every component of the exceptional divisor $E$ 
is contracted to a curve. Let $C:= f(E)$.
Assume furthermore that over a general point of every component of $C$,
$f$ coincides with the blow up along $C$ and $-E$ is $f$-ample.
Then we say that $f$ is an extremal contraction of type $E_1$.
\item 
Assume that $\dim Y=3$, $\text{excep}  f$ is an irreducible divisor
$E$ and $f(E)$ is a point. 
Then we say that $f$ is an extremal contraction of type $E_{\geq 2}$.
\item
If $\dim Y=2$ and every fiber is $1$-dimensional,
then we say that $f$ is an extremal contraction of type $C$.
\item
If $\dim Y= 1$ and $f^{-1} (Q) _{\text{red}}$ is irreducible, then
we say that $f$ is an extremal contraction of type $D$.  
\endroster
\enddefinition

\proclaim {Proposition 2.1 (Flipping contraction)}
Let $X$ be an analytic $3$-fold \newline
with only index $2$ terminal singularities
and $f: X\to (Y,Q)$ a flipping contraction to a germ $(Y,Q)$.
Let $C$ be its exceptional curve.
(Since $(Y, Q)$ is a germ, $C$ is connected.)
Then 
\roster
\item $C \simeq \Bbb P^1$ and there is only one
index $2$ singularity on $C$ and $-K_X.C= \frac 12$;
\item
let $P$ be the unique index $2$ singularity on $C$.
Then locally analytically
$(P\in C \subset X)\simeq
(o\in \{ x_2=x_3=x_4=0 \} \subset 
\{ x_1x_2 +p({x_3}^2,x_4)=0\} / \Bbb Z _2 (1,1,1,0))$;
\item 
let $p(0,x_4)=a{x_4}^k$, where $a$ is a unit in $\Bbb C\{x_1,x_2,x_3,x_4\}$
and $k\in \Bbb N$ (note that $k=$aw$(X,P)$). 
Then there is a deformation $\frak f : \frak X \to \frak Y$
of $f$ over a $1$-dimensional disk $(\Delta,0)$
such that for $t\not =0$, $\frak X_t$ has only 
$k\frac 12(1,1,1)$-singularities
and
$\frak f _t : \frak X _t \to \frak Y _t$ 
is a bimeromorphic morphism which is localized to $k$ flipping contractions.
\item
assume that $P$ is a $\frac 12(1,1,1)$-singularity.
Then we can construct the flip of $f$ as follows:

Let $g:X_1 \to X$ be the blow up of $P$ and $E_1$ the exceptional divisor.
Let $h: X_2 \to X_1$ be the blow up along the strict transform $C_1$ of $C$ 
on $X_1$
and $E_2$ the exceptional divisor.
Then $E_2 \simeq \Bbb P^1 \times \Bbb P^1$ and we can blow it down to
another direction.
Let $i:X_2 \to {X_1}^+$ be the blow down and ${E_1}^+$ the strict transform
of $E_1$ on ${X_1}^+$.
Then ${E_1}^+\simeq \Bbb F_1$ and we can blow it down to the ruling direction.
Let $j: {X_1}^+ \to X^+$ be the blow down.
Then $X\dashrightarrow X^+$ is the flip.
\item
If $X$ is projective and $f$ is an algebraic flipping contraction,
then $(-K_{X^+})^3 = (-K_X)^3 -\frac n2$, where $n=\sum$aw$(X,P)$
and the summation is taken over the non Gorenstein points on flipping curves. 
\endroster
\endproclaim

\demo{Proof}
As for (1), (2) and (4), see [KoMo, (4.2) and (4.4.5)].
We will prove (3).
Construct $Y'$ as in [ibid. (4.3)].
Then $Y'=\{ y_1y_3 +y_2 p({y_2}^2 , y_4) =0\}$ as in [ibid. (4.4.2)].
Then $f$ is obtained by blowing-up of $Y'$ along $\{ y_2=y_3=0\}$
and dividing by the $\Bbb Z_2$ action.
Let $\frak Y '=\{ y_1y_3 +y_2 (p({y_2}^2 , y_4)+ty_4) =0\}$ be 
a deformation of $Y'$ over a $1$-dimensional disk $(\Delta,0)$.
Then by blowing-up of $\frak Y '$ along $\{ y_2=y_3=0\}$
and dividing by the induced $\Bbb Z_2$ action,
we obtain the desired $\frak f$.
Next we prove (5).
If we compactify $\frak X$ in (3), then (5) holds by (4) and the invariance
of $(-K)^3$ in a flat family. Since $(-K_X)^3 - (-K_{X^+})^3$ can be
expressed by an intersection number of the pull back of $(-K_X)$ with
exceptional divisors on a simultaneous resolution of $X^+$ and $X$
(and hence it is determined locally around flipping curves),
the general case follows.
\qed
\enddemo

\proclaim {Proposition 2.2 (Contraction of type $E_1$)}
Let $X$ be an analytic $3$-fold \newline
with only index $2$ terminal singularities
and $f: X\to (Y, Q)$ an extremal contraction of type $E_1$ to a germ $(Y,Q)$.
Let $E$ be the exceptional divisor and $C:= f(E)$. 
Let $l$ be the fiber over $Q$.
Then the following holds:

\roster
\item
Assume that $l$ contains no index $2$ point.
Then $Q$ is a smooth point and $f$ is the blow up along $C$;
\item 
Assume that $l$ contains an index $2$ point.
Then $l$ contains only one index $2$ point 
(we will denote it by $P$) and 
every component $l'$ of $l$ passes through $P$ and
satisfies $-K_X.l' = \frac 12$. 
\endroster
\comment 
\item "(a)" $C$ is a smooth curve;
\item "(b)"  $(Q\in Y) \simeq (o\in ((xy+zw=0) \subset \Bbb C^4))
\ \text{or} \  (o\in ((xy+z^2 + w^3=0) \subset \Bbb C^4))$;
\item "(c)"  in the former case of (b), $l$ is a reducible conic
and in the latter case of (b), $l$ is a double line. 
In any case, for every component $l'$ of $l$,
$-K_X.l' = \frac 12$;
\endroster
\endcomment
Assume furthermore that $X$ is projective and $\rho(X/Y)=1$.
Then the following formula holds:
$$(-K_E)^2 = 8(1-g(\o{C})) -2m ,$$
where $\o{C}$ is the normalization of
$C$ and $m$ is a non-negative integer.
\endproclaim

\demo{Proof}
See [Mo1, Theorem 3.3] for (1).
\comment
In case (2), we first prove that $C$ is smooth.
Let $\t{X} :=\pmb{Spec}  (\Cal O_X \oplus \Cal O_X(K_X))$, where
we define a ring structure of $\Cal O_X \oplus \Cal O_X(K_X)$ by a
smooth general element $G$ of $|-2K_X|$. Let $\t{E}$ be the pull back of $E$.
Note that $\t{X}$ is smooth.
Then there is a natural crepant contraction of $\t{E}$ from $\t{X}$
which contracts $\t{E}$ to a curve $\t{C}\simeq C$.
Note that $\t{E}$ is negative for exceptional curves of 
the crepant contraction and 
the contraction coincides with 
the blow up of $\t{C}$ at a general point of $\t{C}$.
By these and the proof of [W1, Theorem 2.2] and [W2, Proposition 3.1],
we know that $\t{C}$ is smooth.
By this, we know that the case 1 in [KoMo, (4.8.3)] does not occur by
[KoMo, Proposition 4.10.1] and (b) and (c) follow from
[KoMo, (4.8.4) and (4.8.5)].
We also know by [ibid.] that in the former case of (b),
$\text{Sing}  E \cap l = \{ P \}$ and $P$ is an ordinary double point of $E$
and in the latter case of (b),
$\text{Sing}  E \cap l = \{ P, P' \}$ and $P, P'$ are
ordinary double points of $E$.
Hence if there is one singular point of $Y$ on $C$,
then $(-K)^2$ is less than one of a geometrically ruled surface over $C$
by $2$. So we obtain (d).
\endcomment
Assume that $X$ is projective and $\rho(X/Y)=1$.
Let $\mu :\o{E} \to E$ be the normalization
and define a $\Bbb Q$-divisor $Z$ by $K_{\o{E}}=\mu^* K_{E} -Z$.
Then $Z$ is effective and its support is contained in fibers.
Hence $Z.(-K_{\o{E}})\geq 0$ and $(-K_E)^2 \leq (-K_{\o{E}})^2
\leq 8(1-g(\o{C})).$
Since $-K_X -E \sim f^* (-K_Y) -2m$, 
$(-K_E)^2 = (-K_X -E)^2 E =2(2E^3 -2f^* (-K_Y)E^2) \in 2\Bbb Z$.
Hence we have the formula as above.
\qed
\enddemo

\proclaim {Proposition 2.3 (Contraction of type $E_{\geq 2}$)}
Let $X$ be a $3$-fold with only index $2$ terminal singularities
and $f: X\to (Y, Q)$ a divisorial contraction to a germ $(Y,Q)$
which contracts a divisor $E$ to $Q$.
Then the following holds:

\roster
\item
Assume that $E$ contains no index $2$ point.
Then one of the following holds:

$$(E_2): (E,-E|_{E}) \simeq (\Bbb P^2, \Cal O_{\Bbb P^2}(1))
\teb{and} Q \teb{is a smooth point};$$

$$(E_3): (E,-E|_{E}) \simeq (\Bbb P^1 \times \Bbb P^1, 
\Cal O_{\Bbb P^3}(1)|_{\Bbb P^1 \times \Bbb P^1}) \teb{and}
(Y,Q)\simeq (((xy+zw=0) \subset \Bbb C^4),o);$$

$$(E_4): (E,-E|_{E}) \simeq (\Bbb F_{2,0}, 
\Cal O_{\Bbb P^3}(1)|_{\Bbb F_{2,0}}) \teb{and}
(Y,Q) \simeq (((xy+z^2 + w^k=0) \subset \Bbb C^4),o)
(k\geq 3);$$

$$(E_5): (E,-E|_{E}) \simeq (\Bbb P^2, \Cal O_{\Bbb P^2}(2)) \teb{and}
Q \teb{is a} \frac 12(1,1,1)\text{-singularity}.$$

Furthermore for all cases, $f$ is the blow up of $Q$.

\item 
Assume that $E$ contains an index $2$ point.
Then one of the following holds:

$$\multline
(E_6): (E,-E|_{E}) \simeq (\Bbb F_{2,0}, l)\teb{,
where} l \teb{is a ruling of} \Bbb F_{2,0}. \\
Q \teb{is a smooth point and} f 
\teb{is a weighted blow up with weight} (2,1,1). \\
\teb{In particular we have} K_X = f^*K_Y+ 3E;
\endmultline$$

$$(E_7): K_X=f^* K_Y +E \ \text{and} \ Q \teb{is a
Gorenstein singular point.} E^3 = \frac 12;$$
$$(E_8): K_X=f^* K_Y +E \ \text{and} \ Q \teb{is a
Gorenstein singular point.} E^3 = 1;$$
$$(E_9): K_X=f^* K_Y +E \ \text{and} \ Q \teb{is a
Gorenstein singular point.} E^3 = \frac 32;$$
$$(E_{10}): K_X=f^* K_Y +E \ \text{and} \ Q \teb{is a
Gorenstein singular point.} E^3 = 2;$$

$$\multline
(E_{11}): (E,-E|_{E}) \simeq ((\{ xy+w^2 =0 \}\subset \Bbb P (1,1,2,1)),
\Cal O (2)). \\
(Y,Q)\simeq 
(((xy+z^k + w^2=0) \subset \Bbb C^4/\Bbb Z_2 (1,1,0,1)),o).\\
f \teb{is a weighted blow up with a weight} (\frac 12, \frac 12, 1,
\frac 12).\\
\teb{In particular we have} K_X = f^*K_Y+ \frac 12 E;
\endmultline$$

$$\multline
(E_{12}): (E,-E|_{E}) \simeq (\Bbb F_{2,0}, 3l).\\
Q \teb{is a} \frac 13 (2,1,1)\text{-singularity and}\ 
f \teb{is a weighted blow up
with a weight} \frac 13 (2,1,1).\\
\teb{In particular we have} K_X = f^*K_Y+ \frac 13 E;
\endmultline$$

\endroster
\endproclaim

\demo{Proof}
See [Mo1, Theorem 3.4 and Corollary 3.5], [Cu] for (1) and
[Lu, Corollary 2.5 and Theorem 2.6] for (2) and $Q$ is a
non Gorenstein point.
We will prove here that if $Q$ is a Gorenstein point,
$f$ is of type $E_6 \sim E_{10}$.
Let $a$ be the discrepancy for $E$. Since $Q$ is assumed to be Gorenstein,
$a$ is a positive integer.

First assume that $a\geq 2$. Let $L:=-2E$.
Then $L$ is free by [AW] 
since $K_X + \frac a 2 L \equiv 0$ and $\frac a 2 \geq 1$.
Let $D$ be a general member of $|L|$ and $C:=E|_D$.
Since $-K_D\equiv -(a-2)E|_D$ is nef and big, $C$ is a tree of $\Bbb P^1$
by KKV vanishing theorem.
Let $\mu : \t{E} \to E$ be the normalization of $E$. 
If $C$ is reducible, then $\mu ^* C$ is not connected, a contradiction
to the ampleness of $\mu ^ * C$. Hence $C\simeq \Bbb P^1$.
By this we know that $E$ is normal since $E$ satisfies $S_2$ condition.
Since $C$ is ample and isomorphic to $\Bbb P^1$,
$E\simeq \Bbb P^2, \Bbb F_n (n\geq 1) \ \text{or}  \ \Bbb F_{n,0} (n\geq 2)$
by a classical result (see for example [Ba]).
But if former $2$ cases occur, $X$ is smooth, a contradiction 
to the assumption of (2).
Hence $E\simeq \Bbb F_{n,0} (n\geq 2)$. 
We will prove that $n=2$.
Let $v$ be the vertex of $E$. 
Then $v$ is the unique singularity on $E$ and hence it is of index $2$.
If $E$ is Cartier at $v$, then for a exceptional divisor $F$ over $v$
with discrepancy $\frac 12$, the discrepancy of $F$ for $K_Y$ 
is not an integer, a contradiction.
Hence $K_X + E$ is a Cartier divisor and hence $K_E$ is Cartier at $v$.
So $n$ must be $2$.
Furthermore by $K_{E} = (a+1)E|_E$, $a=3$ since $a\geq 2$ and 
$E\simeq \Bbb F_{2,0}$. By taking the canonical cover near $v$ of $X$,
we know that $v$ is a $\frac 12 (1,1,1)$-singularity.
We will prove that $Q$ is smooth and $f$ is a weighted blow up
with a weight $(2,1,1)$.
We see that $\o{E}$ is contracted to a curve and let $\o{X} \to \o{X'}$
the contraction. Then next we can contract the strict transform of $F$
to a smooth point, which is no other than $Q$.
We can easily show that a weighted blow up
with a weight $(2,1,1)$ is decomposed into contractions as above.
So we are done.

Next we assume that $a=1$.
Let $P$ be an index $2$ point on $X$.
If $E$ is Cartier at $P$, then for a exceptional divisor $F$ over $P$
with discrepancy $\frac 12$, the discrepancy of $F$ for $K_Y$ 
is not an integer,
a contradiction. Hence $E$ is not Cartier at $P$ whence
$M:=-K_X -E$ is an ample Cartier divisor.
So $E$ is a Gorenstein (possibly non normal) del Pezzo surface
since $-K_E=M|_E$. Since $\chi(\Cal O_E)=1$ by [Sak, Theorem (5.1)] and
[RM5, Corollary 4.10],
$\text{Pic}  E$ is torsion free. So $-K_X +E|_E \sim 0$ and
hence $-K_X + E\sim 0$ by $\text{Pic}  X \simeq \text{Pic}  E$.
So we note that $M\sim -2K_X$.
Since $(-K_E)^2 = 4E^3 \geq 2$, $|-K_E|$ is free by 
[RM5, Corollary 4.10]
and [Fu2, Corollary 1.5].
By the exact sequence 
$$\ex{X}{-2E-K_X}{X}{-E-K_X}{E}{-K_E} $$
and the KKV vanishing theorem,
$|M|$ is also free.
Let $G$ be a general member of $|M|$, $l:=E|_G$ and $G' :=f(G)$. 
Then $Q$ is a minimally elliptic singularity of $G'$ 
by the formula $K_G = {f|_G}^*K_{G'} - l$ and [La, Theorem 3.4].
On the other hand, the embedded dimension of $Q$ is at most $4$
since $Q$ is a cDV singularity on $Y$.
Hence we have $-(l^2)_G \leq 4$ by [La, Theorem 3.13]
whence $(-K_E)^2 = -2(l^2)_G = 2,4,6,8$.
These correspond to type $E_7\sim E_{10}$ respectively.
\comment
Let $\t{X} :=\pmb{Spec}  (\Cal O_X \oplus \Cal O_X(K_X))$, where
we define a ring structure of $\Cal O_X \oplus \Cal O_X(K_X)$ by $G$. 
Let $\t{E}$ the pull back of $E$ on $\t{X}$. Note that $\t{E}$ is irreducible.
Then there is a natural crepant contraction
of $\t{E}$ from $\t{X}$. Note that $\t{X}$ is smooth.
Assume that $E$ is non normal.
Then $\t{E}$ is also non normal and hence by [G, Theorem 5.2],
$\t{E}^3 = 2E^3$ must be $7$ but this contradicts above.

Hence $E$ is normal. Assume that $E$ is not rational.
Then $E$ is a cone over an elliptic curve (cf. [Ba]).
Let $\mu:\o{X} \to X$ be the blow up at the vertex $v$ of $E$,
$\o{E}$ the strict transform of $E$ and $F$ the exceptional divisor.
Then $m:=\o{E} \cap F$ is an elliptic curve in $F\simeq \Bbb P^2$
whence $(m^2)_{\o{E}}=-6$. 
Then $(-K_E)^2 = (-K_{\o{E}})^2 - (m^2)_{\o{E}} = 6$.
So $E^3 = \frac 32$.

Assume that $E$ is a rational surface.
If $(-K_E)^2 = 8$, then $E\simeq \Bbb P^1 \times \Bbb P^1 \ \text{or} \
\Bbb F_{2,0}$. The former case forces $X$ to be smooth, a contradiction.
The second case also derive a contradiction by investigating the blow up
of the vertex of $E$.
Hence we obtain the required bound of $E^3$.
\endcomment
\qed
\enddemo

\proclaim {Proposition 2.4 (Contraction of type $C$)}
Let $X$ be an analytic $3$-fold with only index $2$ terminal singularities 
and $f: X\to (Y, Q)$ an extremal contraction of type $C$
to a germ of surface. Let $l$ be the fiber over $Q$.
Then $Q$ is a smooth point or an ordinary double point.
Furthermore the following description holds:
 
\roster
\item 
if $l$ contains no index $2$ point,
$Q$ is a smooth point and $f$ is a usual conic bundle; 
\item 
if $l$ contains an index $2$ point and $Q$ is a smooth
point, $l$ contains only one index $2$ point and 
every component $l'$ of $l$ passes through it.
Furthermore $-K_X.l'=\frac 12$;
\item if $l$ contains an index $2$ point and $Q$ is an ordinary
double point, $f$ is analytically isomorphic to one of the following:

\item "(3-1)"
Let $\Bbb P^1 \times (\Bbb C^2,o) \to (\Bbb C^2,o)$ be the natural projection.
Define the action of the group $\Bbb Z_2$ on 
${\Bbb P^1}_{x_0,x_1} \times {\Bbb C^2}_{u,v}$:
$$(x_0,x_1;u,v) \to (x_0,-x_1;-u, -v).$$
Let $X=\Bbb P^1 \times \Bbb C^2 / \Bbb Z_2$ and $(Y,Q)=(\Bbb C^2/\Bbb Z_2,o)$.

In particular $X$ has two $\frac 12(1,1,1)$-singularities on $l$
and $l_{\text{red}}\simeq \Bbb P^1$ and $-K_X.l_{\text{red}}=1$.

\item "(3-2)"
Let $X'$  be a hypersurface in 
${\Bbb P^2}_{x_0,x_1,x_2} \times {\Bbb C^2}_{u,v}$ defined by the equation
${x_0}^2 + {x_1}^2 + {x_2}^2 \phi (u,v) =0$,
where $\phi (u,v)$ has no multiple factors and contains only monomials of
even degree.
Let $f' : X' \to \Bbb C^2$ be the natural projection.
Define the action of the group $\Bbb Z_2$ on $X'$ as follows:
$$(x_0,x_1,x_2;u,v) \to (-x_0,x_1,x_2;-u, -v).$$
Let $X:=X'/ \Bbb Z_2$ and $(Y,Q)=(\Bbb C^2/\Bbb Z_2,o)$.

In particular
$P$ is the unique index $2$ point and aw$(X,P)=2$.
If $\text{mult}_{(0,0)} (\phi)=2$, then $(X,P)$ is a $cA /2$ point
or if $\text{mult}_{(0,0)} (\phi)\geq 4$, then $(X,P)$ is a $cAx /2$ point.
\endroster
\endproclaim

\demo{Proof}
See [Mo1, Theorem 3.5] for (1) and
[Pr, Theorems 3.1, 3.15 and Examples 2.1 and 2.3] 
for (2) and (3).
\qed
\enddemo

\proclaim{Proposition 2.5 (Contraction of type $D$)}
Let $X$ be an analytic $3$-fold \newline
with only index $2$ terminal singularities
and $f: X\to (C,Q)$ be an extremal contraction of type $D$
to a germ of a curve.
Let $F$ be the fiber over $Q$.
Then $Q$ is a smooth point and the following description holds:
\roster
\item 
if $F$ contains no index $2$ point,
then all fibers are irreducible and reduced and
(possibly non normal) Gorenstein del Pezzo surfaces.
Furthermore
if $(-K_F)^2 =9$, we can write $-K_X \sim 3A$ for some relatively ample
divisor $A$ and $X=\Bbb P(f_*\Cal O_X (A))$ which is a $\Bbb P^2$-bundle;

if $(-K_F)^2 =8$, we can write $-K_X \sim 2A$ for some relatively ample
divisor $A$ and $X$ is embedded in $\Bbb P^3$-bundle
$\Bbb P(f_*\Cal O_X (A))$ as a quadric bundle
(the last means all fibers are quadrics in $\Bbb P^3$);

the case $(-K_F)^2 = 7$ does not occurs.
\item 
if $F$ contains an index $2$ point,
then $F$ is irreducible and reduced  
or $F=2F_{\text{red}}$ and $F_{\text{red}}$ is irreducible.
$F_{\text{red}}$ is a del Pezzo surface of Gorenstein index $2$.
\endroster
\endproclaim  

\demo{Proof}
See [Mo1, Theorem 3.5] for (1). (2) follows 
from the existence of a section [Co].
\qed
\enddemo

\head 3. Takeuchi's theory 
\endhead

\definition{Definition 3.0}
Let $X$ be a $\Bbb Q$-Fano variety.
We say that
a birational morphism $f:Y\to X$ is a weak $\Bbb Q$-Fano blow up
if the following hold:
\roster
\item
$Y$ is a weak $\Bbb Q$-Fano variety;
\item $f$ is an extremal divisorial contraction such that
$f$-exceptional divisor is a prime $\Bbb Q$-Cartier divisor. 
\endroster
\enddefinition

In this section, 
we consider a $\Bbb Q$-factorial 
$\Bbb Q$-Fano $3$-fold $X$ with the following properties:

\definition{Assumption 3.1}
\roster 
\item the Picard number $\rho(X)$ is $1$;
\item 
there is a weak $\Bbb Q$-Fano blow up $f:Y\to X$.
Let $E$ be the $f$-exceptional divisor. 
\endroster
\enddefinition

We fix $f$ as in Assumption 3.1 and
set $\alpha := (-K_Y)^2 E$.

Consider the extremal ray $R$ of $Y$ other than the ray associated to $f$.
If $R$ is a ray associated to a non small contraction, 
denote by $f' : Y\to X'$ the contraction associated to $R$.
If $R$ is a flopping ray, then after the flop 
$Y_0:=Y \overset g_0 \to \dashrightarrow Y_1$,
another extremal ray of $Y_1$ is $K_{Y_1}$-negative because
$K_{Y_1}$ is not nef and $\rho(Y_1)=2$.
Hence we can start the minimal model program from $Y$ or $Y_1$ and 
because the canonical bundle can not become nef while the program,
we obtain the following diagram:
$$\matrix
\ & Y & \dashrightarrow & {Y'} \\
\ & {f\swarrow} & \ &  {\searrow f'} \\
X & \ & \ & \ & X' & ,
\endmatrix $$
where $Y \dashrightarrow Y'$ is an isomorphism or a composition
of possibly one flop and flips and
$f'$ is the first non small contraction.
Let $E$ be the exceptional divisor of $f$.
We do the similar calculations as Kiyohiko Takeuchi did in [T1] in the below.
The following lemma is basic for our computations:

\proclaim{Lemma 3.2}
Assume that $Y \dashrightarrow Y'$ is not an isomorphism.
Let $$Y_0:=Y \overset g_0 \to \dashrightarrow Y_1 \dots
\overset g_{k-1} \to \dashrightarrow Y_k:=Y'$$ be the decomposition
of $Y \dashrightarrow Y'$ into flops and flips.
Let $l_i$ be an irreducible component of the flipping (or flopping) curve
for $g_i$ and $E_i$ the strict transform of $E$ on $Y_i$.
Then
\roster
\item 
there is at most one flop in the above decomposition
and if there is, the flop is $g_0$;
\item 
$E_i.l_i >0$;
\item 
if $Y \overset g_0  \to \dashrightarrow Y_1$ is a flop, then
$(-K_{Y_1})^3 = (-K_Y)^3$, $(-K_{Y_1})^2 E_1=(-K_Y)^2 E$,
$(-K_{Y_1}) {E_1}^2=(-K_Y) E^2$ and 
$e:= E^3 - {E_1}^3 \in \frac {\Bbb N}s$,
where $s$ is the minimum positive integer such that $sE$ is a Cartier divisor;
\item 
if $Y_i \overset g_i \to \dashrightarrow Y_{i+1}$ is a flip, let
$d_i := (-K_{Y_i})^3 -(-K_{Y_{i+1}})^3$. Then $d_i >0$
and 
$(-K_{Y_{i+1}})^2 E_{i+1}=(-K_{Y_i})^2 E_i-a_i d_i$,
$(-K_{Y_{i+1}}) {E_{i+1}}^2=
(-K_{Y_i}) {E_i}^2 - {a_i}^2 d_i$ 
and ${E_{i+1}}^3={E_i}^3 - {a_i}^3 d_i$,
where $a_i:= \frac{{E_i}.l_i}{(-K_{Y_i}).l_i}$
(note that this number $a_i$ is well defined since flipping curves are
numerically proportional); 
\item 
if $Y_i \overset g_i \to \dashrightarrow Y_{i+1}$ and
$Y_{i+1} \overset g_{i+1} \to \dashrightarrow Y_{i+2}$
are flips, then $a_{i+1} < a_i$; 

\item
$\text{Pic} Y'$ is torsion free.
\endroster
\endproclaim

\demo{Proof}
\roster
\item 
This is clear from the above consideration;

\item
we prove this by induction for $i$.
For $i=0$, assume that $E_0.l_0 \leq 0$.
Then $E_0$ is non positive for two extremal rays of $Y$ and hence
$E_0$ is non positive for all effective curves on $Y$ since $\rho(Y)=2$.
But this contradicts the effectivity of $E$.
Assume that $E_i.l_i >0$ is proved. Then $E_{i+1}.{l_i}^+<0$,
where ${l_i}^+$ is the flipped curve corresponding to $l_i$.
Hence we can prove $E_{i+1}.l_{i+1} >0$ by the same way as 
proving $E_0.l_0 >0$;
\item
let 
$$\matrix
\ & \ & Z & \ \\
\ & {p\swarrow} & \ &  {\searrow q} \\
Y & \ & \ & \ & Y_1 & ,
\endmatrix $$
the common resolution of $Y$ and $Y_1$.
Then by the negativity lemma ([FA, Lemma 2.19]), 
we can easily see that
$p^* K_Y = q^* K_{Y_1}$ (for example, see [Kol1, Proof of Lemma 4.3]
or below argument).
By this, former $3$ equalities follows.
We prove that $e\in \frac {\Bbb N}s$.
Since $sE_1$ is Cartier by Proposition 1.5, 
we have $e\in \frac{\Bbb Z}s$.  
Let 
$$p^{-1} E = p^* E-R = q^* E_1  -R' ,$$
where $R$ and $R'$ are effective divisors 
which are exceptional for $p$ and $q$.
Rewrite this as
$$-p^* E = -q^* E_1 + R'-R.$$
Then since $-q^* E_1$ is $p$-nef by (3),
we see that $R'-R >0$ and $p_* (R'-R) \not = 0$
by $E.l_0 >0$ and the negativity lemma.
Hence we can write $p^* E= q^* E_1 -F$, where $F:=R' -R$ is an effective
divisor.
Consider the identity $(p^* E)(q^* E_1)^2 = (q^* E_1 -F)(q^* E_1)^2$.
Its right side is equal to ${E_1}^3$.
Its left side is equal to $(p^* E)(p^* E +F)^2 
= E^3 + E. p_* (F^2)$.
By $p_* F \not =0$, we know that 
$-p_* (F^2)$ is a non zero effective $1$-cycle.
Hence $E. p_* (F^2)<0$ and we are done;
\comment
We will show that $p_* (F^2) \equiv \sum \alpha _i m_i$,
where $m_i$ is a flopping curve and $\alpha _i$ is some non negative integer.
Let $G$ be a $p$-exceptional curve whose image on $Y$ is a curve.
If $G$ appears with non integer coefficient in $F$,
then $q(G)$ is a non Gorenstein point on $Y_1$
and the discrepancy of $G$ for $K_{Y_1}$ is not an integer.
On the other hand 
the discrepancy of $G$ for $K_Y$ is an integer, a contradiction to
$p^* K_Y = q^* K_{Y_1}$. 
Hence $p_* (F^2)$ is a $1$-cycle with integer coefficient. 
By restricting $q$ to a general hyperplane section at a general point of $m_i$,
we can see that $-q_*(F^2) \equiv \sum m_i$, where $m_i$'s are flopped curves. 
Hence we have $e:=-\sum (E_1.m_i) = E^3 - {E_1}^3$ and 
$e$ is an positive integer by (3).
Furthermore if $E.l=1$ for any flopping curve $l$, then clearly
$e$ is the number of flopping curves.
\endcomment 

\item
the proof is very similar to one of (4).
Let $$\matrix
\ & \ & Z & \ \\
\ & {p\swarrow} & \ &  {\searrow q} \\
Y_i & \ & \ & \ & Y_{i+1} & ,
\endmatrix $$
the common resolution of $Y_i$ and $Y_{i+1}$.
By the definition of $a_i$,
$$H_i := a_i (-K_{Y_i}) - E_i \tag a$$ is numerically trivial 
for the flipping curves.
Let ${H_i}^+$ be the strict transform of $H_i$.
By the negativity lemma, we can easily see that
$p^* H_i = q^* {H_i}^+$
and 
$p^*(-K_{Y_i})
=q^*(-K_{Y_{i+1}}) -G$, where $G$ is 
an effective divisor
which is exceptional for $p$ and $q$.
$d_i >0$ can be proved similarly to the proof of positivity of $e$.
Consider the following identities:
$$\multline
(-K_{Y_i})^2 H_i
= (p^*(-K_{Y_i}))^2 p^*H_i= \\
(q^*(-K_{Y_{i+1}})-G)^2 q^* {H_i}^+
=(-K_{Y_{i+1}})^2 {H_i}^+  
\endmultline
\tag b$$
and similarly
$$(-K_{Y_i}){H_i}^2
=(-K_{Y_{i+1}}){{H_i}^+}^2 \tag c $$
and
$${H_i}^3 = {H^+_i}^3 . \tag d$$

By (a)$\sim$(d) and the definition of $d_i$, we obtain the assertion;
\item
let $l_{i}^+$ be a flipped curve on $Y_{i+1}$.
By $(a_i(-K_{Y_{i+1}}) -E_{i+1}).l_{i}^+ =0$ and
$(a_i(-K_{Y_{i+1}}) -E_{i+1}).m >0$ for a general curve $m$ on $Y_{i+1}$,
we have $(a_i(-K_{Y_{i+1}}) -E_{i+1}).l_{i+1}>0$.
On the other hand we have
$(a_{i+1}(-K_{Y_{i+1}}) -E_{i+1}).l_{i+1}=0$.
Hence we are done;
\item
it is easy to see by Riemann-Roch theorem that
$\text{Pic} Y$ is torsion free since $Y$ is a weak $\Bbb Q$-Fano $3$-fold.
Since $Y\dashrightarrow Y'$ is a composition of a flop or flips
and linear equivalence is preserved under a flop and a flip,
$\text{Pic} Y'$ is also torsion free.
\endroster
\qed
\enddemo

We will define $e$, all $a_i$'s and $n_i$'s to be $0$ if $Y=Y'$.
If $Y\dashrightarrow Y'$ is not an isomorphism,
we will define $e$ to be $0$ if $Y\dashrightarrow Y_1$ is not a flop 
and $a_i$ and $n_i$ to be $0$ 
if $Y_i \dashrightarrow Y_{i+1}$ is not a flip.  
From now on, we divide $f'$ into cases. For this, we have

\proclaim{Claim}
If $f'$ is a crepant contraction, then
$Y=Y'$ and $\dim X' = 3$. 
\endproclaim

\demo{Proof}
The fact that $Y=Y'$ is clear by consideration above Lemma 3.2.
Since $-K_Y$ is a supporting divisor of $f'$ and 
$-K_Y$ is nef and big, $\dim X' =3$.
\qed
\enddemo

Hence we have the following cases:

\definition{Case 1}
$f'$ is an extremal contraction of type $E_1$.
\enddefinition

\definition{Case 2}
$f'$ is an extremal contraction of type $E_2 \sim E_{11}$.
\enddefinition

\definition{Case 3}
$f'$ is an extremal contraction of type $C$.
\enddefinition

\definition{Case 4}
$f'$ is an extremal contraction of type $D$.
\enddefinition

\definition{Case 5}
$f'$ is a crepant divisorial contraction.
\enddefinition

\proclaim{Claim 3.3}
$\t{E}$ and $-K_{Y'}$ are numerically independent.
\endproclaim

\demo{Proof}
For $Y'$, the numerical equivalence is equal to the
$\Bbb Q$-linear equivalence by Lemma 3.2 (6).
So the assertion follows since no multiple of $\t{E}$ moves and
$-K_{Y'}$ is big.
\qed
\enddemo

In Case 1, 2 or 5,
let $E'$ be the exceptional divisor of $f'$,
$\t{E}$ the strict transform of $E$ on $Y'$ and $\t{E'}$
the strict transform of $E'$ on $Y$.
By Claim 3.3 and $\rho (Y') = 2$, 
we can write $$E'\equiv z(-K_{Y'})-u\t{E}.\tag 3-0-1$$ 

In Case 3 or 4, let $L$ be the pull back of the ample generator of 
$\text{Pic} X'$ and $\t{L}$ the strict transform of $L$ on $Y$.
By Claim 3.3 and $\rho (Y') =2$, 
we can write $$L \equiv z(-K_{Y'})-u\t{E}.\tag 3-0-2$$

\definition{Assumption 3.4}
In the below we further assume that
$P:=f(E)$ is a point of index $r$ and
$-K_Y = f^*(-K_X) -\frac 1r E$ and write
$-K_X\equiv q S$, where $S$ is the positive generator of $Z^1(X)/\equiv$
and $q$ is a positive integer.
\enddefinition

Then we have

\proclaim{Claim 3.5}
$z \in \Bbb N / q$ and $u$ is a positive rational number such that 
$z + ru \in \Bbb N$.
\endproclaim

\demo{Proof}
On $X$, $f(\t{E'})\equiv zq S$ in Case 1, 2 or 5
(resp. $f(\t{L})\equiv zqS$ in Case 3 or 4).

So by Assumption 3.4, $z \in \frac{\Bbb N}q$. 

If $u\leq 0$, sufficient multiple of $E'$ must be move in Case 1, 2 or 5
(resp. $L$ must be big in Case 2 or 3), a contradiction.

Let $\t{E'}$ be the strict transform of $E'$ on $Y$.
By (3-0-1) and $-K_Y = f^*(-K_X) - \frac 1r E$,
we have $\t{E'} \equiv zf^*(-K_X) -(\frac zr + u)E$ in Case 1, 2 or 5
(resp. by (3-0-2) and $-K_Y = f^*(-K_X) - \frac 1r E$,
we have $\t{L} \equiv zf^*(-K_X) -(\frac zr + u)E$ in Case 3 or 4).
Hence $\frac zr + u \in \frac {\Bbb N}r$. 
\qed
\enddemo

\definition {Case 1}
Let $C:=f'(E')$. 
\enddefinition

\proclaim{Claim 3.6}
$z+1=uk$ for some $k\in \Bbb N$.
\endproclaim

\demo{Proof}
By (3-0-1) and $-K_{Y'} = {f'}^*(-K_{X'}) - E'$,
we have $(z+1)E' \equiv z{f'}^*(-K_{X'}) -u\t{E}$.
Since $f'(\t{E})$ is a Cartier divisor along $C$ 
outside a finite set of points,
$\frac {z+1} u$ is an integer.
\qed
\enddemo

We have the following:

Recall that $\alpha := (-K_Y)^2 E$.

$$\multline
(-K_{Y'}+E')^2 (-K_{Y'})= \\
(z+1)^2(-K_{Y'})^3 -2u(z+1)(-K_{Y'})^2\t{E}+u^2(-K_{Y'})\t{E}^2 =(-K_{X'})^3.
\endmultline
\tag {3-1-1}$$

$$\multline
(-K_{Y'}+E')^2 E'=\\
z(z+1)^2(-K_{Y'})^3 -u(z+1)(3z+1)(-K_{Y'})^2\t{E}+u^2(3z+2)(-K_{Y'})\t{E}^2-
u^3 \t{E}^3 = 0. 
\endmultline
\tag {3-1-2}$$

$$\multline
(-K_{Y'}+E')E'(-K_{Y'})=\\
(z+1)z(-K_{Y'})^3 -u(2z+1)(-K_{Y'})^2\t{E}+u^2(-K_{Y'})\t{E}^2 =(-K_{X'}.C).
\endmultline
\tag {3-1-3}$$

$$\multline
(-K_{Y'}-E')^2 E'= 4\{ (-K_{X'})^3 -(-K_{Y'})^3 -2(-K_{X'} .C) \} \\
z(z-1)^2(-K_{Y'})^3 -u(z-1)(3z-1)(-K_{Y'})^2\t{E}+u^2(3z-2)(-K_{Y'})\t{E}^2-
u^3 \t{E}^3 = \\
(-K_{E'})^2 \leq 8 (1-g(\o{C})),
\endmultline \tag {3-1-4}$$
where $\o{C}$ is the normalization of $C$.
(the last inequality of (3-1-4) can be proved similarly 
to the proof of Proposition 2.2.)

Hence by Lemma 3.2, we obtain the following:

$$
\{k^2 (-K_Y)^3 - (2k+r) \alpha -
\sum d_i (a_i -k)^2 \}u^2 =(-K_{X'})^3.
\tag {3-1-1'}$$

$$\multline
(uk-1) k^2 (-K_Y)^3 - \{ur^2 +(3uk-1)r +k(3uk-2)\} \alpha +
 \sum d_i \{u (a_i -k)^3 + (a_i -k)^2 \}+eu= 0. 
\endmultline
\tag {3-1-2'}$$

$$\{k(uk-1) (-K_Y)^3 -(2uk -1+ur)\alpha -
\sum d_i (a_i -k)(a_i u -ku +1) \}u =(-K_{X'}.C).
\tag {3-1-3'}$$
The positivity of the left hand side gives some information.
By (3-1-1') and (3-1-2'), we have the following:
$$e+\sum d_i a_i (a_i -k)^2 = (r+k)^2 \alpha - \frac{uk-1}{u^3} (-K_{X'})^3.
\tag {3-1-5'}$$

The following claim is useful for solving the equations:

\proclaim{Claim 3.7}
If $Y_i \dashrightarrow Y_{i+1}$ is a flip, then
$k <a_i$.
\endproclaim

\demo{Proof}
Note that ${f'}^* (-K_{X'}) \equiv u\{k(-K_{Y'}) - \t{E}\}$.
Hence $k(-K_{Y_i}) - E_i$ is $\Bbb Q$-effective for any $i$.
If $Y_i \dashrightarrow Y_{i+1}$ is a flip and
$k \geq a_i$ for some $i$, then 
$(k(-K_{Y_i}) - E_i).l_i \geq 0$ and hence
$(k(-K_{Y_{i+1}}) - E_{i+1}).{l^+_i} \leq 0$.
By the $\Bbb Q$-effectivity of $k(-K_{Y_{i+1}}) - E_{i+1}$ and 
$\rho (Y_{i+1}) =2$, $k(-K_{Y_{i+1}}) - E_{i+1}$ is positive for 
another extremal ray of $Y_{i+1}$. So 
$k(-K_{Y'}) - \t{E}$ is positive for a fiber of $f'$. But this is absurd.
\qed
\enddemo

\definition{Case 2}
We note that $Y\not = Y'$ since otherwise
$E\not = E'$ and the nonempty intersection curve $E\cap E'$ is contracted
by two extremal contractions $f$ and $f'$, a contradiction.
\enddefinition
Let $\frac {d}{r'}$ be the discrepancy of $E'$ for $K_{X'}$,
where $r'$ is the index of $P':=f'(E')$.
By the similar way to the proof of Claim 3.6, we have the following claim:
\proclaim{Claim 3.8}
$zd+r'=uk$ for some $k\in \Bbb N$.
\endproclaim

\demo{Proof}
By (3-0-1) and $-K_{Y'} = {f'}^*(-K_{X'}) - \frac {d}{r'} E'$,
we have $(zd+r')E' \equiv r'z{f'}^*(-K_{X'}) -ur' \t{E}$.
Since $r'f'(\t{E})$ is Cartier divisor at $P'$,
$\frac {zd+r'} u$ is an integer.
\qed
\enddemo

We have the following:

$$
z^3(-K_{Y'})^3 -3z^2 u(-K_{Y'})^2\t{E}+3zu^2(-K_{Y'})\t{E}^2-
u^3 \t{E}^3 = (E')^3. 
\tag {3-2-1}$$

$$
z^2(-K_{Y'})^3 -2zu(-K_{Y'})^2\t{E}+u^2(-K_{Y'})\t{E}^2
=(-K_{Y'})(E')^2.
\tag {3-2-2}$$

$$
z(-K_{Y'})^3 -u(-K_{Y'})^2\t{E}=(-K_{Y'})^2 E'. 
\tag {3-2-3}$$

Hence by Lemma 3.2, we obtain the following:

$$
z^3 (-K_Y)^3 -u \alpha (u^2 r^2 +3zur +3z^2) +\sum d_i (u a_i -z)^3 + u^3 e 
=(E')^3. 
\tag {3-2-1'}$$

$$z^2 (-K_Y)^3 -u\alpha (2z+ur) -\sum d_i (u a_i -z)^2 
= (-K_{Y'})(E')^2.
\tag {3-2-2'}$$

$$z(-K_Y)^3 -u \alpha +\sum d_i (u a_i -z) = (-K_{Y'})^2(E'). \tag {3-2-3'}$$

By (3-2-1') and (3-2-2'), we have the following:

$$
\sum d_i a_i (a_i u -z)^2 +u^2 e = \alpha (z +ur)^2 + \frac k{r'} (E')^3.
\tag {3-2-4'}$$

By (3-2-2') and (3-2-3'), we have the following:

$$\alpha (z +ur)+
\sum d_i a_i (a_i u -z)= \frac {dk}{{r'}^2} (E')^3.
\tag {3-2-5'}$$

Similarly to Claim 3.7, we have the following:

\proclaim{Claim 3.9}
If $Y_i \dashrightarrow Y_{i+1}$ is a flip, then
$k <d a_i$.
\endproclaim

\demo{Proof}
Note that ${f'}^* (-K_{X'}) \equiv \frac u{r'} \{k(-K_{Y'}) - d \t{E}\}$.
Hence $k(-K_{Y_i}) - d E_i$ is $\Bbb Q$-effective for any $i$.
The rest is similar to the proof of Claim 3.7.
\qed
\enddemo

\definition{Case 3}
By [Pr, Lemma 1.10], $X'$ has only cyclic quotient singularities.
By the general theory of the conic bundle, 
$-4K_{X'} \equiv {f'}_*(-K_{Y'}^2)+\Delta$, where $\Delta$ is the discriminant
divisor of $f'$. Hence $-K_{X'}.A >0$ for any ample divisor $A$ on $X'$
since $-K_{Y'}$ is big.
Hence $X'$ is a log del Pezzo surface with $\rho(X')=1$.
\enddefinition

We have the following:

$$L^3 =
z^3(-K_{Y'})^3 -3z^2 u(-K_{Y'})^2\t{E}+3zu^2(-K_{Y'})\t{E}^2-
u^3 \t{E}^3 = 0. \tag {3-3-1}$$

$$
z^2(-K_{Y'})^3 -2zu(-K_{Y'})^2\t{E}+u^2(-K_{Y'})\t{E}^2 =
(-K_{Y'})L^2.
\tag {3-3-2}$$

$$z(-K_{Y'})^3 -u(-K_{Y'})^2\t{E}=(-K_{Y'})^2 L. 
\tag {3-3-3}$$

We set $u = mz$ and $l = f'_* L$.
By Lemma 3.2, we obtain the following:

$$
(-K_Y)^3 -m\alpha (m^2 r^2 +3mr +3) +\sum d_i (m a_i -1)^3 + m^3 e =0. 
\tag {3-3-1'}$$

$$z^2 \{(-K_Y)^3 -m\alpha (2+mr) -\sum d_i a_i (m a_i -1)^2 \} 
= 2 l^2.
\tag {3-3-2'}$$

$$z\{(-K_Y)^3 -m\alpha +\sum d_i (m a_i -1)\} = (-K_{Y'})^2 L. \tag {3-3-3'}$$
If $l$ is free, then 
$(-K_{Y'})^2 L = 8(1-g(l))-\Delta .l +4l^2$.

By (3-3-1') and (3-3-2'), we have the following:
$$z^2 \{\sum d_i m a_i (m a_i -1)^2 +m^3 e\} =z^2 m\alpha (mr +1)^2 -2 l^2.
\tag {3-3-4'}$$

\definition{Case 4}
By the edge sequence of the Leray spectral sequence
$0 \to H^1(X', \Cal O_{X'}) \to H^1(Y', \Cal O_{Y'})$
(exact) and $H^1(Y', \Cal O_{Y'})=0$, we have $H^1(X', \Cal O_{X'})=0$, i.e.,
$X' \simeq \Bbb P^1$.
Hence $L= {f'}^*\Cal O_{\Bbb P^1}(1)$
\enddefinition

We calculate the following:

$$(-K_{Y'})L^2=
z^2(-K_{Y'})^3 -2zu(-K_{Y'})^2\t{E}+u^2(-K_{Y'})\t{E}^2 =0.
\tag {3-4-1}$$

$$\t{E} L^2=
z^2(-K_{Y'})^2\t{E} -2zu(-K_{Y'})\t{E}^2+u^2\t{E}^3 =0.
\tag {3-4-2}$$

$$(-K_{Y'})^2 L=
z(-K_{Y'})^3 -u(-K_{Y'})^2\t{E}= \ \text{deg} F, 
\tag {3-4-3}$$
where $F$ is a general fiber of $f'$ and $\text{deg}  F := (-K_F)^2$.

We set $u=mz$. We obtain the following:

$$(-K_X)^3 = \frac {\alpha}r + m\alpha (2+mr) + \sum d_i (m a_i -1)^2.
\tag {3-4-1'}$$

$$(mr+1)^2 \alpha = \sum d_i a_i (m a_i -1)^2 + m^2 e. \tag {3-4-2'}$$

$$z\{m\alpha (1+mr) 
+\sum d_i ma_i(m a_i -1)\} = \ \text{deg}  F. \tag {3-4-3'}$$

The following claim is useful for solving the equations:
\proclaim{Claim 3.10}
In Case 3 or 4,
if $Y_i \dashrightarrow Y_{i+1}$ is a flip, then
$m a_i >1$.
\endproclaim

\demo{Proof}
Note that $L \equiv z(-K_{Y'} -m \t{E})$.
The proof is similar to the one of Claim 3.7.
\qed
\enddemo

\definition {Case 5}
Since $-K_Y.l=0$ and $E'.l= -2$ for a general fiber $l$ of $E'$, 
we have $u(E.l)=2$. 
By $(-K_Y)^2 E' =0$,
we have $z(-K_Y)^3 =u\alpha$.
\enddefinition

By an additional geometric assumption that $|-K_Y -E| \not = \phi$,
the relation of $u$ and $z$ is restricted as follows:

\proclaim{Claim 3.11}
If $|-K_Y -E| \not = \phi$, then $z \leq u$.
Furthermore in Case 3, $m=1 \ \text{or} \ 2$
and in Case 4, $m = 1$ or $m=2$ and $F\simeq \Bbb P^1 \times \Bbb P^1$
or $m=\frac 32 \ \text{or} \ 3$ and $F\simeq \Bbb P^2$.
\endproclaim

\demo{Proof}
By (3-0-1), we have 
$$E' \equiv (z-u)(-K_{Y'})+u(-K_{Y'}-\t{E}) \tag a$$ in Case 1, Case 2
or case 5
(resp. by (3-0-2),
$$L \equiv (z-u)(-K_{Y'})+u(-K_{Y'}-\t{E}) \tag b$$ in Case 3 or Case 4).
By the assumption,
$|-K_{Y'} -\t{E}|\not = \phi$.
Hence if $z>u$, sufficient multiple of $E'$ must move by (a)
(resp. if $z>u$, $\kappa (L)$ must be $3$ by (b)), a contradiction.
So $z\leq u$.

In Case 3, for a general fiber $C$, we have
$\t{E}.C = \frac {2z}u \in \Bbb N$. So $\frac {2z}u = 1 \ \text{or} \ 2$
since $z\leq u$. 
In Case 4, let $C$ be a $(-1)$-curve in $F$ if $F \not \simeq 
\Bbb P^1 \times \Bbb P^1, \Bbb P^2$ or a ruling if 
$F \simeq \Bbb P^1 \times \Bbb P^1$ or a line if
$F \simeq \Bbb P^2$.
By calculating $\t{E}.C$, we obtain the assertion similarly to Case 3.
\qed
\enddemo

\head 4. Existence of a weak $\Bbb Q$-Fano blow up for a 
$\Bbb Q$-Fano $3$-fold with $I(X)=2$ \endhead

In this section, we give a sufficient condition for
the existence of a weak $\Bbb Q$-Fano blow up for a 
$\Bbb Q$-Fano $3$-fold with $I(X)=2$.  

\proclaim{Theorem 4.0}
Let $X$ be a weak $\Bbb Q$-Fano $3$-fold with only log terminal singularities.
Assume the following:
\roster 
\item $I(X) \leq 2$;
\item there are only a finite number of non Gorenstein points on $X$;
\item $(-K_X)^3 \geq 1$ and $h^0 (-K_X)\geq 1$.
\endroster
Then $|-2K_X|$ is free.
\endproclaim

\demo{Proof}
By replacing $X$ by its anti-canonical model, we can
assume that $X$ is a $\Bbb Q$-Fano $3$-fold 
with only log terminal singularities.
By [Am, Theorem 1.2], $S$ has only log terminal singularities.
By the exact sequence 
$0\to \Cal O_X \to \Cal O_X(-2K_X) \to \Cal O_S(-2K_X |_S) \to 0$
and $h^1(\Cal O_X)=0$,
we have $|-2K_X|_S|=|-2K_X||_S$ and 
$\text{Bs} |-2K_X| = \text{Bs} |-2K_X |_S|$.
Note that $-K_X |_S = K_S$.
Hence it suffices to prove that $|K_S + K_S|$ is free.
Assume that $|2K_S|$ is not free.
Let $y$ be a base point of $|K_S+K_S|$.
Assume that $y$ is worse than canonical.
By [KT, Theorem 9], $y$ is a cyclic quotient singularity
of index $2$. So Kawachi's invariant $\delta '$ defined in [KT] is $\frac 12$
at $y$.
On the other hand 
by the assumption that $(-K_X)^3 \geq 1$, ${K_S}^2 \geq 2$ holds.
So ${K_S}^2 > \delta _y$ holds ($\delta _y$ is defined in [KaMa]).
But by (1), 
we have $K_S .C = -K_X.C \geq \frac 12$ for any curve $C$
whence by [ibid.], $y$ cannot be a base point of $|2K_S|$, 
a contradiction.
So we may assume that $S$ does not contain a non Gorenstein point of $X$
by (2) and has only canonical singularities.
Let $\mu: \t{S} \to S$ be the minimal resolution.
Since $h^0 (K_{\t{S}})=h^0 (K_S)=h^0 (-K_X)\geq 1$, $|2K_{\t{S}}|$ is free
by [Fr] and hence $|2K_S|$ is free, a contradiction again. 

Hence $|K_S + K_S|$ is free and also $|-2K_X|$ is free.
\qed
\enddemo

\comment
\proclaim{Theorem 4.1}
Let $X$ be a weak $\Bbb Q$-Fano $3$-fold with only log terminal singularities.
Assume the following:
\roster 
\item $I(X) \leq 2$;
\item there are only a finite number of non Gorenstein points on $X$;
\item $(-K_X)^3 > 2$.
\endroster
Then $|-2K_X|$ is free.
\endproclaim

\demo{Proof}
By replacing $X$ by its anti-canonical model, we can
assume that $X$ is a $\Bbb Q$-Fano $3$-fold 
with only log terminal singularities.
By [Am, Theorem 1.2], $S$ has only log terminal singularities.
By the exact sequence 
$0\to \Cal O_X \to \Cal O_X(-2K_X) \to \Cal O_S(-2K_X |_S) \to 0$
and $h^1(\Cal O_X)=0$,
we have $|-2K_X|_S|=|-2K_X||_S$ and 
$\text{Bs} |-2K_X| = \text{Bs} |-2K_X |_S|$.
Note that $-K_X |_S = K_S$.
Hence it suffices to prove that $|K_S + K_S|$ is free.
Assume that $|2K_S|$ is not free.
Let $y$ be a base point of $|K_S+K_S|$.
By the assumption that $(-K_X)^3 >2$, ${K_S}^2 > 4$ holds.
So ${K_S}^2 > \delta _y$ holds ($\delta _y$ is defined in [KaMa].)

Assume that $y$ is worse than canonical.
By [KT, Theorem 9], $y$ is a cyclic quotient singularity
of index $2$. So Kawachi's invariant $\delta '$ defined in [KT] is $\frac 12$
at $y$. By (1), 
we have $K_S .C = -K_X.C \geq \frac 12$ for any curve $C$.
Hence by [ibid.], $y$ cannot be a base point of $|2K_S|$, 
a contradiction.
So we may assume that $S$ does not contain a non Gorenstein point of $X$
by (2).
In particular for all curve $C\subset S$, we have $K_S . C \geq 1$. 
Hence 
if $y$ is a canonical singularity, then 
again by [ibid.], $y$ cannot be a base point of $|2K_S|$, 
a contradiction.

Hence we may assume that $S$ has only canonical singularities and
is smooth at $\text{Bs} |-2K_X|$. 
Let $\mu: \t{S} \to S$ be the minimal resolution.
Since $y$ is a smooth point of $S$, 
we denote by $y$ the point of $\t{S}$ corresponding to $y$.
Note that $y$ is a base point of $|2K_{\t{S}}|$.
By [RI, Theorem 1] (or [KaMa, Theorem (1)]), 
there is an effective divisor $C$ through $y$ such that
$0\leq K_{\t{S}}.C <2$. But since $y \in C$, $K_{\t{S}}.C >0$ and 
hence $K_{\t{S}}.C=1$. 
Furthermore by [RI, Theorem 1] (or [KaMa, Theorem (2)]),
$0\leq C^2 \leq \frac 1 {K_{\t{S}}^2} < \frac 14$
and hence $C^2 =0$.
But $K_{\t{S}} .C=1$ and $C^2 =0$ contradict the genus formula
$2p_a(C) -2 = C.(K_{\t{S}} + C)$.
Hence $|K_S + K_S|$ is free and also $|-2K_X|$ is free.
\qed
\enddemo
By the proof of Theorem 4.1, we have the following:

\proclaim{Corollary 4.2}
Replace (3) by the assumption that $(-K_X)^3 \geq 1$ in Theorem 4.2.
Then a general member of $|-2K_X|$ has only canonical singularities.
Furthermore assume that $h^0 (-K_X)\geq 1$. Then $|-2K_X|$ is free.
\endproclaim
\endcomment

\proclaim{Proposition 4.1}
Let $X$ be a weak $\Bbb Q$-Fano $3$-fold with $I(X)=2$ such that $|-2K_X|$
is free.
Let $P$ be an index $2$ point such that there is no curve $l$ through $P$
such that $-K_X.l=0$.
Let
$f:Y\to X$ an extremal divisorial contraction from a $3$-fold
with only terminal singularities such that
\roster
\item
$f$-exceptional divisor is a prime $\Bbb Q$-Cartier divisor.
We call it $E$; 
\item
$P:=f(E)$ and
$-K_Y = f^*(-K_X) -\frac 12 E$;
\item $(-K_Y)^3 >0$.
\endroster
Then $Y$ is a weak $\Bbb Q$-Fano $3$-fold.
\endproclaim

\demo{Proof}
By the assumption that there is no curve $l$ through $P$
such that $-K_X.l=0$,
$\text{Bs} |-2K_X - P|$ is a finite set of points near $P$.
So by $H^0 (-2K_Y) \simeq H^0 (\Cal O (-2K_X) \otimes m_P)$, we know 
$-K_Y$ is nef.  So by (3), it is also big and we are done. 
\qed
\enddemo

\head 5. Solution of the equations of Diophantine type for
a $\Bbb Q$-Fano $3$-fold with $I(X)=2$ \endhead

\proclaim{Theorem 5.0}
Let $X$ be a $\Bbb Q$-factorial
$\Bbb Q$-Fano $3$-fold with the following properties:
\roster 
\item $\rho(X)=1$;
\item $I(X)=2$;
\item $F(X)=\frac 12$;
\item $h^0 (-K_X) \geq 4$;
\item there exists an index $2$ point $P$ such that
$$(X,P)\simeq (\{xy +f(z^2,u)=0\} / \Bbb Z _2 (1,1,1,0), o)$$
with $\text{ord}f(Z,U)=1$.
\endroster
Let $f:Y \to X$ be the weighted blow up at $P$ with weight
$\frac 12 (1,1,1,2)$.
Then $Y$ is a weak $\Bbb Q$-Fano $3$-fold with $I(Y)=2$
(and hence we can ran the program as in Section 3).
We use the notation as in there.
Then $z\leq u$ and there is at most one flip while 
$Y \dashrightarrow Y'$ and $a_i=2$ for $i$ such that 
$Y_i \dashrightarrow Y_{i+1}$ is a flip.
Furthermore
we figure out the solutions of equations in Section 3
as in the following tables:

See the tables in the main theorem for the explanation about the notation.

$$\vbox{
\offinterlineskip
\halign{\strut#&&\vrule#&\quad\hfil#\hfil\quad\cr
\multispan7 \hfil Table 1. $f'$ is of type $E_1$ and $u=z+1$ \hfil \cr
\noalign{\hrule}
&&$h$&&$(-K_X)^3$&&$N$&&$e$&&$n$&&$z$&&$l_C$&&$X'$&\cr
\noalign{\hrule}
&&$\bigcirc 6$&&$7$&&$2$&&$7$&&$0$&&$4$&&$35$&&$[5]$&\cr
\noalign{\hrule}
&&$\bigcirc 6$&&$\frac{15}2$&&$3$&&$7$&&$0$&&$2$&&$9$&&$[2]$&\cr
\noalign{\hrule}
&&$\bigcirc 6$&&$\frac{15}2$&&$3$&&$6$&&$1$&&$4$&&$30$&&$[5]$&\cr
\noalign{\hrule}
&&$\bigcirc 7$&&$\frac{17}2$&&$1$&&$6$&&$0$&&$3$&&$36$&&$\Bbb P^3$&\cr
\noalign{\hrule}
&&$\bigcirc 7$&&$9$&&$2$&&$6$&&$0$&&$2$&&$18$&&$[3]$&\cr
\noalign{\hrule}
&&$\bigcirc 7$&&$9$&&$2$&&$5$&&$1$&&$3$&&$32$&&$\Bbb P^3$&\cr
\noalign{\hrule}
&&$\bigcirc 7$&&$\frac{19}2$&&$3$&&$5$&&$1$&&$2$&&$15$&&$[3]$&\cr
\noalign{\hrule}
&&$\bigcirc 7$&&$\frac{19}2$&&$3$&&$4$&&$2$&&$3$&&$28$&&$\Bbb P^3$&\cr
\noalign{\hrule}
&&$\bigcirc 8$&&$\frac{21}2$&&$1$&&$6$&&$0$&&$1$&&$6$&&$B_3$&\cr
\noalign{\hrule}
&&$\bigcirc 8$&&$\frac{21}2$&&$1$&&$5$&&$0$&&$2$&&$27$&&$10$&&$Q_3$&\cr
\noalign{\hrule}
&&$\bigcirc 8$&&$11$&&$2$&&$4$&&$1$&&$2$&&$24$&&$Q_3$&\cr
\noalign{\hrule}
&&$8$&&$\frac{23}2$&&$3$&&$3$&&$2$&&$2$&&$21$&&$Q_3$&\cr
\noalign{\hrule}
&&$\bigcirc 9$&&$\frac{25}2$&&$1$&&$5$&&$0$&&$1$&&$10$&&$B_4$&\cr
\noalign{\hrule}
&&$\bigcirc 10$&&$\frac{29}2$&&$1$&&$4$&&$0$&&$1$&&$14$&&$B_5$&\cr
\noalign{\hrule}
&&$\bigcirc 10$&&$15$&&$2$&&$3$&&$1$&&$1$&&$12$&&$B_5$&\cr
\noalign{\hrule}
}}$$

\comment
$$m=0.$$
\endcomment

$$\vbox{
\offinterlineskip
\halign{\strut#&&\vrule#&\quad\hfil#\hfil\quad\cr
\multispan7 \hfil Table 2. $f'$ is of type $E_1$ and $z=u=1$ \hfil \cr
\noalign{\hrule}
&&$(-K_X)^3$&&$N$&&$e$&&$l_C$&&$X'$&\cr
\noalign{\hrule}
&&$\frac72$&&$3$&&$10$&&$1$&&$V_6$&\cr
\noalign{\hrule}
&&$4$&&$4$&&$8$&&$2$&&$V_8$&\cr
\noalign{\hrule}
&&$\frac92$&&$5$&&$6$&&$3$&&$V_{10}$&\cr
\noalign{\hrule}
&&$5$&&$6$&&$4$&&$4$&&$V_{12}$&\cr
\noalign{\hrule}
&&$\frac{11}2$&&$7$&&$2$&&$5$&&$V_{14}$&\cr
\noalign{\hrule}
&&$6$&&$8$&&$0$&&$6$&&$V_{16}$&\cr
\noalign{\hrule}
}}$$

$$h=4 \teb{and} n=0.$$

$$\vbox{
\offinterlineskip
\halign{\strut#&&\vrule#&\quad\hfil#\hfil\quad\cr
\multispan7 \hfil Table 3. $f'$ is of type $E_2 \sim E_{12}$  \hfil \cr
\noalign{\hrule}
&&$h$&&$(-K_X)^3$&&$N$&&$e$&&$n$&&type of $f'$ and $X'$&\cr
\noalign{\hrule}
&&$\bigcirc 4$&&$\frac 52$&&$1$&&$15$&&$0$
&&$E_5$ or $E_{11}$, $(-K_{X'})^3 =\frac 52$, $I(X')=2$&\cr
\noalign{\hrule}
&&$\bigcirc 4$&&$3$&&$2$&&$12$&&$0$&&$E_9$, $V_4$ &\cr
\noalign{\hrule}
&&$4$&&$4$&&$4$&&$9$&&$3$&&$E_2$, $V_{10}$&\cr
\noalign{\hrule}
&&$4$&&$\frac92$&&$5$&&$12$&&$3$&&$E_6$, $V_{16}$&\cr
\noalign{\hrule}
}}$$

$$z=u=1.$$

$$\vbox{
\offinterlineskip
\halign{\strut#&&\vrule#&\quad\hfil#\hfil\quad\cr
\multispan7 \hfil Table 4. $f'$ is of type $C$ \hfil \cr
\noalign{\hrule}
&&$h$&&$(-K_X)^3$&&$N$&&$e$&&$n$&&$\text{deg} \Delta$&\cr
\noalign{\hrule}
&&$5$&&$\frac{11}2$&&$3$&&$8$&&$0$&&$8$&\cr
\noalign{\hrule}
&&$5$&&$6$&&$4$&&$7$&&$1$&&$6$&\cr
\noalign{\hrule}
&&$5$&&$\frac{13}2$&&$5$&&$6$&&$2$&&$4$&\cr
\noalign{\hrule}
&&$5$&&$7$&&$6$&&$5$&&$3$&&$2$&\cr
\noalign{\hrule}
&& $5$&&$\frac{15}2$&&$7$&&$4$&&$4$&&$0$&\cr
\noalign{\hrule}
&&$\bigcirc 6$&&$\frac{13}2$&&$1$&&$7$&&$0$&&$7$&\cr
\noalign{\hrule}
&&$6$&&$7$&&$2$&&$6$&&$1$&&$6$&\cr
\noalign{\hrule}
&&$6$&&$\frac{15}2$&&$3$&&$5$&&$2$&&$5$&\cr
\noalign{\hrule}
&&$6$&&$8$&&$4$&&$4$&&$3$&&$4$&\cr
\noalign{\hrule}
&&$6$&&$\frac{17}2$&&$5$&&$3$&&$4$&&$3$&\cr
\noalign{\hrule}
&&$6$&&$9$&&$6$&&$2$&&$5$&&$2$&\cr
\noalign{\hrule}
&&$6$&&$\frac{19}2$&&$7$&&$1$&&$6$&&$1$&\cr
\noalign{\hrule}
&&$6$&&$10$&&$8$&&$0$&&$7$&&$0$&\cr
\noalign{\hrule}
&&$10$&&$\frac{29}2$&&$1$&&$6$&&$0$&&$0$&\cr
\noalign{\hrule}
}}$$

$$\tea{If} h=5, \teb{then}
z=u=2 \teb{and} X'\simeq \Bbb F_{2,0}.$$
$$\tea{If} h=6, \teb{then}
z=u=1 \teb{and} X'\simeq \Bbb P^2.$$
$$\tea{If} h=10, \teb{then}
z=1, u=2 \teb{and} X'\simeq \Bbb P^2.$$

$$\vbox{
\offinterlineskip
\halign{\strut#&&\vrule#&\quad\hfil#\hfil\quad\cr
\multispan7 \hfil Table 5. $f'$ is of type $D$ \hfil \cr
\noalign{\hrule}
&&$h$&&$(-K_X)^3$&&$N$&&$e$&&$n$&&$\text{deg} F$&\cr
\noalign{\hrule}
&&$4$&&$\frac 92$&&$5$&&$9$&&$0$&&$6$&\cr
\noalign{\hrule}
&&$4$&&$5$&&$6$&&$8$&&$1$&&$8$&\cr
\noalign{\hrule}
&&$\bigcirc 5$&&$\frac92$&&$1$&&$9$&&$0$&&$3$&\cr
\noalign{\hrule}
&&$5$&&$5$&&$2$&&$8$&&$1$&&$4$&\cr
\noalign{\hrule}
&&$5$&&$\frac{11}2$&&$3$&&$7$&&$2$&&$5$&\cr
\noalign{\hrule}
&&$5$&&$6$&&$4$&&$6$&&$3$&&$6$&\cr
\noalign{\hrule}
}}$$
$$z=u=2\teb{in case} h=4.$$
$$z=u=1\teb{in case} h=5.$$

If $f'$ is a crepant divisorial contraction, then 
$$h=4,(-K_X)^3 =\frac 52, N=1, z=1 \teb{and} u=2.$$ 
\endproclaim

\definition{Remark}
We discuss the geometric realization in Section 6.
\enddefinition

\demo{Proof}
By (4) and Corollary 1.4, we have $(-K_X)^3 >2$.
Furthermore $(-K_Y)^3 =(-K_X)^3 -\frac 12 >0$.
Hence by Proposition 4.1, $Y$ is a weak $\Bbb Q$-Fano $3$-fold.
We can easily check that $I(Y)=2$ by calculating the weighted blow up. 

We run the program as in Section 3.
By the assumption that $h^0(-K_X) \geq 4$ and
the exact sequence 
$$0\to \Cal O_Y (-K_Y -E) \to \Cal O_Y (-K_Y) \to \Cal O_E (1) \to 0,
$$ we have $|-K_Y-E|\not = \phi$.
Hence by Claim 3.11, we have $z\leq u$.

First assume that Case 5 occurs.
Since $u \in \frac {\Bbb N} 2$ and $E.l \in \Bbb N$,
we have $u=\frac 12, 1, 2$ by $u(E.l)=2$.
Furthermore since $z(-K_Y)^3 =u$, $(-K_Y)^3 > \frac 32$ and $z\leq u$,
we have $z=1$, $u=2$ and $(-K_Y)^3 =2$.
Hence we are done in this case.

\proclaim{Claim 5.1} 
$E_i$ is a Cartier divisor for any $i$.
In particular $a_i$ is an even integer.
\endproclaim

\demo{Proof}
Assume that $g_0$ is a flop. By Proposition 1.5,
$E_1$ is a Cartier divisor since $E$ is a Cartier divisor.
The latter half follows from Proposition 2.1 (1).
\comment
it suffices to prove that
the property being a Cartier divisor is preserved by the flop.
Furthermore by decomposing the flop to analytic flops, 
we may assume that the flopping curve is irreducible. Let $l$ be the
flopping curve.

If there are only Gorenstein points on $l$, the assertion is clear
by the fact that the analytic structure is unchanged by the flop
[Kol1, Theorem 2.4 and Remark 5.1].

If there are $\frac 12 (1,1,1)$-singularities on $l$,
there are the following possibilities by [Kol1, pp 27--28]:

Let $h:Y\to Z$ be the flopping contraction of $l$, $h^+: Y_1 \to Z$
its flop and $Q:=h(l)$. 
In the below, $x, y, z, u$ are coordinates of $\Bbb C^4$.
 
$$(Q\in Z) 
\simeq (o \in (xy + z^2 - u^{4s} = 0 \subset \Bbb C^4/\Bbb Z_2(1,1,0,1)).
\tag a$$
$h$ is the blow up along $(x=z-u^{2s}=0)$ and
$h^+$ is the blow up along $(y=z-u^{2s}=0)$;

$$(Q\in Z) 
\simeq (o \in (x^2+y^2 -z^2 - u^{4s} = 0 
\subset \Bbb C^4/\Bbb Z_2(1,0,1,1)). \tag b$$
$h$ is the blow up along $(x-z=y-u^{2s}=0)$ and
$h^+$ is the blow up along $(x+z=y-u^{2s}=0)$.

For these two cases, we can see by direct computations
that there is a Weil divisor $D$ on $Y$ such that $D.l=\frac 12$.
Let $L$ be a Cartier divisor on $Y$ and $a:=L.l (\in \Bbb Z)$.
Since we consider analytically locally, 
$\text{Pic} Y \simeq \text {Pic} l\simeq \Bbb Z$.
Hence $L\sim a(2D)$ since $2D$ is a Cartier divisor.
After the flop, the linear equivalence is preserved
and hence $L^+\sim a(2D^+)$, where $L^+$ (resp. $D^+$) is the strict
transform of $L$ (resp. $D$). 
So $L^+$ is also Cartier.
\endcomment
If $g_i$ is a flip, there is no non Gorenstein point on the flipped curves.
Hence $E_i$ is Cartier by induction for $i$.
\qed
\enddemo

Note that by Lemma 3.2 (5) and Claim 5.1, 
once we prove that $a_i =2$ if $a_i>0$,
we see that there is at most one flip while 
$Y \dashrightarrow Y'$. 
\enddemo
\definition{Case 1}
In this case we first show that $F(X')\geq 1$.
 
In fact by (3-0-1), we have $-K_{X'}\equiv \frac u z f'(\t{E})$.
In this, $f'(\t{E})$ is Cartier and $u\geq z$. Hence the assertion holds.
Furthermore by [I3] and [San2], $F(X')=1, \frac 32, 2, \frac 52, 3$ or $4$. 

We note that 
by Proposition 2.2, we have $(-K_{E'})^2 = 8(1-g(\o{C}))-2m$ with some
non negative integer $m$.

By $z+1=uk$ and $z\leq u$, we have $z+1=u$ or $z=u=1$. 

First assume that $z+1=u$. 
Define $a\in \Bbb Z$ 
by the formula $f(\t{E})=aH$, where $H$ is a primitive Cartier
divisor of $X'$. Then $F(X')= a\frac {z+1} z$.
Hence $z= 1,2,3,4$ and if
$z=1$, then $F(X') =2$ or $4$,
if $z=2$, then $F(X') = \frac 32$ or $3$,
if $z=3$, then $F(X')=4$,
or if $z=4$, then $F(X') =\frac 52$.
\comment
Note that if $X'$ is smooth, then $Y'$ must be smooth.
For if $Y'$ is not smooth, then every $\frac 12(1,1,1)$-singularity
is contained in $E'$.
But the image of a fiber containing $\frac 12(1,1,1)$-singularity
must be singular by Proposition 2.2, a contradiction.
\endcomment
But we will prove that
the case that $z=1$ and $F(X')=4$ does not occur.
For otherwise, let $H'$ be the strict transform of ${f'}^* H$ on $Y$.
Then we have $-K_Y \equiv 2H' +E$ and hence $-K_X \equiv 2f(H')$,
a contradiction to $F(X)=\frac 12$. 

Assume $a_i \geq 4$ for some $i$.
Note that $a_i u > z$ by $u\geq z$.
By (3-1-5'),
$e\leq (k+2)^2 - 2(4-k)^2 <0$, a contradiction.

Set $n:=\sum n_i$. We obtain the following:

$$(-K_X)^3=\frac{9+n}2 +\frac 1{u^2}(-K_{X'})^3 \tag {5-1-1}$$ 
obtained by (3-1-1').
 
$$e+n=9-\frac{u-1}{u^3}(-K_{X'})^3 \tag {5-1-2}$$
obtained by (3-1-5').

$$(-K_{X'}.C)=\frac {u-1}u (-K_{X'})^3 -(3+n)u \tag {5-1-3}$$
obtained by (3-1-1') and (3-1-3').

$$-2(-K_{X'})^3+4(-K_{X'} . C) + 2(-K_X)^3 -n+3=4g(\o{C})+m
\tag {5-1-4}$$ 
obtained by (3-1-4).

Use (5-1-4) for the bound of $n$.

By (3-0-1), we have $\t{E}.l = 1$ for a general fiber $l$ of $E'$.
If $E'$ contain an index $2$ point, 
then there is a component $l'$ of a fiber 
such that $-K_{Y'}.l' = \frac 12$ by Proposition 2.2. 
So we have $\t{E}.l'= \frac 12$. But this contradicts the fact that 
$\t{E}$ is a Cartier divisor. Hence $E'$ contains no index $2$ point.
This fact and information from $X'$ determine $N$.
Hence we can easily figure out the solutions.

Next assume $z=u=1$.
By Claim 3.7 and Claim 5.1, $a_i \geq 4$ if $a_i >0$.
Assume that $a_i \geq 6$ for some $i$.
By (3-1-5'),
$e\leq (k+2)^2 - 3(6-k)^2 <0$, a contradiction.
Hence we must have $a_i = 4$ for all $i$ such that $Y_i
\dashrightarrow Y_{i+1}$ is a flip.
By setting $n:=\sum n_i$, we rewrite (1-1) $\sim$ (1-4) as follows:

$$e+8n=16-(-K_{X'})^3 \tag {5-1-1'}$$
obtained by (3-1-5').

$$(-K_X)^3=6-\frac 14 e -\frac 32 n .\tag {5-1-2'}$$
obtained by (3-1-2').

$$(-K_{X'}. C)=6-\frac 12 e -6n . \tag {5-1-3'}$$
obtained by (5-1-2') and (3-1-3').

$$2(-K_X)^3 -5 =4g(\o{C})+m. \tag {5-1-4'}$$

By (5-1-3') and $(-K_{X'}. C)>0$, we must have $n=0$, i.e.,
there is no flip while $Y\dashrightarrow Y'$.
 
By (1-1) and (1-2), we deduce that 
$(-K_{X'})^3 =16-e>0$.
By (1-1) and (1-3), we have 
$(-K_{X'}. C)=\frac 12 (-K_{X'})^3 -2 >0$.
By these, $(-K_{X'})^3 = 6,8,10,12,14,16$.
\enddefinition

\proclaim{Claim 5.2}
$h^0(-K_X)=4$.
\endproclaim

\demo{Proof}
By $\t{E} \equiv -K_{Y'} - E'$, we have
$E \equiv -K_Y - \t{E'}$, where $\t{E'}$ is the strict transform of $E'$.
Since $E - (-K_Y - \t{E'})$ is a Cartier divisor,
we must have $E \sim -K_Y - \t{E'}$ by Lemma 3.2 (6).
Hence $h^0(-K_Y -E)=1$.
But by the exact sequence
$$0 \to \Cal O_Y(-K_Y -E)\to \Cal O_Y(-K_Y)\to \Cal O_E(-K_Y |_E) 
\to 0,$$
we have $$h^0(-K_X)=h^0 (-K_Y) \leq h^0(-K_Y -E) + h^0 (-K_Y |_E)=4.$$
So $h^0(-K_X) =4$.
\qed
\enddemo

Hence we have $N=\frac{16-e}2$.

We will prove that $X'$ is Gorenstein.
Assume that $X'$ is non Gorenstein.
If $F(X')=1$, then by [San2], $N-1\geq 8$, a contradiction.
Since $(-K_{X'})^3=16-e$, $F(X')>1$ does not hold by [San1].
Hence $X'$ is Gorenstein.

Next we prove that $F(X')=1$. By $(-K_{X'})^3=16-e$, we only have to disprove
that $F(X')=2$. Assume that $F(X')=2$.  
Let $H$ be the ample generator of $\text{Pic} X'$
and $H':= {f'}^* H$. This is a Cartier divisor on $Y'$ and so is
the strict transform $H''$ on $Y$ since $n=0$.
Since $H''\equiv \frac 12 (-K_Y +\t{E'}) \equiv (-K_Y)- \frac 12 E$, we have
$f^* f_* H'' = H'' + E$. By this, we know $f_* H''$ is a Cartier divisor on $X$
([KMM, Lemma 3-2-5 (2)]).
On the other hand, $f_* H'' \equiv -K_X$ and so $F(X)$ must be
an integer, a contradiction to the assumption of Theorem 5.0.
Hence we have $F(X')=1.$

So we obtain the solutions as in the table.
 
\definition{Case 2}
By Proposition 2.3, we obtain the following data:

In the below equations, the right sides are the values of the left sides
in case $f'$ is of type $E_i$ ($i=2\sim 12$).

$$(E')^3 =\overset{E_2}\to {1} \ \overset{E_3, E_4}\to{2} \ 
\overset{E_5} \to{4}
 \ \overset{E_6} \to {\frac 12}\ 
\overset{E_7}\to {\frac 12} \overset{E_8}\to {1} 
\overset{E_9}\to {\frac 32} \overset{E_{10}}\to {2}
\overset{E_{11}}\to {4} \overset{E_{12}}\to {\frac 92}.$$ 

$$(-K_{Y'})(E')^2 =\overset{E_2}\to {-2} \ \overset{E_3, E_4}\to{-2} \ 
\overset{E_5} \to{-2}
 \ \overset{E_6} \to {-\frac 32}\ 
\overset{E_7}\to {-\frac 12} \overset{E_8}\to {-1} 
\overset{E_9}\to {-\frac 32} \overset{E_{10}}\to {-2}
\overset{E_{11}}\to {-2} \overset{E_{12}}\to {-\frac 32}.$$

$$(-K_{Y'})^2 E'=\overset{E_2}\to {4} \ \overset{E_3, E_4}\to{2} \ 
\overset{E_5} \to{1}
 \ \overset{E_6} \to {\frac 92}\
\overset{E_7}\to {\frac 12} \overset{E_8}\to {1} 
\overset{E_9}\to {\frac 32} \overset{E_{10}}\to {2} 
\overset{E_{11}}\to {1} \overset{E_{12}}\to {\frac12}.$$

Assume that $f'$ is of type $E_2$.
By (3-2-5'),
we have $z+2u\leq 2k$.
On the other hand, we have $1+2z=uk \geq zk$.
Hence $z=u=1$ and $k=3$.
By (3-2-5') again, $\sum d_i a_i (a_i-1)=3$.
Since $a_i\geq 2$ if $a_i >0$, we have $a_i =2$ if $a_i>0$.
By setting $n:=\sum n_i$, we have $n=3$.
We can easily see that $e=9$, $(-K_X)^3 =4$ and
$(-K_{X'})^3 =10$.
By the assumption (4), we have $N=4$.
This also prove that $X'$ is Gorenstein and hence
$X'$ is $V_{10}$.

We will prove that $f'$ cannot be of type $E_3$ or $E_4$.
Assume that $f'$ is of type $E_3$ or $E_4$.
Similarly to the above case, we have $k=2$ and $\sum d_1 a_i (a_i -1)=1$
using (3-2-5'). But by Claim 3.9, we have $a_i\geq 4$ if $a_i>0$, 
a contradiction.

If $f'$ is of type $E_5 \sim E_{11}$,
then we can figure out the solution similarly.
\enddefinition

By these we can obtain the solutions.

\definition{Case 3}
By Proposition 2.4, $X'$ has at worst ordinary double points as singularities.
Hence $X'\simeq \Bbb P^2 \ \text{or} \ \Bbb F_{2,0}$ and
$L= {f'}^*\Cal O_{\Bbb P^2}(1)$ if $X' \simeq \Bbb P^2$ or
$L= {f'}^*(\Cal O_{\Bbb P^3}(1)|_{\Bbb F_{2,0}})$ if $X' \simeq \Bbb F_{2,0}$.

Assume $a_i \geq 4$ for some $i$.
Note that $a_i u > z$ by $u\geq z$.
By (3-3-2'), 
$m^2 e < (2m+1)^2 - 2(4m - 1)^2 <0$, a contradiction.
Hence $a_i=2$ for all $i$ such that $Y_i \dashrightarrow Y_{i+1}$ is a flip.

By setting $n:=\sum n_i$, we have the following:

$$(-K_X)^3 =\frac 12 + m(4m^2 +6m +3)-\frac n2 (2m-1)^3 -m^3 e. \tag {5-3-1}$$

$$m z^2\{ (2m+1)^2 - n(2m-1)^2 -m^2 e\} =2l^2=
\overset{\Bbb P^2} \to {2} \ \overset{\Bbb F_{2,0}} \to {4}. \tag {5-3-2}$$

$$mz\{ 2(2m+1)(m+1)-2n (2m-1)(m-1)-m^2 e\}
= \overset{\Bbb P^2} \to {12 - \Delta.l} \ 
\overset{\Bbb F_{2,0}} \to {16 - \Delta.l}. 
\tag {5-3-3}$$
By Claim 3.11, we have $m=1$ or $2$.

If $m=2$, we can easily figure out the solution.

Assume that $m=1$. Then we have $3$ sequences of solutions:
\roster
\item $X'\simeq \Bbb P^2$, $z=1$,
$n+e=7$, $\Delta .l=e$ and $h^0(-K_X)=\frac{25+n-N}4$;

\item $X'\simeq \Bbb F_{2,0}$, $z=1$,
$n+e=5$, $\Delta .l=4+e$ and $h^0(-K_X)=\frac{29+n-N}4$;

\item $X'\simeq \Bbb F_{2,0}$, $z=2$,
$n+e=8$, $\Delta .l=2e-8$ and $h^0(-K_X)=\frac{23+n-N}4$.
\endroster

If $X'\simeq \Bbb P^2$ and $Y'$ has an index $2$ point 
(resp. If $X'\simeq \Bbb F_{2,0}$ and $\text{aw}(Y')>2$),
then there is a fiber containing a component $l$ 
such that $-K_{Y'}.l=\frac 12$ by Proposition 2.4.
But these cases does not occur.
For otherwise we have $\t{E}.l=\frac z{2u} <1$, a contradiction to 
that $\t{E}$ is a Cartier divisor.
Hence for (1) and (2) (resp. (3)), we have
$N-n=1$ (resp. $N-n=3$) since $\text{aw}(Y') =\text{aw}(Y) -n =N-n-1$.
But if (2) and $N-n=1$ hold, $Y'$ must be Gorenstein,a contradiction
to Proposition 2.4. 
Hence we figure out the solutions as in the table.
\enddefinition

\definition{Case 4}
Similarly to Case 3, we can prove 
that $a_i=2$ for all $i$ such that $Y_i \dashrightarrow Y_{i+1}$ is a flip
using (3-4-2').

By setting $n:=\sum n_i$, we rewrite (3-4-1')$\sim$(3-4-3') as follows:

$$(-K_X)^3=\frac 12 +2m(m+1)+\frac 12 n(2m-1)^2. \tag {5-4-1}$$

$$(2m+1)^2 =n(2m-1)^2 +m^2 e. \tag {5-4-2}$$

$$z\{ m(2m+1)+ nm(2m-1) \}= \text{deg}  F. \tag {5-4-3}$$
By Claim 3.11, we have $m=1, \frac 32, 2$ or $3$.
 
We can easily see that there is no solution for $m=\frac 32, 2$ or $3$.

If $m=1$, then we have 
$n+e=9$, $(-K_X)^3=\frac{n+9}2$ and $z(3+n)=\text{deg} F$.
Since $h^0 (-K_X)=3+\frac{9+n-N}4 \geq 4$, we have $N-n = 1 \ \text{or} \ 5$.
If $N-n = 1$, then $Y'$ is Gorenstein. Hence by the primitivity of $L$,
$z=1$.
If $N-n = 5$ and $u=z=1 \ \text{or}\ 3$,
$L\not \sim z(-K_{Y'}-\t{E})$ since the right side is not Cartier. 
By Riemann-Roch theorem,
$\chi (\Cal O(L)) - \chi (O(z(-K_{Y'}-\t{E}))) = \frac 12$,
a contradiction.
Hence if $N-n = 5$, then $z=2$ and so $n=0$ or $1$ by 
$z(3+n)= \text{deg}  F$.

We will prove $n\leq 3$.
If $n=4$, then $\deg \ F = 7$, a contradiction to Proposition 2.5.
If $n=5$, then $Y' \to X'$ is a quadric bundle over a $\Bbb P^1$
by Proposition 2.5.
But then $(-K_{Y'})^3$ must be a multiple of $8$, a contradiction.
If $n=6$, then $Y' \to X'$ is a $\Bbb P^2$-bundle over a $\Bbb P^1$
by Proposition 2.5.
But then $(-K_{Y'})^3$ must be $54$, a contradiction.

Hence we obtain the solutions as in the table.
\qed
\enddefinition
Now we complete the proof of Theorem 5.0.

In the next section, we give examples for the cases with mark $\bigcirc$
in the table and prove the non-existence of the cases
$N=7,8$ in Table 2 and $h=6$ and $N=8$ in Table 4.
This will complete the proof of the main theorem.

\head 6. Completion of the proof of the main theorem and examples \endhead

\comment
Next two lemmas are needed for this purpose.

\proclaim{Lemma 5.8}
Let $X$ be as in Theorem 4.1. 
Assume that $g_0$ is a flop. (see Lemma 3.2 for this notation.) 
Then $|-2K_Y|$ is free.
\endproclaim

\demo{Proof}
Let $h: Y \to Z$ be
the flopping contraction associated to $g_0$.
If $N\geq 2$,
then by the tables of Theorem 4.1,
$Z$ is a $\Bbb Q$-Fano $3$-fold
which satisfies the assumptions of Theorem 5.0.
Hence $|-2K_Y|$ is free.
Assume that $N=1$. Then by the tables of Theorem 4.1,
$Z$ is Gorenstein and $(-K_Z)^3 \geq 2$.
So $|-2K_Z|$ is free by [I3, the proof of Chap.1, Theorem 6.3] and 
also $|-2K_Y|$ is free.
\qed
\enddemo
\endcomment

We state a theorem proved by T. Minagawa which we need.

\proclaim{Theorem 6.0 (T. Minagawa)}
Let $X$ be a $\Bbb Q$-Fano $3$-fold 
(resp. weak $\Bbb Q$-Fano $3$-fold) with $I(X)=2$.
Assume that there exists a smooth member of $|-2K_X|$.
Then there exists a flat family $\frak f : \frak X \to (\Delta , 0)$ 
over a $1$-dimensional disk $(\Delta , 0)$
such that $X\simeq \frak f ^{-1} (0)$ and
$\frak f^{-1} (t)$ is a $\Bbb Q$-Fano $3$-fold
(resp. a weak $\Bbb Q$-Fano $3$-fold) with only
ODP's, QODP's or $\frac 12 (1,1,1)$-singularities as its singularities 
for $t\in \Delta \backslash \{ 0 \}$,
where ODP (resp. QODP) means a singularity analytically isomorphic to
$\{ xy +z^2 +u^2 =0 \subset \Bbb C^4 \}$
(resp. $\{xy +z^2 +u^2=0\} \subset \Bbb C^4 / \Bbb Z _2 (1,1,1,0)$.)
\endproclaim
\demo{Proof}
See [Mi2, Theorem 2.4].  
\qed
\enddemo

From now on we assume that $X$ is a Fano $3$-fold as in the main theorem.
\definition{Table 1}
First we construct examples for the case that $f'$ is of type $E_1$
and $u=z+1$. We can treat all the cases at a time except
$h=8$ and $N=3$. We don't know whether the case that $h=8$ and $N=3$
occurs or not.
 
Let $S$ be a smooth Cartier divisor in $X'$ such that
$S\equiv \frac z{z+1} (-K_{X'})$. 
We can take such a $S$ by [San2, Remark 4.1]. $S$ is a del Pezzo surface.
We represent $S$ as blowing up at $r$ points of $\Bbb P^2$ in general position,
where $r:=e+n$.
Let $E_1, \dots , E_r$ be its exceptional curves and $l$ the total transform
of a line in $\Bbb P^2$.
Let $D:=l+E_1+\dots +E_n$ and $C:=-K_{X'}|_{S}-D$.
Then we will show that $|C|$ is free.
By computing intersection numbers with $(-1)$-curves,
we can check that $C$ is nef in any case in the table
except the case that $h=8$ and $N=3$.
Let $M:= C- K_S$. Check that $M^2 >4$.
Hence if $C$ is not free, there is an effective divisor $l$
such that $M.l =1$ and $l^2 =0$ whence $-K_S .l =1$ by Reider's theorem [RI].
But $l.(K_S + l)=-1$ is a contradiction. So $|C|$ is free.

Hence we denote a general smooth curve in $|C|$ also by $C$.
Let $f' :Y' \to X'$ be the blow up along $C$ and $E'$ the exceptional divisor.
Let $\t{E}$ be the strict transform of $S$ and $B:= 2(-K_{Y'})-\t{E}$.
Let ${E_i}'$ be the inverse image of $E_i$ for $i=1 \dots n$.
We will check that $B$ is nef and big.
Let $A$ be a Cartier divisor numerically equivalent to
$\frac {z+2}{z+1} (-K_{X'})$.
Since $B\sim {f'}^* A - E'$, we have only to show $|A-C|$ is free.
Since $A-S$ is free and $C\subset S$, $\text{Bs}  |A-C| \subset S$.
By the exact sequence 
$$0\to \Cal O_{X'}(A-S) \to \Cal O_{X'}(A)\to \Cal O_{S}(A|_S)\to 0$$ 
and the KKV
vanishing theorem, we see that $\text{Bs}  |A-C| = \text{Bs} |A|_S - C|$.
We can check that $|A|_S - C|$ is free by the same way as 
the check of the freeness of $|C|$. We also know that $A|_S - C$
is numerically trivial only for $E_i$'s ($i=1\dots n$).

Hence $B$ is free and ${E_i}'$'s ($i=1\dots n$)
 are numerically trivial for $B$.
In particular two extremal rays of $Y'$ are generated by 
the class of a fiber of $f'$ and the class of ${E_i}'$.
Let $R_1$ be the extremal ray generated by the class of ${E_i}'$.
Then we will show that $\supp R_1 = \cup {E_i}'$.
Since $\t{E}.{E_i}' = -2 <0$, $\supp R_1 \subset \t{E}$.
By $-K_{Y'}.{E_i}' =-1 <0$, it is enough to show that
$\text{Bs}  |-K_{Y'}||_{\t{E}} = \cup {E_i}'$.
For this, we have only to see that
$\text{Bs} |-K_{X'} - C||_{S} = \cup E_i$.
By the exact sequence 
$$\ex{X'}{-K_{X'}-S}{X'}{-K_{X'}}{S}{-K_{X'}|_S} $$ 
and the KKV vanishing theorem, we see that $\text{Bs}  |-K_{X'}-C||_{S} = 
\text{Bs}  |-K_{X'}|_S -C|= \text{Bs} |D|$.
By the definition of $D$, it is clear that
$\text{Bs} |D| = \cup E_i$. So we are done.

By this, $R_1$ is a flipped ray. Observe that $\Cal N_{{E_i}' / Y'} \simeq
\Cal O_{\Bbb P^1}(-1) \oplus \Cal O_{\Bbb P^1}(-2)$.
Let $Y' \dashrightarrow Y_1$ be the inverse of the flip
and ${E_i}'$ the strict transform of $E_i$ on $Y_1$ for $i=n+1\dots r$.
Then we can easily show that 
$\cup _{i=n+1}^{r} E_i$ coincides with the support of the flopped ray of $Y_1$.
Let $Y_1 \dashrightarrow Y$ be the inverse of the flop
and $E$ the strict transform of $\t{E}$ on $Y$.
Then $(E, -K_Y |_E) \simeq (\Bbb P^2, \Cal O_Y (1))$
and we can contract it. Set $X$ the target of the  contraction.
Then $X$ is what we want.
\enddefinition

\definition{Table 2}
\enddefinition

\comment
\proclaim{Claim 1}
Assume that the Fano index of $X'$ is $1$.
(We need this assumption only in case that $h=4$ and $N=8$.)
Let $l$ be a line on $X'$ intersecting $C$
and $l'$ the strict transform of $l$ on $Y'$.
Then one of the following holds:
\roster
\item 
$-K_{Y'}.l'= 0 \teb{and} E'.l'=1$;
\item
$-K_{Y'}.l'= \frac 12 \teb{and} E'.l'=\frac 12$.
$\t{E} \cap l' = \phi$.
\endroster
Furthermore $C\cap l$ is one point and 
$C$ and $l$ are not tangent.
\endproclaim

\demo{Proof}
by the formula $-K_{Y'}.l' = -K_{X'}.l' - E'.l'$ and the fact that
$-K_{Y'}$ is nef (note that $n=0$),
we have $-K_{Y'}.l'= 0 \teb{and} E'.l'=1$ or
$-K_{Y'}.l'= \frac 12 \teb{and} E'.l'=\frac 12$.
Assume that the latter case occurs. Then we have $\t{E}.l' =0$.
By Proposition 7.1, there is no flopped curve containing 
a $\frac 12 (1,1,1)$-singularity. 
Hence there is no $\frac 12 (1,1,1)$-singularity on $\t{E}$.
So $l' \not \subset \t{E}$ and hence $\t{E} \cap l' = \phi$.

The last assertion clearly holds in case (2)
and it holds by Proposition 7.1 in case (1).
\qed
\enddemo
\endcomment

\comment
\definition{$N=3$}
First we construct an example for the case that $N=3$.
Consider a weighted projective space
$\Bbb P(1^4,2^2)$.
Let $x_i$ ($i=0,1,2,3$) be the coordinates of weight $1$
and $y_i$ ($i=1,2$) be the coordinates of weight $2$.
Take a weighted complete intersection
$(((3,4)\subset \Bbb P (1^4,2^2))$ which is defined by the following
$2$ equations:

$$f=l_1 y_1 + l_2 y_2 + q x_3 $$ and 
$$g= a y_1 y_2 + q_1 y_1 + q_2 y_2 + r x_3 ,$$
where $a$ is a constant,
$l_i$ is a linear form, $q$ and $q_i$ are polynomials of degree $2$
and $r$ is a polynomial of degree $3$.
Note that $Z:= (f=0) \cap (g=0)$ contains $\o{E} :=(y_1 = y_2 = x_3 =0)$
which is $\Bbb P^2$.
By taking $a$, $l_i$, $q$, $q_i$ and $r$ generally,
we assume that $Z$ is smooth outside $\o{E}$ and its $\frac 12 (1,1,1)$-
singularities.
Then we can see that 
$Z$ has exactly $10$ ordinary double points on $\o{E}$ by easy computations.
Let $Y:= \pmb{Proj}  \bigoplus _{i=0}^{\infty} \Cal O_Z (m\o{E})$.
Then $Y$ is smooth outside its $\frac 12 (1,1,1)$-singularities
and the inverse image $E$ of $\o{E}$ remains $\Bbb P^2$.
We can contract $E$ to a $\frac 12(1,1,1)$-singularity
and let $X$ be the target. Then $X$ is an example as we want.
($\rho (X) = 1$????)
\enddefinition
\endcomment

\comment
\definition{$N=5,6,7$}
Next we consider the case that $N=5,6,7$.
Let $P$ be a Gorenstein singularity on $C$.
Let $g: Z\to X'$ be the blow up of $P$ and $F$ the exceptional divisor.
Since $(P\in X') \simeq (o\in ((xy+zw=0) \subset \Bbb C^4))
\ \text{or} \  (o\in ((xy+z^2 + w^3=0) \subset \Bbb C^4))$ by Proposition 2.2
and $X'$ is $\Bbb Q$-factorial, $\rho(Z)=2$.
Since $-K_{X'}$ is very ample and $-K_Z = g^*(-K_{X'}) - F$,
$|-K_Z|$ is free and $(-K_Z)^3 >0$. Hence $Z$ is a weak Fano $3$-fold.
Starting by $g$, we consider the similar diagram as in the section 3 and
do similar calculations as there:
\enddefinition

$$\matrix
\ & Z & \dashrightarrow & {Z'} \\
\ & {g\swarrow} & \ &  {\searrow g'} \\
X' & \ & \ & \ & X'' & .
\endmatrix $$

Let $\t{F}$ be the strict transform of $F$ on $Z'$.
Then by a similar way,
we have $(-K_{Z'})^2 \t{F} = 2$, $(-K_{Z'}) (\t{F})^2 = -2$ 
and $(\t{F})^3 = 2-e'$, where $e'$ is a non negative integer.
Set $d:=(-K_{Z'})^3$ and input these into (1-1)$\sim$(5-1).
Then we obtain the following:

If $N=5$, then $e'=6$ and $g'$ is a conic bundle over $\Bbb P^2$
with $\text{deg} \ \Delta ' = 6$, where $\Delta '$ is the discriminant 
divisor for $g'$.

If $N=6$, then $e'=5$ and
$g'$ is of type $E_1$ and $X''\simeq \Bbb P^3$.
Let $F'$ be the exceptional divisor of $g'$. Then $F' \sim
3(-K_{Z'})-4\t{F}$.
For the center $C'$ of $g'$, $\text{deg}  C' = 8$ and $p_a (C') = 6$.

If $N=7$, then $e'=4$ and 
$g'$ is of type $E_1$ and $X''\simeq Q_3$.
(Since $X''$ is $\Bbb Q$-factorial, $X''$ is a smooth quadric.)
Let $F'$ be the exceptional divisor of $g'$. Then $F' \sim
2(-K_{Z'})-3\t{F}$.
For the center $C'$ of $g'$, $\text{deg}  C' = 8$ and $p_a (C') = 4$.

In case $N=7$, $Z'$ has $5$ Gorenstein singular points 
on the strict transform of $C$.
Since $X''$ is smooth, $C'$ must have
$5$ singular points by [Cu, Theorem 4].
But by $p_a (C') = 4$, this is impossible.

In case $N=6$, IN PREPARATION.
\endcomment

\comment
In case $N=5$, we consider the flopping contraction
$h:Z\to W$ (note that $e'>0$).
Then $W$ is a $(2,2,2)$ complete intersection in $\Bbb P^6$.
Let $\o{F}$ be the image of $F$ and
Let $\o{C}$ be the strict transform of $C$ on $W$.

\proclaim{Claim 2}
\roster
\item
Let $m$ be a flopping curve for $h$.
Then $m$ is the strict transform of a line on $X'$.
In particular 
$\o{F}\simeq F$, i.e., $\o{F}$ is a quadric in a linear space $\Bbb P^3$ in
$\Bbb P^6$.
\item
$\o{C}$ is a conic in a plane in $\Bbb P^6$
and no point which is
an image of a flopping curves is on it.
Furthermore $\o{F} \cap \o{C}$ is one point and the intersection is simple.
\endroster
\endproclaim

\demo{Proof}
The proof of (1) is similar to Lemma 5.9.
We will prove only (2) here.
For a flopping curve $m$, $g(m)$ is a line and hence by
Claim 1, $g^{-1} C \cap m = \phi$.
This proves the former half of (2).
The latter half of (2) is clear.
\qed
\enddemo
Let $x_0, x_1, \dots x_6$ the coordinates of $\Bbb P^6$.
Fix a $\Bbb P^3$ (resp. $\Bbb P^2$)
containing $\o{F}$ (resp. $\o{C}$) and denote it by $P$ (resp. $P'$).
We may assume that $P=(x_4 = x_5 = x_6 =0)$ and
$W$ is defined by $3$ homogeneous polynomials
$$\align
f_1= aQ+ x_4 l_{41} + x_5 l_{51} + x_6 l_{61} \\
f_2= bQ+ x_4 l_{42} + x_5 l_{52} + x_6 l_{62} \\
f_3= cQ+ x_4 l_{43} + x_5 l_{53} + x_6 l_{63}, 
\endalign$$
where $Q:=Q(x_0,\dots x_3)$ is a defining equation of $\o{F}$
and $l_{ij}$ is a linear form. 
 
By Claim 2 (2),
$P\cap P'$ is a line or a point.
 
First assume that $P\cap P'$ is a point.
Then we may assume that $P' =(x_0 = x_1 = x_2 = x_6 =0)$
and $P\cap P'$ is $(0:0:0:1:0:0:0)$.
Furthermore
$$\align
f_1= aQ+a'R+x_4 l_{41} + x_5 l_{51} + x_6 l_{61} \\
f_2= bQ+b'R+x_4 l_{42} + x_5 l_{52} + x_6 l_{62} \\
f_3= cQ+c'R+x_4 l_{43} + x_5 l_{53} + x_6 l_{63}, 
\endalign$$
where $l_{4i}$ and $l_{5i}$ are linear forms of $x_0$, $x_1$ and $x_2$.
By calculating the Jacobian matrix, we see that $P\cap P'$ is 
a singular point, a contradiction.
\endcomment

\definition{N=8}
We deny the case that $N=8$. 
Assume that $N=8$. By Theorem 6.0, we may assume that
any index $2$ point is a $\frac 12 (1,1,1)$-singularity
or a QODP.
In this case $Y=Y'$ holds since $e=0$.
By $F(X')=1$ and the $\Bbb Q$-factoriality of $X'$, 
there exists a line $l$ intersecting $C$.
Let $l'$ be the strict transform of $l$ on $Y$.
By $-K_{Y}.l' = -K_{X'}.l - E'.l'$ and the fact that
$-K_{Y}$ is nef,
we have $-K_{Y}.l'= 0 \teb{and} E'.l'=1$ or
$-K_{Y}.l'= \frac 12 \teb{and} E'.l'=\frac 12$.
But the latter case does not occur since $e=0$.
In the former case $E \cap l' = \phi$ by $E.l'=0$. 
Hence $K_X.f(l') = \frac 12$,
which in turn show that 
for a $\Bbb Q$-Fano blow up whose center is an index $2$ point on $f(l')$, 
the resulting weak $\Bbb Q$-Fano $3$-fold 
is not a $\Bbb Q$-Fano $3$-fold. But by the tables in Section 5,
we again fall into table 2 for the new choice of a $\Bbb Q$-Fano blow up, 
a contradiction (the new $e$ must be $0$).
\enddefinition

\comment
\definition{$N=8$}
Finally we deny the case that $N=8$. In this case $Y=Y'$ holds.
We first show that the Fano index of $X'$ must be $2$.

Assume that the Fano index of $X'$ is $1$.
Then there exists a line $l$ intersecting $C$.
By Claim 1, $-K_Y.l'= \frac 12 \teb{or} 0$,
where $l'$ is the strict transform of $l$.
But the latter case does not occur since $e=0$.
In the former case $E \cap l' = \phi.$ 
Hence $K_X.f(l') = \frac 12$,
which in turn show that if we blow up $X$ at a $\frac 12(1,1,1)$-singularity on
$f(l')$, the resulting $3$-fold is a weak $\Bbb Q$-Fano $3$-fold
which is not a $\Bbb Q$-Fano $3$-fold. But it is impossible by the Table 2.
($e$ must be $0$.)

\enddefinition
\endcomment

\definition{Table 3}
If $h=4$ and $N=1$,
then $(X, -K_X) \simeq (((5)\subset \Bbb P (1^4, 2)), \Cal O(1))$ is
an example.

If $h=4$ and $N=2$,
then $(X, -K_X) \simeq (((3,4)\subset \Bbb P (1^4,2^2)), \Cal O(1))$ is an 
example.

\comment
If $h=5$ and $N=4$, then for a flopped curve $l$,
we have $E' .l = -\frac 32 \t{E} .l = \frac 32$ since
by Lemma 5.9, $g (E) \simeq E$, where $g$ is the flopping contraction from $Y$.
But this contradicts the fact that $E'$ is a Cartier divisor.
\endcomment
\enddefinition

\definition{Table 4}
\enddefinition
\definition{$h=6$ and $N=1$}
About the case that $h=6$ and $N=1$, we know that 
an example exists by Corollary 8.1 below and
the existence of the case that $h=6$ and $N=2$ (see Table 1).
\enddefinition

\definition{$h=6$ and $N=8$}
We will show that if $h=6$ and $N=8$, then $F(X)=1$.
So we will exclude this case.

In this case, $f'$ is a $\Bbb P^1$-bundle associated to some vector bundle
$\Cal E$ of rank $2$ on $\Bbb P^2$. 
Let $T$ be its tautological divisor.
By the adjunction formula
$-K_{Y'}\sim 2T- (c_1(\Cal E)-3)L$,
we have $6=(-K_{Y'})^3 = 8T^3 - 6{c_1(\Cal E)}^2 + 54$
and hence $c_1(\Cal E)$ is an even.
Hence $H' := 3T -(\frac 32 c_1(\Cal E) - 4)L$ is an integral Cartier
divisor. Note that $H'\equiv -K_{Y'} + \frac 12 \t{E}$.
Hence for a flipped curve ${l_i}^+$ on some $Y_i$
and the strict transform $H_i$ of $H'$ on $Y_i$,
we have $H_i . {l_i}^+ = -2$.
Hence the strict transform $H$ of $H'$ on $Y$ is a Cartier divisor
numerically equivalent to $-K_Y + \frac 12 E$.
Note that $H$ is $f$-numerically trivial.
So by [KMM, Lemma 3-2-3 (2)], $f(H)$ is a Cartier divisor
and clearly numerically equivalent to $-K_X$.
\enddefinition
\comment 
We construct an example for the case that $h=6$ and $N=8$.
Let $\t{E}$ be a smooth del Pezzo surface of degree $2$.
We represent $\t{E}$ as a surface which is obtained by $7$ points blow up
from $\Bbb P^2$ and let $E_i$ be the exceptional curves ($i= 1, \dots 7$).
and $l$ the total transform of a line in $\Bbb P^2$.
On the other hand, $\t{E}$ is realized by a double cover of $\Bbb P^2$
whose branch locus is a smooth quartic. Let $p: \t{E} \to \Bbb P^2$
be such a double cover and 
$\Cal E := p_* (\Cal O_{\t{E}} (2l + 2 K_{\t{E}}))$.
Since $p$ is flat, $\Cal E$ is a locally free sheaf of rank $2$.
Set $Y':= \Bbb P (\Cal E)$ and 
$f'$ the natural projection $\Bbb P(\Cal E) \to \Bbb P^2$.
Then $\t{E}$ is naturally embedded in $Y'$
and $H|_{\t{E}} = 2l + 2 K_{\t{E}}$, where $H$ is the tautological line bundle
for $\Cal E$. Let $L:= {f'}^* \Cal O(1)$.

We show that $\t{E}=2H+2L$ and $-K_{Y'} = 2H+3L$.
We recall that $-K_{Y'}= 2H -(c_1(\Cal E) -3)L$.
Let $\t{E} = 2H + aL$. Compute $H.L.\t{E}$ in two ways. Then
$$2=(2l +2K_{\t{E}}).(-K_{\t{E}})=H.L.\t{E}= 2c_1(\Cal E) +a. \tag a$$
On the other hand, by the adjunction formula, we have
$-K_{\t{E}}= -K_{Y'} -\t{E} |_{\t{E}} = -(a+c_1(\Cal E) -3)L|_{\t{E}}$.
Since $-K_{\t{E}} = L|_{\t{E}}$, we obtain 
$$a+c_1(\Cal E)-3 =-1. \tag b$$
By (a) and (b), we obtain $c_1(\Cal E)=0$ and $a=2$.
So we are done.

Next we show that all $E_i$'s belong to the same numerical class
and the numerical class generate an extremal ray of $Y'$.
Note that $L.E_i = -K_{\t{E}}.E_i = 1$.
On the other hand $H.E_i = (2l+2 K_{\t{E}}).E_i = -2$.
Since $\rho (Y') =2$, it show that all $E_i$'s belong to the same
numerical class.
Let $m$ be the numerical class and $n$ the numerical class of a fiber of $f'$.
Let $C$ be an irreducible curve and $[C]=am+bn$ for some $a,b \in \Bbb Q$.
By $L.C\geq 0$, we have $a\geq 0$. We want to show $b\geq 0$.
Note that $\t{E}.m = (2H+2L).m=-2$ and $\t{E}.n=2$.
If $\t{E}.C =-2a+2b \geq 0$, then $b\geq a \geq 0$ and we are done.
Assume that $\t{E}.C <0$.
Then $C\subset \t{E}$ and hence we can write as
 $C\sim \alpha l - \sum \alpha _i E_i$ in $\t{E}$.
Clearly we have $\alpha \geq 0$.
So
$$
2b=(-2K_{Y'} -\t{E}).C=4l.(\alpha l - \sum \alpha _i E_i)=
4\alpha \geq 0.$$
Hence $b\geq 0$ which in turn show $m$ generates an extremal ray.

Let $R_1$ be the extremal ray generated by $m$.
$-2K_{Y'} -\t{E}$ is a supporting divisor for $R_1$.
We will show that $|-2K_{Y'} -\t{E}|$ is free.
Since $-2K_{Y'} -\t{E}= \t{E}+2L$,
$\text{Bs} |-2K_{Y'} -\t{E}| \subset \t{E}$.
Consider the exact sequence
$$0 \to \Cal O_{Y'}(-2K_{Y'}-2\t{E})\to \Cal O_{Y'}
(-2K_{Y'} -\t{E})\to \Cal O_{\t{E}}
(-2K_{Y'} -\t{E}|_{\t{E}}) \to 0. $$  
Since $(-3K_{Y'} -2 \t{E}).m =1$ and $(-3K_{Y'} -2 \t{E}).n = 2$,
$-3K_{Y'} -2\t{E}$ is ample. Hence by
KKV vanishing theorem and the above exact
sequence, $\text{Bs} |-2K_{Y'} -\t{E}|= \text{Bs} |-2K_{Y'} -\t{E}|_{\t{E}}|$.
But $-2K_{Y'} -\t{E}|_{\t{E}} = 4l$ is free and hence
$|-2K_{Y'} -\t{E}|$ is also free.

We can see that $\supp R_1 = \cup E_i$.
In fact if $[C] \in R_1$ for an irreducible curve $C$,
$C\subset \t{E}$ by $\t{E}.C<0$.
Write $C\sim \alpha l - \sum \alpha _i E_i$.
By $(-2K_{Y'} -\t{E}).C=0$, we have $\alpha =0$.
If $C\not =E_i$, then $\alpha _i \geq 0$, a contradiction.
So $C= E_i$ for some $i$. By this $R_1$ is a flipped ray.

Observe that $\Cal N _{E_i / Y'} \simeq \Cal O(-1)\oplus \Cal O(-2)$.
Hence for the inverse of the flip 
(we will denote it by $Y' \dashrightarrow Y$),
$Y$ has seven $\frac 12 (1,1,1)$-singularities on the flipping curve
${E_i}'$. We can see that $Y$ is a $\Bbb Q$-Fano $3$-fold with $\rho (Y)=2$
and the strict transform $E$ of $\t{E}$ is $\Bbb P^2$ and can be
contracted to a $\frac 12(1,1,1)$-singularity.
Let $X$ be the target of the contraction.
Then $X$ is what we want.
We need to check that the Fano index of $X$ is $\frac 12$.
If the Fano index of $X$ is $1$,
by Sano's classification, $\rho(X) =2$ since
$(-K_X)^3 = 10$, a contradiction. Hence we are done.
\endcomment
\comment
\definition{$h=10$ and $N=1$}
Next we prove that the case that $h=10$ and $N=1$ does not occur.
Assume that this case occurs.
Since $\text{deg} \Delta =0$ and $\t{E}.m =1$ for a fiber $m$ of $f'$,
we have $Y' = \Bbb P ({f'}_* \Cal O_{Y'} (\t{E}))$.
We examine $\Cal E := {f'}_* \Cal O_{Y'} (\t{E})$.

We see first that $\t{E}$ is a del Pezzo surface of degree $3$.
In fact, for a flopped curve $l$,
$-K_{\t{E}}.l= \t{E}.l >0$ and for a fiber $m$ of $f'$,
$(-K_{Y'} -\t{E}).m= 1 >0$, which in turn show that $-K_{\t{E}}$ is ample.

Furthermore $\t{E}$ is normal.
Assume that $\t{E}$ is non normal. 
Then by [RM5, Theorem 1.1],
we have $3$ possibilities for the normalization ${\t{E}}^n$ of $\t{E}$:

\roster
\item ${\t{E}}^n \simeq \Bbb F_{3,0}$.
\item ${\t{E}}^n \simeq \Bbb F_1$ and the inverse image of the non normal locus
is a smooth member in $|C_0 + m|$.
\item ${\t{E}}^n \simeq \Bbb F_1$ and 
the inverse image of the non normal locus is a reducible member in $|C_0 + m|$.
\endroster
Note that there is a birational morphism
$\t{E} \to E\simeq \Bbb P^2$. But for the $3$ possibilities, it is impossible,
a contradiction.
Hence $\t{E}$ is normal.

Let $\mu:{\t{E}}^{\mu} \to \t{E}$ be the minimal resolution,
$p:=f' |_{\t{E}}$ and $p' := p \circ \mu$.
Let $D:= \t{E} |_{\t{E}}$ and $\o{D}:= p_* D$.
By the negativity lemma, we can write as 
$\mu ^* D = {p'}^* \o{D} -G$, where $G$ is an effective exceptional curve
for $p'$.
Hence we have an injection 
$0\to p_* \Cal O_{\t{E}} (D) \to \Cal O_{\Bbb P^2}(\o{D})$
and hence $p_* \Cal O_{\t{E}} (D)$ is an invertible sheaf on $X'$.
By Riemann-Roch Theorem and KKV vanishing theorem,
$h^0 (\Cal O_{\t{E}} (D)) = 0$.
So the degree of $p_* \Cal O_{\t{E}} (D)$ is negative.
Consider the exact sequence
$$0\to \Cal O_{X'} \to \Cal E \to p_* \Cal O_{\t{E}} (D) \to 0.$$
This sequence is split.
Hence $\Cal E = \Cal O_{\Bbb P^2} \oplus \Cal O_{\Bbb P^2} (-a)$
for some positive integer $a$.
By the definition of $\Cal E$, $\t{E}$ is the tautological divisor
for $\Cal E$.
Since for any line $l$ in $\Bbb P^2$, 
$\t{E}|_{{f'}^{-1} l}$ is the negative section of ${f'}^{-1} l$,
$\t{E} \simeq \Bbb P^2$, a contradiction.
Hence we exclude this case.
\enddefinition

\definition{Table 5}
\enddefinition
\definition{$h=5$ and $N=1$}
If $h=5$ and $N=1$,
then $(X, -K_X) \simeq (((3,3)\subset \Bbb P (1^5, 2)), \Cal O(1))$ is
an example. (See also Theorem 8.2 for this two case.)
Furthermore by [T2, p.9 (22)],
$Y'$ is embedded in 
$\Bbb P(\Cal O_{\Bbb P_1} ^{\oplus 3} \oplus \Cal O_{\Bbb P_1} (1))$
as a divisor linearly equivalent to $3H+F$,
where $H$ is the tautological divisor and $F$ is a fiber. 
\enddefition
\comment
If $h=5$ and $N=2$, we start from a smooth subvariety $Y'$ in
$\Bbb P(\Cal O_{\Bbb P_1} \oplus 
{\Cal O_{\Bbb P_1} (1)} ^{\oplus 3} \oplus \Cal O_{\Bbb P_1} (2))$
which is a complete intersection of two divisors $V_i \sim 2H-F$ ($i=1,2$),
where $H$ is the tautological divisor and $F$ is a fiber. 
\endcomment

\comment
\definition{$h=5$ and $N=2$}
If $h=5$ and $N=2$, we construct an example.

Consider a weighted projective space
$\Bbb P(1^5, 2)$.
Let $x_i$ ($i=0,1,2,3,4$) be the coordinates of weight $1$
and $y$ be the coordinate of weight $2$.
Take a weighted complete intersection
$(((3,3)\subset \Bbb P (1^5, 2))$ which is defined by the following
$2$ equations:

$$f=l(x)y + q_3 (x) x_3 + q_4 (x) x_4 $$ and 
$$g=L(x)y + Q_3 (x) x_3 + Q_4 (x) x_4 ,$$
where 
$l$ and $L$ are linear forms, $q_i$ and $Q_i$ are polynomials of degree $2$.
Note that $Z:= (f=0) \cap (g=0)$ contains $\o{E} :=(y = x_3 = x_4 =0)$
which is $\Bbb P^2$.
By taking $l$, $L$, $q_i$ and $Q_i$ generally,
we assume that $Z$ is smooth outside $\o{E}$ and their $\frac 12 (1,1,1)$-
singularities.
Then we can see that 
$Z$ has exactly $8$ ordinary double points on $\o{E}$ by easy computations.
Let $Y:= \pmb{Proj}  \bigoplus _{i=0}^{\infty} \Cal O_Z (m\o{E})$.
Then $Y$ is smooth outside one $\frac 12 (1,1,1)$-singularity and 
the inverse image $E$ of $\o{E}$ remains $\Bbb P^2$.
We can contract $E$ to a $\frac 12(1,1,1)$-singularity
and let $X$ be the target. Then $X$ is an example as we want.
($\rho (X) = 1$????)
\enddefinition
\endcomment

Now we completes the proof of the main theorem.

We close this section after proving some corollaries to the main theorem.

\proclaim{Corollary 6.1}
Let $X$ be a $\Bbb Q$-factorial
$\Bbb Q$-Fano $3$-fold with the following properties:
\roster 
\item $\rho(X)=1$;
\item $I(X)=2$;
\item $F(X)=\frac 12$;
\item $h^0 (-K_X) \geq 4$.
\endroster
Then $(-K_X)^3$ and $\text{aw}(X)$ are effectively bounded
as in the main theorem.
\endproclaim

\demo{Proof}
By the main theorem and Theorem 6.0, we obtain the assertion
since $(-K)^3$ and $\text{aw}$ are invariant under a deformation. 
\qed
\enddemo

\proclaim{Lemma 6.2}
Let $X$ be a $\Bbb Q$-factorial $\Bbb Q$-Fano $3$-fold with $\rho (X)=1$,
$I(X)=2$ and $F(X)=\frac 12$
and $f:Y \to X$ a weak Fano blow up with $I(Y)=2$
such that $f$-exceptional divisor $E$ is contracted to a point by $f$.
If $|-2K_Y|$ is free and $h^0 (-K_Y-E) >0$,
then $H^0(\Cal O_Y(-2K_Y)) \to H^0(\Cal O_E(-2K_Y|_E))$ is surjective. 
\endproclaim

\demo{Proof}
We are inspired by [RM1, p.29, Step 4].
It suffices to prove that \newline
$h^1 (\Cal O_Y (-2K_Y -E))=0$. 
Take a general member $T\in |-2K_Y|$.
Then by the exact sequence
$$0 \to \Cal O_Y(-E)\to \Cal O_Y(-2K_Y -E)\to \Cal O_T(-2K_Y -E|_T)
\to 0$$
and $h^i(\Cal O_Y(-E))=0$ for $i=1,2$
(these vanishing easily follow from
$$0\to \Cal O_Y(-E) \to \Cal O_Y \to \Cal O_E \to 0$$
since by Proposition 2.3, $h^1 (\Cal O_E)=0$),
we obtain $h^1 (\Cal O_Y (-2K_Y -E)) =
h^1 (\Cal O_T (-2K_Y -E|_T))$.
By Serre duality, we have
$h^1 (\Cal O_T (-2K_Y -E|_T))=h^1 (\Cal O_T (2K_T -E|_T))
=h^1 (\Cal O_T (K_Y +E|_T))$.
We prove that $h^1 (\Cal O_T (K_Y +E|_T))=0$ in the below.
Take a member $F\in |-K_Y -E| \not =\phi$.
Then since $\rho (X) =1$ and $-K_X$ is 
a positive generator of $Z^1 (X)/\equiv$,
we can write $F= F' + rE$,
where $F'$ is a prime divisor and $r$ is a nonnegative integer.
Since $|-2K_Y|$ is free and $T$ is general,
we may assume that $F'|_T$ and $E|_T$ is irreducible.
Note that $(F'+rE)|_T .E|_T = (-K_Y -E). E. (-2K_Y)>0$ and $(E|_T)^2 <0$.
Hence if $r>0$,
for every integer $b\in [1,r]$,
we have
$(F'|_T + (r-b) E|_T).bE|_T >0$,
which means $F|_T$ is numerically $1$-connected.
So by [RC, Lemma 3],
we have $H^0 (\Cal O_{F|_T}) \simeq \Bbb C$.
Hence by the exact sequence  
$$0\to \Cal O_T(-F|_T) \to \Cal O_T \to \Cal O_{F|_T} \to 0,$$
we have $h^1(\Cal O_T(-F|_T)) = 0$
which is exactly what we want.
\qed
\enddemo

\proclaim{Corollary 6.3}
Let $X$ be a $\Bbb Q$-factorial $\Bbb Q$-Fano $3$-fold 
with the following properties:
\roster 
\item $\rho(X)=1$;
\item $I(X)=2$;
\item $F(X)=\frac 12$;
\item $h^0 (-K_X) \geq 4$.
\endroster
Then for any index $2$ point $P$,
there exists a smooth rational curve $l$ through $P$ such that
$-K_X.l=\frac 12$.
\endproclaim

\demo{Proof}
First we treat the case that any index $2$ point is of type as in
Theorem 5.0 (5).
By Table 1 $\sim$ Table 5 and nonexistence of the case that
$h=4$ and $N=8$ and the case that $h=6$ and $N=8$,
$e$ is positive or $f'$ is a crepant divisorial contraction 
for any choice of an index $2$ point $P$.
Let $g:Y\to Z$ be the anti-canonical model.
Let $l'$ be a flopping curve if $g$ is a flopping contraction
or a general fiber of $E'$ if $g$ is a crepant divisorial contraction. 
Then by Lemma 6.2, $g (E)\simeq E$ whence $E. l' = 1$.
Hence $l:=f(l')$ is what we want.

Next we treat the general case. 
Let $\frak f: \frak X \to \Delta$ be a flat family as in Theorem 6.0.
By [KoMo, Corollary 12.3.4], $\rho (\frak X _t) =1$
and furthermore by Table 1 $\sim$ 5, [San1] and [San2],
$\frak X _t (t\not = 0)$ satisfies the assumptions of Theorem 5.0.
Let $P$ be an index $2$ point on $X$ and 
$P_t$ an index $2$ point on $\frak X _t$ which specializes to $P$.
By the first part of this proof,
there is a curve $l_t$ on $\frak X _t (t\not =0)$ such that
$l_t \simeq \Bbb P^1$, $P_t \in l_t$ and $-K_{\frak X_t}.l_t =\frac 12$.
Since there are only countably many components of relative Hilbert scheme
$\text{Hilb} (\frak X / \Delta)$,
we may assume that they form a flat family over $\Delta$.
Furthermore by the properness of a component of relative Hilbert scheme,
this family extends over $0$. Let $l$ be its fiber over $0$. 
Then $l$ is what we want.
\qed
\enddemo

\comment
\proclaim{Corollary 6.3}
Let $X$ be a $\Bbb Q$-Fano $3$-fold with the assumptions of the main theorem.
Then for any $\frac 12(1,1,1)$-singularity $P$??,
there exists a smooth rational curve $l$ through $P$ such that
$-K_X.l=\frac 12$.
\endproclaim

\demo{Proof}
By Table 1 $\sim$ Table 5 and nonexistence of the case that
$h=4$ and $N=8$ and the case that $h=6$ and $N=8$,
we have $e>0$ for $X$ and any $P$.
Let $l'$ be a flopping curve. 
Then by Lemma 6.2, $E. l' = 1$.
Hence $l:=f(l')$ is what we want.
\qed
\enddemo
\endcomment

\comment
The next lemma is used in Section 6.
 
\proclaim{Proposition 5.11}  
Let $X$ be a weak $\Bbb Q$-Fano $3$-fold with canonical singularities.
Assume the following:
\roster
\item
$X$ has only $\frac 12 (1,1,1)$-singularities as its non Gorenstein points;
\item
$|-K_X|$ has a member with only canonical singularities
(we will call such a member a canonical elephant after Miles Reid);
\item
for any $\frac 12 (1,1,1)$-singularity, there is no curve $l$ which passes
through it such that $-K_X .l =0$;
\item
$h^0(\Cal O_X(-K_X)) \geq 2$ and $(-K_X)^3 \geq 1$.
\endroster

Let $P$ be a $\frac 12 (1,1,1)$-singularity
and $f: Y \to X$ the blow up of $P$.
Then $Y$ is a weak $\Bbb Q$-Fano $3$-fold.
\endproclaim

\demo{Proof}
We have $$K_Y= f^*K_X + \frac 12 E, \tag 1$$
where $E$ is the exceptional divisor of $f$.
Let $D$ be a canonical elephant of $|-K_X|$.
Then $D$ has an ordinary double point at $P$.
Let $\t{D}$ be the strict transform of $D$ on $Y$.
Write $$f^*D \equiv  \t{D} + kE. \tag 2$$
Then we have $K_{\t{D}} = (f|_{\t{D}})^*K_D + (\frac 12 - k) E$
by (1) and (2).
But since $f|_{\t{D}}$ is the blow up of an ordinary double point,
we have $k= \frac 12$. In particular we see that $\t{D} \in |-K_Y|$.
To check the nefness of $-K_Y$, we have only to compute intersection
numbers of $-K_Y$ and curves in $\text{Bs} |-K_Y|$.
Let $\t{l}$ be an irreducible curve in $\text{Bs} |-K_Y|$.
If $\t{l} \subset E$, then $-K_Y.\t{l} >0$. Hence we may assume that
$\t{l} \not \subset E$. Denote the curve $f(\t{l})$ by $l$.
We clearly assume that $\t{l} \cap E \not = \phi$.
Since $f(\text{Bs} |-K_Y|) \subset \text{Bs}|-K_X|$, 
$l \subset \text{Bs} |-K_X|$.
By the exact sequence 
$$0\to \Cal O_X \to \Cal O_X(-K_X) \to \Cal O_D(-K_X|_{D}) \to 0$$
and $h^1(\Cal O_X) = 0$, 
we have $h^0 (\Cal O_D(-K_X|_{D})) = h^0(\Cal O_X(-K_X))-1$ and 
hence $C:=-K_X|_{D}$ is effective by the assumption (4).
We also know that
$\text{Bs} |-K_X||_{D} = \text{Bs} |C|$.
Since $D$ is a K3 surface and $l\subset \text Bs |C|$, 
$l\simeq \Bbb P^1$ by [Al, Corollary 1.5].
By the same reason, $\t{l}\simeq \Bbb P^1$.
Hence we have $$\t{l}=(f|_{\t{D}})^*l -\frac 12 F, \tag 3$$
where $F$ is the exceptional curve of $f|_{\t{D}}$ and 
$F= E|_{\t{D}}$. 
By (1) and (3), we obtain 
$-K_Y.\t{l}= -K_X.l - \frac12$.
On the other hand $-K_X.l \geq \frac 12$ by the assumption (3) and 
the assumption that $X$ has only singularities of index $2$.
Hence $-K_Y.\t{l} \geq 0$ and the nefness of $-K_X$ is proved.

By (1), $(-K_Y)^3 = (-K_X)^3 - \frac 12$.
Hence by the assumption that 
$(-K_X)^3 \geq 1$ and $-K_Y$ is nef,
we see that $-K_Y$ is big.
\qed
\enddemo
\endcomment

\head 7. Existence of an anti-canonical divisor 
with only canonical singularities
for a $\Bbb Q$-Fano $3$-fold with $I(X)=2$
\endhead

\proclaim{Proposition 7.0}
Let $X$ be as in Theorem 5.0 and assume furthermore that
$X$ has only $\frac 12 (1,1,1)$-singularities as its non Gorenstein points.
Then 
$|-K_Y|$ has no base curve containing a $\frac 12 (1,1,1)$-singularity.
The similar assertion holds also for $X$. 
\endproclaim

\demo{Proof}
Consider the exact sequence
$$\ex{Y}{-K_Y -E}{Y}{-K_Y}{E}{-K_Y |_E} .$$
Assume that we prove $$h^0 (\Cal O_Y(-K_Y -E)) = h-3. \tag {7.1.1}$$
Then the map $H^0 (\Cal O_Y (-K_Y)) \to H^0 (\Cal O_E (-K_Y |_E))$
is surjective. Since $|-K_Y |_E|$ is free, we know that
$|-K_Y|$ has no base curve intersecting $E$.
By this, we know that $|-K_X|$ has no base curve through $f(E)$.
Since $f(E)$ is any $\frac 12 (1,1,1)$-singularity,
it means $|-K_X|$ has no base curve containing a $\frac 12(1,1,1)$-singularity.
By this, $|-K_Y|$ has also
no base curve containing a $\frac 12(1,1,1)$-singularity.

So it suffices to show that (7.1.1) holds.
We note that (7.1.1) is equivalent to
$$h^0 (\Cal O_{Y'}(-K_{Y'} -\t{E})) = h-3.\tag {7.1.2}$$
We will prove this using the data of the tables in Theorem 5.0.
\enddemo
\definition{Table 1}
We have $-K_{Y'} -\t{E} \sim {f'}^* D$,
where $D$ is a primitive ample Weil divisor (we can easily see that the linear
equivalent class of $D$ is unique).
Hence $h^0(-K_{Y'} -\t{E})=h^0 (D)$. $h^0 (D)=h-3$ is easy to see. 
\enddefinition

\definition{Table 2 or Table 3}
We have $-K_{Y'} -\t{E} \sim E'$ whence
$h^0(-K_{Y'} -\t{E})=1 =h-3$.
\enddefinition

\definition{Table 4}
Since $-K_{Y'} -\t{E}-K_{Y'}$ is nef and big, we can compute
$h^0(-K_{Y'} -\t{E})$ by Riemann-Roth theorem and we are done.

But if $h=6$, then $L \sim -K_{Y'} -\t{E}$ and hence
we can see that $h^0(-K_{Y'} -\t{E})=h^0 (L)=3=h-3$ more easily.
\enddefinition

\definition{Table 5}
Since $-K_{Y'} -\t{E}-K_{Y'}$ is nef and big, we can compute
$h^0(-K_{Y'} -\t{E})$ by Riemann-Roth theorem and we are done.
But 
if $h=5$, then $L \sim -K_{Y'} -\t{E}$ and hence
we can see that $h^0(-K_{Y'} -\t{E})=h^0 (L)=2=h-3$ more easily.
\enddefinition
\qed

\proclaim{Proposition 7.1}
Let $X$ be a weak $\Bbb Q$-Fano $3$-fold with log terminal singularities
and satisfies the following condition:
\roster
\item $|-K_X|\not = \phi$;
\item there are a finite number of non Gorenstein points on $X$;
\item there is a member of $|-K_X|$ which is normal near
non Gorenstein points.
\endroster
Then $|-K_X|$ has a member which is normal and has only canonical singularities
outside non Gorenstein points of $X$.
\endproclaim

\demo{Proof}
The proof is almost the same as one of [Am, Main Theorem].
So we only give an outline of the proof.
Let $U:= \{ x| x \teb{is a Gorenstein point of} X \}$.
Let $S$ be a general member of $|-K_X|$.
Let $\gamma := \max \{ t | K_X + tS|_U \teb{is log canonical} \}$.
As Ambro did, it suffices to prove that there is no element of 
$\text{CLC} (K_X + \gamma S |_U)$ contained in $\text{Bs} |-K_X|$. 
Assume the contrary and let $Z$ be a minimal element of
$\text{CLC} (K_X + \gamma S |_U)$ contained in $\text{Bs} |-K_X|$.
By the assumption (3), $Z$ is a complete variety.
Hence by using Theorem 1.0, we know that it suffices to prove
$H^0 (\Cal O _Z (-K_X |_Z)) \not = 0$.
It is done by Adjunction Theorem and a nonvanishing argument.
\qed
\enddemo
 
\proclaim{Corollary 7.2}
Let $X$ be as in Theorem 5.0 and assume furthermore that
$X$ has only $\frac 12 (1,1,1)$-singularities as its non Gorenstein points.
Then $|-K_X|$ has a member with only canonical singularities.
\endproclaim

\demo{Proof}
Fix a $\frac 12(1,1,1)$-singularity $P$ and the blow up $f:Y\to X$ at $P$.
By Proposition 7.0 and Proposition 7.1,
we can find a member $S\in |-K_Y|$ such that $S$ is normal and
has only canonical singularities outside $\frac 12(1,1,1)$-singularities
of $Y$. Since $f|_S$ is crepant, $f(S)$ has only canonical singularity
outside $\frac 12(1,1,1)$-singularities of $X$ except $P$. 
Since $P$ is any $\frac 12(1,1,1)$-singularity,
we can find a member of $|-K_X|$ with only canonical singularities.
\qed
\enddemo

\comment
\definition{Conjecture 0.2 (Miles Reid)}
Let $X$ be a $\Bbb Q$-Fano $3$-fold with only terminal singularities.
Then a general member of $|-K_X|$ has only canonical singularities.
(M.Reid calls such a member canonical elephant. So this conjecture
is quoted in the literatures as General Elephant Conjecture.)
\enddefinition

\proclaim{Main Theorem}
Let $X$ be a $\Bbb Q$-Fano $3$-fold with only canonical singularities
and the following properties:
\roster 
\item $I(X)=2$;
\item $\text{Bs} |-K_X|$ contains neither 
$1$-dimensional components of $\text{Sing} (X)$
or non pseudo-terminal singular points 
(in particular any index $2$ point is pseudo-terminal);
\item $|-K_X|$ has a member which is irreducible and reduced;
\item $h^0 (-K_X) \geq 4$.
\endroster

Then $|-K_X|$ has a member with only canonical singularities.
\endproclaim

We construct a birational morphism $f:Y \to X$,
where $Y$ is a weak $\Bbb Q$-Fano $3$-fold
with only canonical singularities and $I(Y) \leq 2$.
Let $g:Y \to Z$ be the anti-canonical model.
We prove that 
$Z$ satisfies the same properties as $X$
with $(-K_Z)^3 =(-K_X)^3 - \frac 12$.
Hence the theorem follows from the existence of a canonical elephant for
a Gorenstein Fano $3$-fold with only canonical singularities
by induction on the value of $(-K)^3$.

\proclaim{Corollary}
Let $X$ be a $\Bbb Q$-Fano $3$-fold with only terminal singularities
and assume that $\rho (X)=1$ and $h^0 (-K_X) \geq 4$.
Then $|-K_X|$ has a member with only canonical singularities.
\endproclaim

\demo{Proof}
By the main theorem, we only have to check that $|-K_X|$ has a member
which is irreducible and reduced. But this follows from [Al, Theorem 2.18]
by $rho (X) =1$.
\qed
\enddemo

\definition{Acknowledgement}
I express my hearty thanks to Professor Shigefumi Mori for
giving me useful comments. 
Lemma 1.1 is inspired by a conversation with him.
I am grateful to Professor Takayuki Hayakawa 
who kindly sent me his preprints [Hay1] and [Hay2].
\enddefinition

\definition{Notation and Conventions}
\roster
\item In this paper, we will work over $\Bbb C$, the complex number field;
\item we denote the linear equivalence by $\sim$
and the numerical equivalence by $\equiv$.
The equality $=$ in an adjunction formula means the $\Bbb Q$-linear
equivalence;
\item we denote the Hirzebruch surface of degree $n$ by $\Bbb F_n$
and the surface which is obtained by the contraction of the negative section
of $\Bbb F_n$ by $\Bbb F_{n,0}$.
\endroster
\enddefinition

\head 1. A partial resolution for a $3$-dimensional pseudo-terminal singularity
\endhead
\definition{Definition 1.0}
Let $(X, p)$ be a germ of a $3$-dimensional singularity.
We say that $(X, p)$ is pseudo-terminal if it is canonical and
its canonical cover has only cDV singularities.
\enddefinition
Pseudo-terminal singularities are classified by [HT] and [KS]
(see also [S]).

The following lemma can be checked by reading [KY2], [Hay1] and [Hay2] but
we contain the proof for reader's convenience.

\proclaim{Lemma 1.1}
Let $(X, p)$ be a germ of a $3$-dimensional pseudo-terminal singularity.
Then there is a birational morphism $f: Y \to X$ with the following properties:
\roster
\item $Y$ has only canonical singularities with indices $\leq 2 ? r?$;
\item $-K_Y$ is $f$-nef;
\item $\text{Bs} |-K_Y|$ contains neither 
$1$-dimensional components of $\text{Sing} (Y)$
or non pseudo-terminal singular points 
(in particular any index $2$ point of $Y$ is pseudo-terminal).
\endroster

In fact we construct $f$ as a composition of divisorial extractions
with minimal discrepancies at pseudo-terminal points with index $> 2$.
\endproclaim

\demo{Proof}
In the following, we construct a sequence of weighted blow ups:
$$Y:=Y_k \overset f_k\to\longrightarrow Y_{k-1}, \dots ,
Y_1 \overset f_1\to\longrightarrow X,$$
where $f_i :Y_i \to Y_{i-1} (i\geq 2)$ (resp. 
$f_1 :Y_1 \to X$) is a weighted blow up at the unique
non Gorenstein point $p_{i-1}$ of index $>2$ on $Y_{i-1}$
(resp. $p_0:=p$).
The uniqueness of $p_i$ follows from the below construction. 
We denote by $E_i$ the $f_i$-exceptional divisor. 
In the below, we choose the blow ups
so that $E_i$'s are reduced (but possibly reducible).  
The check of (1) can be done as [Hay1] and [Hay2].
So we do not mention this below. 
The check of (2) is done as follows:

Let $C_i$ be a curve in $|-K_{Y_i} | _{E_i}|$ and ${C_i}'$(resp. ${E_i}'$)
be the strict transform of $C_i$(resp. $E_i$) on $Y_j$ for any $j \geq i+1$
(for notational simplicity, we use such notation independent of $j$). 
We see that $|-K_{Y_i} | _{E_i}| \not = \phi$ and
$C_i$ can be chosen so that ${C_i}_{\text{red}}$ is irreducible
and $p_{i+1} \in {C_i}'$. 
Assume that $-K_k$ is nef over $X$ for $k\leq i$.
We claim here that to check the nefness of $-K_{Y_{i+1}}$ over $X$, 
we have only to show that $-K_{Y_{i+1}}.{C_i}' \geq 0$.
Let $l$ be an irreducible curve in ${E_k}'$ for $k\leq i$ such that
$l \not \subset {E_j}'$ for $k < j \leq i+1$.
If $p_i \not \in l$, then $-K_{Y_{i+1}}.l \geq 0$ by inductive assumption.
If $p_i \in l$, 
then ${C_k}'.l \leq -K_{Y_{i+1}}.l$ since 
$-K_{Y_{i+1}} |_{{E_k}'}\equiv \sum a_j {E_j}'|_{{E_k}'}+{C_k}'$ 
for some $a_j \geq 0$. Hence if $l \not = {C_k}'_{\text{red}}$,
$-K_{Y_{i+1}}.l \geq 0$ by ${C_k}'.l\geq 0$. 
For $k\leq i-1$, we have $-K_{Y_{i+1}}.{C_k}' \geq 0$ by inductive assumption
and $p_k \in {C_k}'$. Hence we only have to check $-K_{Y_{i+1}}.{C_i}' \geq 0$.

\roster
\item "(cA/2)" 
$$(X,p)\simeq (\{xy +f(z^2,u)=0\} / \Bbb Z _2 (1,1,1,0), o).$$
Let $f_1$ be the weighted blow up with weight $\frac 12(2k-1,1,2,1)$,
where $k:=\text{ord}f(Z, U)$. 
Then $E_1 \simeq (\{xy +f_k(z^2,u)=0\} \subset \Bbb P (2k-1,1,2,1))$ and
$p_1$ is 
the cyclic quotient singularity of $\frac 1{2k-1} (-2,2,1)$-type.
Let $C_1$ be a general curve in $|\Cal O_{E_1} (1)|$.
Let $f_i$ ($2\leq i \leq k$) be the weighted blow up with weight 
$\frac 1{2(k-i)+3} (2(k-i)+1,2,1)$ at $p_{i-1}$ 
which is the cyclic quotient singularity of $\frac 1{2(k-i)+3} (-2,2,1)$ type.
Let $C_i$ ($2 \leq i \leq k$) be the unique curve in $|\Cal O_{E_i} (1)|$,
where $E_i \simeq \Bbb P (2(k-i)+1,2,1)$.
Then we can easily check that 
$-K_i.{C_{i-1}}' =0$ and $p_i \not \in {C_{i-1}}'$.
We also see that $Y_i$ has only pseudo-terminal singularities
and $E_i \cap \teb{Sing} Y_i$ is $0$-dimensional.
Hence (3) holds.

\item "(cAx/2)"
$$(X,p)\simeq (\{x^2 +y^2 +f(z,u)=0\} / \Bbb Z _2 (0,1,1,1), o).$$
Let $\text{ord}f =2k$. We only treat the case that $k$ is even since
if $k$ is odd, the below argument holds by exchanging $x$ and $y$.
Let $f_1$ be the weighted blow up with weight $\frac 12(k,k+1,1,1)$.
Then $E_1 \simeq (\{x^2 +f_{\text{deg} =2k}(z,u)=0\} 
\subset \Bbb P (k,k+1,1,1))$ and
$p_1$ is 
the cyclic quotient singularity of $\frac 1{k+1} (-1,1,1)$-type.
Note that $E_1$ is possibly reducible.
Let $C_1$ be a general curve in $|\Cal O_{E_1} (1)|$.
Let $f_i$ ($2\leq i \leq k$) be the weighted blow up with weight 
$\frac 1{k+3-i} (k+2-i,1,1)$ at $p_{i-1}$ 
which is the cyclic quotient singularity of $\frac 1{k+3-i} (-1,1,1)$ type.
Let $C_i$ ($2 \leq i \leq k$) be a general curve in $|\Cal O_{E_i} (1)|$,
where $E_i \simeq \Bbb P (k+2-i)+1,1,1)$.
Then we can easily check that 
$-K_i.{C_{i-1}}'=0$ and $p_i \not \in {C_{i-1}}'$.
We also see that $Y_i$ has only pseudo-terminal singularities
and $E_i \cap \teb{Sing} Y_i$ is $0$-dimensional.
Hence (3) holds.

To treat the cases (cD/2) and (cE/2), we need to examine the cases
(cAx/4) and (cD/3) and prove the following:

\item "(cAx/4)" 
$$(X,p)\simeq (\{x^2+y^2 +f(z,u)=0\} / \Bbb Z _4 (1,3,1,2), o).$$
Give $z$ and $u$ weights as wt $z=\frac 14$ and wt $u=\frac 12$.
Let wt $f =\frac {2k+1}2$. We only treat the case that $k$ is even since
if $k$ is odd, the below argument holds by exchanging $x$ and $y$.
Let $f_1$ be the weighted blow up with weight $\frac 14(2k+1,2k+3,1,2)$.
Then $E_1 \simeq (\{x^2 +f_{\text{wt} ={\frac {2k+1}2}}(z,u)=0\} 
\subset \Bbb P (2k+1,2k+3,1,2))$ and
$p_1$ is 
the cyclic quotient singularity of $\frac 1{2k+3} (2k+1,1,2)$-type.
Note that $E_1$ is possibly reducible.
Let $C_1$ be the unique curve in $|\Cal O_{E_1} (1)|$.
Let $f_i$ ($2\leq i \leq k+2$) be the weighted blow up with weight 
$\frac 1{2(k-i)+7} (2(k-i)+5,1,2)$ at $p_{i-1}$ 
which is the cyclic quotient singularity of $\frac 1{2(k-i)+7} (-2,1,2)$ type.
Let $C_i$ ($2 \leq i \leq k+2$) be the unique curve in $|\Cal O_{E_i} (1)|$,
where $E_i \simeq \Bbb P (2(k-i)+5,1,2)$.
Then we can easily check that 
$-K_i.{C_{i-1}}'=0$ and $p_i \not \in {C_{i-1}}'$.
We also see that $Y_i$ has only pseudo-terminal singularities
and $E_i \cap \teb{Sing} Y_i$ is $0$-dimensional.
Hence (3) holds.

\item "(cD/3-1) (resp. (cD/3-2))"
$$(X,p)\simeq (\{u^2 +x^3+yz(y+z)=0\} / \Bbb Z _3 (2,1,1,0), o).$$
(resp. 
$$(X,p)\simeq (\{u^2 +x^3+yz^2+xy^4\lambda (y^3) +y^6 \mu (y^3)=0\} / 
\Bbb Z _3 (2,1,1,0), o).)$$
Let $f_1$ be the weighted blow up with weight $\frac 13(2,4,1,3)$.
Then $E_1 \simeq (\{ u^2+x^3+yz^2 =0\} 
\subset \Bbb P (2,4,1,3))$ and
$$(Y_1,p_1) \simeq 
(\{{u'}^2 +{x'}^3+{z'}({y'}+{z'})=0\} / \Bbb Z _4 (2,1,1,3), o).$$ 
(resp.
$$(Y_1,p_1) \simeq 
(\{{u'}^2 +{x'}^3+{z'}^2+{x'}{y'}^4\lambda ({y'}^4) +{y'}^6 \mu ({y'}^4)\} 
/ \Bbb Z _4 (2,1,1,3), o).)$$ 
Let $C_1$ be the unique curve in $|\Cal O_{E_1} (1)|$.
Let $f_2$ be the weighted blow up with weight 
$\frac 14 (2,1,1,3)$ (resp. $\frac 14 (2,1,5,3)$) at $p_1$.
Note that these weighted blow ups coincide 
one in the treatment of case (cDx/4) above. 
Then we can easily check that 
$-K_2.{C_1}'=0$ and $P_2 \not \in {C_1}'$.
The definition of $C_i$ and 
the check $-K_{Y_{i+1}}.{C_i}' \geq 0$ for $i\geq 2$ are done 
in the treatment of the case (cDx/4). So we omit them.
We also see that $Y_i$ has only pseudo-terminal singularities
and $E_i \cap \teb{Sing} Y_i$ is $0$-dimensional.
Hence (3) holds.

\item "(cD/3-3)"
$$(X,p)\simeq 
(\{u^2 +x^3+y^3+xyz^3\alpha (z^3) +xz^4 \beta (z^3) +yz^5 \gamma(z^3)
+z^6 \delta (z^3)=0\} / \Bbb Z _3 (2,1,1,0), o).$$
Let $f_1$ be the weighted blow up with weight $\frac 13(2,4,1,3)$.
Then $E_1 \simeq (\{ u^2+x^3+xz^4 \beta (0) +z^6 \delta (0) =0\} 
\subset \Bbb P (2,4,1,3))$ and
$$(Y_1,p_1) \simeq 
(\{{u'}^2 +{x'}^3+{y'}^2+{x'}{y'}{z'}^3\alpha ({y'}{z'}^3) 
+{x'}{z'}^4 \beta ({y'}{z'}^3) +{y'}{z'}^5 \gamma({y'}{z'}^3)
+{z'}^6 \delta ({y'}{z'}^3)=0 \} / \Bbb Z _4 (2,1,1,3), o).$$ 
Let $C_1$ be the unique curve in $|\Cal O_{E_1} (1)|$.
We see that $Y_1$ has only pseudo-terminal singularities and
outside $C_1$ (and hence outside $\text{Bs} |-K_{Y_1}|$), 
$E_1 \cap \teb{Sing} Y_1$ is $0$-dimensional.
Let $f_2$ be the weighted blow up with weight 
$\frac 14 (2,1,5,3)$ at $p_1$.
Note that these weighted blow ups coincide 
one in the treatment of case (cDx/4) above. 
Then we can easily check that 
$-K_2.{C_1}'=0$ and $P_2 \not \in {C_1}'$.
The definition of $C_i$ and 
the check $-K_{Y_{i+1}}.{C_i}' \geq 0$ for $i\geq 2$ are done 
in the treatment of the case (cDx/4). So we omit them.
We also see that $Y_i$ has only pseudo-terminal singularities
and outside $\text{Bs} |-K_{Y_i}|$, 
$E_i \cap \teb{Sing} Y_i$ is $0$-dimensional.
Hence (3) holds.

\item "(cE/2)"
$$(X,p)\simeq (\{u^2+x^3 +\phi (x,y,z)=0\} / \Bbb Z _2 (0,1,1,1), o).$$
We give $x,y,z$ and $u$ weights as wt $(x,y,z,u) =\frac 12 (2,3,1,3)$.
We divide (cE/2) into finer cases as follows:
\item "(cE/2-1)"
$$\phi _{\text{wt} \leq 3} = a_{0,4} z^4 x + b_{0,6} z^6 +yz^3.$$
\item "(cE/2-2)"
$$\phi _{\text{wt} \leq 4} = a_{0,4} z^4 x + b_{0,6} z^6 +
yz^3 (a_{1,3} x + b_{1,5} z^2)+z^6 (a_{0,6} x+ b_{0,8} z^2) +y^2 z^2.$$
\item "(cE/2-3)"
$$\phi _{\text{wt} \leq 4} = a_{0,4} z^4 x + b_{0,6} z^6 +
yz^3 (a_{1,3} x + b_{1,5} z^2)+z^6 (a_{0,6} x+ b_{0,8} z^2)$$
and $y^3 z$ or $y^4 \in \phi$.

Let $f_1$ be the weighted blow up with weight $\frac 12(2,3,1,3)$.
Then $E_1 \simeq (\{u^2+x^3+ a_{0,4} z^4 x + b_{0,6} z^6 +yz^3 \} 
\subset \Bbb P (2,3,1,3))$ in (cE/2-1) case and 
$E_1 \simeq (\{u^2+x^3+ a_{0,4} z^4 x + b_{0,6} z^6 \} 
\subset \Bbb P (2,3,1,3))$ in (cE/2-2) case or (cE/2-3) case.

Furthermore
$$(Y_1,p_1) \simeq 
(\{{u'}^2 +{x'}^3+a_{0,4} {z'}^4 x + b_{0,6} {z'}^6 +{z'}^3 +{y'} \psi\} / 
\Bbb Z _3 (2,1,1,0), o),$$
where $\psi \in \Bbb C \{x',y',z'\}$ in (cE/2-1) case,
$$(Y_1,p_1) \simeq 
(\{ {u'}^2 +{x'}^3+a_{0,4} {z'}^4 {x'} + b_{0,6} {z'}^6 +
{y'}{z'}^3 (a_{1,3}{x'} + 
b_{1,5} {z'}^2)+{y'}{z'}^6 (a_{0,6}{x'}+ b_{0,8} {z'}^2) +{y'} {z'}^2+
{y'}^2 \psi \} 
/ \Bbb Z _4 (2,1,1,3), o),$$ 
where $\psi \in \Bbb C \{x',y',z'\}$ in (cE/2-2) case,
and
$$(Y_1,p_1) \simeq (\{{u'}^2 +{x'}^3+a_{0,4} {z'}^4 {x'} + b_{0,6} {z'}^6 +
{y'}{z'}^3 (a_{1,3}{x'} + b_{1,5}{z'}^2)
+{y'}{z'}^6 (a_{0,6}{x'}+ b_{0,8}{z'}^2)+ {y'}^2 \psi \} 
/ \Bbb Z _4 (2,1,1,3), o),$$
where $\psi \in \Bbb C \{x',y',z'\}$ and ${z'}$ or ${y'} \in \psi$ in (cE/2-3).

Let $C_1$ be the unique curve in $|\Cal O_{E_1} (1)|$.
From now on we treat three cases separately.
\item "(cE/2-1)"
We see that $Y_1$ has only pseudo-terminal singularities
and $E_1 \cap \teb{Sing} Y_1$ is $0$-dimensional.
Let $f_2$ be the weighted blow up at $p_1$ with weight $\frac 12(2,1,1,3)$ if
$x' \in \psi$ (resp. weight $\frac 12(2,1,4,3)$ if $x' \not \in \psi$. 
We can easily check that $-K_2.{C_1}'=0$ and $p_2 \not \in {C_1}'$.
Hence (1)$\sim$(3) are checked by the result in case (cA/3) or (cD/3).
\item "(cE/2-2)"
We see that $Y_1$ has only pseudo-terminal singularities
and $E_1 \cap \teb{Sing} Y_1$ is $0$-dimensional.
Let $f_2$ be the weighted blow up at $p_1$ with weight $\frac 12(2,4,1,3)$. 
$-K_2.{C_1}'=0$ and $p_2 \not \in {C_1}'$.
Hence (1)$\sim$(3) are checked by the result in case (cD/3).
\item "(cE/2-3)"
We see that $Y_1$ has only canonical singularities,
and outside $C_1$ (and hence outside $\text{Bs} |-K_{Y_1}|$), 
$Y_1$ has only pseudo-terminal singularities and
$E_1 \cap \teb{Sing} Y_1$ is $0$-dimensional.
Let $f_2$ be the weighted blow up at $p_1$ with weight $\frac 12(2,1,4,3)$
if ${z'}\in \psi$(resp. weight $\frac 12(2,4,1,3)$ if ${y'} \in \psi$).
We can easily see that 
$-K_2.{C_1}'=0$ and $p_2 \not \in {C_1}'$.
Hence (1)$\sim$(3) are checked by the result in case (cD/3).

\item "(cD/2)"
$$(X,p)\simeq (\{\phi(x,y,z,u)=0\} / \Bbb Z _2 (1,0,1,1), o),$$
where $\phi$ has one of the following forms:
\item "(cD/2-1)"
$$\phi = x^2 +yzu + y^a +\lambda z^{2b} +\mu u^{2c},$$
where $\lambda, \mu \in \Bbb C, a\geq 3, b,c\geq 2$. 

We have only to treat the cases 
(1) $b=c=2$ and (2) $b\geq 3$ or $\lambda=0$.
Assume that $b=c=2$ (resp. $b\geq 3$ or $lambda=0$).
Let $f_1$ be the weighted blow up with weight $\frac 12(3,2,1,1)$ 
(resp. weight $\frac 12(3,2,1,3)$ if $b\geq 3$).
Then $E_1 \simeq (\{yzu+\lambda z^4 +\mu u^4 =0\} 
\subset \Bbb P (3,2,1,1))$
(resp. $E_1 \simeq (\{x^2+\delta_{a,3}y^3 +\delta_{b,3} \lambda z^6 =0\} 
\subset \Bbb P (3,2,1,3))$).
$p_1$ is 
the cyclic quotient singularity of $\frac 13 (2,1,1)$-type
(resp. 
$$(Y_1,p_1) \simeq (\{{x'}^2 +{u'}^{a-3}{y'}^a+{y'}{z'} + 
\lambda {u'}^{b-3}{z'}^{2b} +\mu {u'}^{3c-3} \} 
/ \Bbb Z _3 (0,2,1,1), o).)$$
We see that $Y_1$ has only pseudo-terminal singularities
and $E_1 \cap \teb{Sing} Y_1$ is $0$-dimensional.
Let $C_1$ be a general curve (resp. the unique curve) in $|\Cal O_{E_1} (1)|$.
Let $f_2$ be the weighted blow up at $p_1$ with weight $\frac 13 (2,1,1)$
(resp. $\frac 13 (3,2,1,1)$ if $c=2$, 
or weight $\frac 13 (3,2,4,1)$ if $c\geq 3$ or $\mu=0$).
Then we can easily check that 
$-K_2.{C_1}'=0$ and $p_2 \not \in {C_1}'$.

\item "(cD/2-2)"
$$\phi = x^2 +y z^2 + \lambda zu^{2a-1} + g(y,u^2),$$
where 
$\lambda \in \Bbb C, a\geq 1, g(y,u^2) 
\in (y^3, y^2 u^2, u^4) \Bbb C \{y,u^2\}$ and $h(y,0)\not \equiv 0$.
Assume that $a=2$ (resp. $a\geq 3$ or $\lambda=0$).
Let $f_1$ be the weighted blow up with weight $\frac 12(3,2,1,1)$ 
(resp. weight $\frac 12(3,4,1,1)$).
Then 
$E_1 \simeq (\{yz^2+\lambda zu^3 +\mu u^4=0\} 
\subset \Bbb P (3,2,1,1))$ 
(resp. $E_1 \simeq (\{x^2+yz^2+\delta_{a,3}\lambda zu^5 +\nu u^6=0\} 
\subset \Bbb P (3,4,1,1))$)
(resp. $E_1 \simeq (\{x^2+\delta_{a,3}y^3 +\delta_{b,3} \lambda z^6 =0\} 
\subset \Bbb P (3,2,1,3))$).
$p_1$ is 
the cyclic quotient singularity of $\frac 13 (2,1,1)$-type
(resp. 
$$(Y_1,p_1) \simeq (\{{x'}^2 +{z'}^2 + 
\lambda {y'}^{a-3}{z'}{u'}^{2a-1} +\nu {u'}^6 +
{y'}(b{u'}^4 +c{u'}^8) +{y'}^2 \psi \} 
/ \Bbb Z _4 (3,2,1,1), o).)$$
We see that $Y_1$ has only pseudo-terminal singularities.
Let $C_1$ be a general curve in $|\Cal O_{E_1} (1)|$.
Let $f_2$ be the weighted blow up at $p_1$ with weight $\frac 13 (2,1,1)$
(resp. $\frac 14 (3,2,1,1)$ if $c=2$, 
or weight $\frac 13 (3,2,4,1)$ if $c\geq 3$ or $\mu=0$).
Then we can easily check that 
$-K_2.{C_1}'=0$ and $p_2 \not \in {C_1}'$.
\endroster
\qed
\enddemo

\head 2. the proof of the main theorem
\endhead

\demo{Proof of the main theorem}
For a weak $\Bbb Q$-Fano $3$-fold $*$ with only canonical singularities,
we consider the following properties:
\roster 
\item"(1$_*$)" $I(*)=2$;
\item"(2$_*$)" $\text{Bs} |-K_*|$ contains neither 
$1$-dimensional components of $\text{Sing} (*)$
or non pseudo-terminal singular points;
\item"(3$_*$)" $|-K_*|$ has a member which is irreducible and reduced;
\item"(4$_*$)" $h^0 (-K_*) \geq 4$.
\endroster

Let $p$ be an index $2$ point of $X$ and $f: Y \to X$ be a birational
morphism localized to the morphism as in Lemma 1.1.
Then $f$ is projective by [Nak, Corollary 1.6].
By $(4_X)$, Riemann-Roch theorem [RM2, Theorem (10.2)] and 
$-K_X .c_2 (X) \leq 24$ [KY, Lemma 2.2 and Lemma 2.3], 
we know that $(-K_X)^3 >2$. 
Hence $|-2K_X|$ is free by Theorem ??.
Since $|-2K_X|$ is ample and free,
$\text{Bs} |-2K_X - p|$ is a finite set of points.
So by $H^0 (-2K_Y) \simeq H^0 (-2K_X -p)$ and Lemma 1.1 (2), we know 
$-K_Y$ is nef. 
By Lemma 1.1, $(-K_Y)^3 = (-K_X)^3 -\frac 12 >0$.
Hence $Y$ is a weak $\Bbb Q$-Fano $3$-fold with only canonical singularities.
We will check that $(2_Y)$ and $(3_Y)$.
By $(2_X)$ and Lemma 1.1 (3), $(2_Y)$ is satisfied.
Since $h^0 (-K_Y)=h^0 (-K_X) \geq 4$, a general member of $|-K_Y|$
contains no irreducible component of $\text{excep} f$.
Hence $(3_Y)$ is satisfied.  

Let $g: Y \to Z$ be the anti-canonical model. We will show that $Z$
satisfies the properties as in the main theorem.
The only unclear thing is to check $(2_Z)$.
Assume that 
$\text{Bs} |-K_Z|$ contains a
$1$-dimensional component $C$ of $\text{Sing} (Z)$
(resp. a non pseudo-terminal singular points $Q$).
Then by $-K_Y =f^* (-K_Z)$ and $(3_Y)$,
there is no $g$-exceptional divisor $F$ such that $g(F) \subset C$
(resp. $g(F) =Q$).
Hence the transform $C'$ of $C$ on $Y$ is contained in
$\text{Sing} (Y) \cap \text{Bs} |-K_Y|$
(resp. $g^{-1} (Q) \subset \text{Bs} |-K_Y|$, and
there is a non pseudo-terminal singular point on $g^{-1} (Q)$
or 
there is a component of $g^{-1} (Q)$ contained in $\text{Sing} Y$), 
a contradiction to $(2_Y)$.

Note that $(-K_Z)^3 = (-K_X)^3 -\frac 12$. If $I(Z)=2$, then we repeat
the similar procedure. After a finite number of such procedures, we obtain
a Gorenstein Fano $3$-fold $W$ with only canonical singularities.
By [RM1, Theorem (0.5)], $|-K_W|$ has a member with only canonical
singularities.
By reversing the procedure, this member is transformed to a member
of $|-K_X|$ with only canonical singularities. So we are done.
\qed
\enddemo
\endcomment

\head 8. Some properties of $\Bbb Q$-Fano $3$-folds 
with only $\frac 12(1,1,1)$-singularities as its non Gorenstein points
\endhead

\proclaim{Theorem 8.0}
Let $X$ be a (not necessarily $\Bbb Q$-factorial)
weak $\Bbb Q$-Fano $3$-fold with $I(X)=2$.
Assume that $X$ has the following properties:

\roster
\item $|-K_X|$ has no fixed component. $|-K_X|$ has no
base curve containing an index $2$ point;
\item $|-K_X|$ has a member with only canonical singularities;
\item there is no divisor contracted to a point by the morphism
defined by $|-mK_X|$ for $m>>0$;
\item $h^0(-K_X)\geq 4$;
\item all non Gorenstein singularities of $X$  are
$\frac 12 (1,1,1)$-singularities.
\endroster

Then $X$ can be deformed to a weak $\Bbb Q$-Fano $3$-fold
with only $\frac 12 (1,1,1)$-singularities as its singularities.
\endproclaim

\demo{Proof}
Let $N$ be the number of $\frac 12(1,1,1)$-singularities.
We will prove this theorem by induction of $N$.
We treat the case that $N=0$ later.
First we prove that if the assertion holds in case of $N-1$,
the assertion holds also in case of $N$
(hence we assume that $N>0$).

By assumption (4), $h^0(-K_X) \geq 4$ and
hence by Riemann-Roch theorem and KKV vanishing theorem,
we have $(-K_X)^3 \geq 2$.
Let $P$ be any $\frac 12(1,1,1)$-singularity
and $f:Y\to X$ the blow up at $P$. Then by Proposition 4.1,
$Y$ is a weak $\Bbb Q$-Fano $3$-fold.

Then we verify that the assumption (1)$\sim$ (5) hold for $Y$.
(5) is clear.
Since $f^{-1} |-K_X| = |-K_Y|$.
(1), (2) and (4) follows.
For (3), we assume that there is a divisor $F$ on $Y$ which
is contracted to a point by the morphism
defined by $|-mK_Y|$.
If $E\cap F \not = \phi$,
then $E\cap F$ is a curve since $E$ is a Cartier divisor.
By the nature of $F$, $E\cap F$ is numerically trivial for $-K_Y$
but by the nature of $E$, $E\cap F$ is numerically negative for $-K_Y$,
a contradiction.
Hence $E\cap F = \phi$.
Then however $f(F)$  is contracted to a point by the morphism
defined by $|-mK_X|$, a contradiction.
Hence $Y$ satisfies (1)$\sim$(5).
By the assumption of the induction,
the Gorenstein points of $Y$ are smoothable.
Let $\Cal Y \to \Delta$ be 
a $1$-parameter smoothing of Gorenstein points of $Y$.
Then by [KoMo, Proposition 11.4],
we obtain the deformation $\Cal X \to \Delta$ of $X$
which satisfies the commutative diagram
 $$\matrix
{\Cal Y} & \to & {\Cal X} \\
{\searrow} & \ &  {\swarrow } \\
 \ & \Delta & \ & .
\endmatrix $$
Then 
$\Cal Y _t \to \Cal X _t$ is an $E_5$ type contraction for $t\in \Delta$
since a contraction of type $E_5$ is stable under a deformation by [Kod2].
Hence $X_t$ is a smoothing of Gorenstein points.

Next we prove the assertion in case $N=0$.
The proof is the same as one of [Na] except the following claim:

\proclaim{Claim}
Let $D$ be a member of $|-K_X|$ with only canonical singularities. 
Then $\text{Pic} X \to \text{Pic} D$ is an injection.
\endproclaim

\demo{Proof}
The proof is similar to one of [Na, Proposition 2] by virtue of the assumption
(3).
 
\qed
\enddemo

\qed
\enddemo

\proclaim{Corollary 8.1}
Let $X$, $Y$, $h$ and $N$ as in the main theorem.
Let $g: Y\to Z$ be the anti-canonical model.
Assume that $X$ has only $\frac 12 (1,1,1)$-singularities as its
non Gorenstein points.
Then if $N>1$ (resp. $N=1$), 
$Z$ can be deformed to a $\Bbb Q$-Fano $3$-fold $Z'$ 
with $\rho (Z')=1$ and $F(Z')=\frac 12$
which has only $N-1$ $\frac 12 (1,1,1)$-singularities as its singularities
and $h^0 (-K_{Z'}) = h$ (resp. a smooth Fano $3$-fold $Z'$ 
with $\rho (Z')=1$, $F(Z')=1$ and $h^0 (-K_{Z'}) = h$.)
\endproclaim

\demo{Proof}
To apply Theorem 8.0, we only have to check
that $Z$ satisfies the assumption (1) and (5).
By Proposition 7.0, any flopping curve does not contain a 
$\frac 12 (1,1,1)$-singularity. So the assumption (1) holds.
Hence the assumption (5) is also satisfied by Proposition 7.0.

By these, we can apply Theorem 8.0 and $Z$ can be deformed
to a $\Bbb Q$-Fano $3$-fold $Z'$ with only $\frac 12 (1,1,1)$-singularities.
We can prove by the proof of [KoMo, Corollary 12.3.4] that
$\rho(Z')=1$. If $N>1$, $F(Z')=\frac 12$ by [San1] and [San2].
If $N=1$, we have clearly $F(Z')=1$. Hence we are done.  
\qed
\enddemo

The next result is a first step for the classification of Mukai's type
[Mu3, Theorem 1.10].

\proclaim{Theorem 8.2 (Embedding into a Weighted Projective Space)}
Let $X$ be a (not necessarily $\Bbb Q$-factorial)
$\Bbb Q$-Fano $3$-fold with canonical singularities and $I(X)=2$.
Assume that $X$ has the following properties:
\roster 
\item $|-K_X|$ is indecomposable, i.e.,
$|-K_X|$ contains no member which is a sum of two movable Weil divisors;
\item $|-K_X|$ has no base curve containing an index $2$ point;
\item $|-K_X|$ has a member with only canonical singularities;
\item for any index $2$ point,
there is a smooth rational curve $l$ through it such that
$-K_X.l=\frac 12$;
\item $h^0(X, \Cal O(-K_X)) \geq 4$;
\item all non Gorenstein singularities of $X$
are $\frac 12 (1,1,1)$-singularities.
\endroster

Then except the following two cases (a) and (b), 
$X$ is embedded into $\Bbb P (1^h, 2^N)$
and $-K_X$ is the restriction of $\Cal O(1)$,
where $h:=h^0 (-K_X)$ and 
$N$ is the number of $\frac 12 (1,1,1)$-singularities:
$$\Phi _{|-K_X|} \teb{is a double cover of} \Bbb P^3  \teb{branched along
a sextic}. \tag a$$
$$\Phi _{|-K_X|}
 \teb{is a double cover of a quadric hypersurface branched along
the intersection with a quartic}. \tag b$$
(Note that in case (a),
$$X\simeq ((6)\subset \Bbb P (1^4,3)).$$
Note also that in case (b),
$$X\simeq ((2,4)\subset \Bbb P (1^5,2))$$
but 
the number of weight $2$ is not equal to the number of non Gorenstein point.)

Furthermore $X$ is an intersection of weighted hypersurfaces of degree
$\leq 6$ in $\Bbb P(1^h, 2^N)$.

If $h=4$ and $N=1$,
then $X \simeq ((5)\subset \Bbb P (1^4,2))$.

If $h=4$ and $N=2$,
then $X \simeq ((3,4)\subset \Bbb P (1^4,2^2))$.

If $h=5$ and $N=1$,
then $X \simeq ((3,3)\subset \Bbb P (1^5, 2))$.
\endproclaim

\demo{Proof}
We prove this by induction of $N$.

In case $N=0$, the assertion follows from
[Mu3, Theorem 6.5 and Proposition 7.8].  

Next we prove that if the assertion holds in case
$X$ has $N-1$ $\frac 12 (1,1,1)$-singularities,
then so does it in case
$X$ has $N$ $\frac 12 (1,1,1)$-singularities.
Let $X$ be a $\Bbb Q$-Fano $3$-fold satisfying the assumptions of this theorem
and with $N$ $\frac 12 (1,1,1)$-singularities.
Let $f:Y\to X$ be the blow up at a $\frac 12 (1,1,1)$-singularity.
Let $E$ be the exceptional divisor of $f$.
Then $Y$ is a weak $\Bbb Q$-Fano $3$-fold by Proposition 4.1.
By the assumption (5), $Y$ is not $\Bbb Q$-Fano $3$-fold.
Let $g: Y \to Z$ be the morphism defined by a sufficient multiple
of $-K_X$ and $\o{E}:=g(E)$.

\proclaim{Claim 1}
$Z$ satisfies the assumption of this theorem
and has $N-1$ $\frac 12 (1,1,1)$-singularities.
\endproclaim

\demo{Proof}
By $-K_Y = g^* (-K_Z)$, if $|-K_Z|$ is decomposable, $|-K_X|$ must be
decomposable, a contradiction. Hence (1) is satisfied.

By (2) for $X$, neither $|-K_Y|$ has 
a base curve containing an index $2$ point.
Hence 
any $g$-exceptional curve does not contain an index $2$ point.
So by $-K_Y = g^* (-K_Z)$,
(2) is satisfied and (6) is also satisfied.

Let $D$ be a member of $|-K_X|$ with only canonical singularities.
Then the strict transform $D'$ of $D$ on $Y$ has the same property
since $D' \to D$ is crepant. Since $D' \to g(D')$ is crepant,
$g(D')$ has also the same property. Hence (3) is satisfied.

By $-K_Y = g^* (-K_Z)$ and $h^0(-K_Y) = h^0 (-K_X)$, 
we know that (5) is satisfied.

We show last that (4) is satisfied.
If $Z$ is Gorenstein, there is nothing to prove.
If $Z$ is non Gorenstein, let $Q$ be any $\frac 12 (1,1,1)$-singularity
and we will denote the corresponding points on $Y$ and $X$ also by $Q$.
Then by (4) for $X$, there is a curve $Q\in l$ on $X$ as stated in (4).
For the strict transform $l'$ of $l$ on $Y$,
we have $-K_Y.l'= \frac 12 \teb{or} 0$. But if the latter case occurs,
$l'$ is a base curve of $|-K_Y|$ 
containing an index $2$ point $Q$, a contradiction.
Hence $-K_Y.l' = \frac 12$ and then $-K_Z. g(l') =\frac 12$.
By blowing up $Q$, $Z$ becomes a weak $\Bbb Q$-Fano $3$-fold
by Proposition 4.1 and (5) for $Z$.
Then $g(l')$ become a curve contracted by an multi-anti-canonical morphism.
Hence $g(l')$ is a smooth rational curve. 
Now we complete the proof of the claim.
\qed
\enddemo

Hence by the assumption of the induction, the following three cases occur:
\definition{Case $\alpha$}
$Z\subset \Bbb P(1^h, 2^{N-1})$ and $-K_Z = \Cal O_Z (1)$;
\enddefinition
\definition{Case $\beta$}
$Z$ is of type (a);
\enddefinition
\definition{Case $\gamma$}
$Z$ is of type (b).
\enddefinition

\proclaim{Claim 2}
$\text{Bs}|-K_X|$ coincides with $\frac 12 (1,1,1)$-singularities as a
set.
\endproclaim

\demo{Proof}
If $N=0$, the assertion follows from 
[Mu3, Theorem 6.5 and Proposition 7.8].
Hence by Claim 1, the assertion follows by induction with respect to
the number of $\frac 12 (1,1,1)$-singularities.
\qed
\enddemo  
 
\definition{Case $\alpha$}
We first show that $\o{E} \simeq E$.
By Claim 2, the similar assertion holds for $|-K_Y|$.
Hence $H^0 (\Cal O_Y (-K_Y) ) \to H^0 (\Cal O_E (-K_Y)) \simeq 
H^0 (\Cal O_{\Bbb P^2}(1))$ is surjective.
Hence  $H^0 (\Cal O_Y (-mK_Y) ) \to H^0 (\Cal O_E (-mK_Y))$ is also
surjective for all $m\geq 0$ 
since $\oplus _{m\geq 0} H^0 (\Cal O_{\Bbb P^2}(m))$ is simply generated.
So $\o{E} \simeq E$ since $g$ is defined by $-mK_Y$ for some $m>0$.

We note here that there is an elementary transformation
$\Bbb P(1^h, 2^N) \dashrightarrow \Bbb P(1^h, 2^{N-1})$
which is decomposed as follows:

Let $\Bbb P$ be the projective bundle over 
$\Bbb P(1^h, 2^{N-1})$ whose vector bundle is $\Cal O \oplus \Cal O(-2)$
and $T$ the effective tautological divisor (which is unique).
Let $a$ be the contraction morphism of $T$.
Then $a(\Bbb P)$ is isomorphic to $\Bbb P(1^h, 2^N)$. 
Let $b: \Bbb P \to \Bbb P(1^h, 2^{N-1})$ be the natural projection.
Then our elementary transformation is $b\circ a^{-1}$.

We seek a natural morphism $Y \to \Bbb P$.
For this, we prove that there is a natural surjection
$g^*(\Cal O_Z \oplus \Cal O_Z(-2)) \to \Cal O_Y(E)$.

There is the natural injection $\Cal O_Y(-E) \to \Cal O_Y$
which represents $\Cal O_Y(-E)$ as the ideal sheaf of $E$.
By Theorem 4.0, there is a member $S \in |-2K_X|$ such that
$f^* S \cap E = \phi$. Associated to $S$, there is an injection
$\Cal O_Y (-f^* S) \to \Cal O_Y$. 
This gives an injection
$\Cal O_Y (-E) \to g^* \Cal O_Z (2)$
since $g*\Cal O_Z (2) \simeq \Cal O_Y (-2K_Y)$,
$-f^*(-2K_X) \sim -(-2K_Y) -E$. 
By these, we can define a injection $\Cal O_Y(-E) \to
g^*(\Cal O_Z \oplus \Cal O_Z(2))$. 
Since $f^* S \cap E = \phi$, the cokernel of this map is locally free
and hence the dual of this map is a surjection.
\comment
Let $d:g^*(\Cal O_Z \oplus \Cal O_Z(-2)) \to \Cal O_Y(E)$ be its dual,
$\Cal L := \text{Im}  d$ and $\Cal C':= \text{coker}  d$.
Since $\text{rank} \Cal C =1$, $\text{rank} \Cal L =1$.
Hence $\Cal C'$ is a torsion sheaf.
Denote the composition
$$\Cal O_Y \hookrightarrow g^*(\Cal O_Z \oplus \Cal O_Z(-2)) \to \Cal O_Y(E)
$$ by $d'$.
We can easily see that $d'$ is not a $0$ map.
By this and $h^0(\Cal O_Y(E))=1$, we have 
$H^0(\Cal L) \simeq H^0 (\Cal O_Y(E))$.
Since $\Cal C'$ is a torsion, it means $\Cal L \simeq \Cal O_Y(E)$.
This is what we want.
\endcomment
Let $\iota: Y \to \Bbb P$ be the morphism over $\Bbb P(1^h, 2^{N-1})$
associated to the surjection
$d:g^*(\Cal O_Z \oplus \Cal O_Z(-2)) \to \Cal O_Y(E)$.
By the definition of $\iota$, we have $\iota(E)= T|_{\iota(Y)}$.
In particular $\iota(E)$ is Cartier on $\iota(Y)$.
Since $\o{E}\simeq \Bbb P^2$, $\iota(E)$ is also $\Bbb P^2$ by the Zariski's
Main Theorem. Hence $\iota(Y)$ is smooth at points of $\iota(E)$.
\enddefinition

\proclaim{Claim 3}
$\iota(Y)$ is normal.
\endproclaim

\demo{Proof}
It suffices to prove that $\iota_* \Cal O_Y = \Cal O_{\iota(Y)}$.
The natural
morphism $\Cal O_{\iota(Y)} \to \iota_* \Cal O_Y$ is injective since
the kernel is at most torsion sheaf.
Let $\Cal C$ be its cokernel.
\comment
By the factorization 
$Y\overset{\iota} \to \rightarrow \iota(Y)\overset{p} \to \rightarrow Z$,
we have a natural morphism  
$\Cal O_Z \to p_* \Cal O_{\iota(Y)} \to g_* \Cal O_Y$.
In this, $\Cal O_Z \to p_* \Cal O_{\iota(Y)}$ and 
$\Cal O_Z \to g_* \Cal O_Y$ is isomorphisms and hence
$p_* \Cal O_{\iota(Y)} \to g_* \Cal O_Y$ is also an isomorphism.
Hence $p_* \Cal K=0$.
\endcomment 
We will prove that $p_* \Cal C=0$.
By the exact sequence
$$0\to \Cal O_{\iota(Y)} \to \iota _* \Cal O_Y \to \Cal C \to 0 ,$$ we have
$$0\to p_*\Cal O_{\iota(Y)} \to p_*\iota _* \Cal O_Y \to p_*\Cal C \to 
R^1 p_* \Cal O_{\iota(Y)}. $$
Since $p_*\Cal O_{\iota(Y)} \to p_*\iota _* \Cal O_Y$ is an isomorphism,
it suffices to prove that $R^1 p_* \Cal O_{\iota(Y)}=0$.
Consider the exact sequence
$$0\to \Cal I_{\iota(Y)} \to \Cal O_{\Bbb P} \to \Cal O_{\iota(Y)} \to 0 .$$  
Since the dimension of a fiber of $p$ $\leq 1$, 
we have $R^2 p_* \Cal I_{\iota(Y)} =0$.
Since $\Bbb P$ is a $\Bbb P^1$-bundle, we have 
$R^1 p_* \Cal O_{\Bbb P} =0$. Thus we obtain $R^1 p_* \Cal O_{\iota(Y)}=0$
and we are done.

Since every fiber of $g: Y \to Z$ intersects
$\iota(E)$ and $\iota(Y)$ is smooth at points of $\iota(E)$,
any fiber is not contained in the singular locus of $\iota(Y)$.
Let $l$ be any $1$-dimensional fiber of $g$.
By the theorem on formal functions,
we have $\Cal C \otimes \Cal O_l =0$
because 
$\dim \ \supp \Cal C \otimes \Cal O_l =0$ (note that $l$ is not contained in
the singular locus of $\iota (Y)$) and $p_* \Cal C=0$.
Hence by Nakayama's lemma, $\Cal C= 0$.
\qed
\enddemo

Hence $\iota : Y \to \iota(Y)$ is finite and birational and $\iota(Y)$
is normal, it is an isomorphism by the Zariski's Main Theorem.
Hence $X\simeq a(\iota(Y))$ is naturally embedded into 
$\Bbb P(1^h, 2^N)$ and $-K_X = \Cal O (1)$.

\definition{Case $\beta$}
Let $g':Y \to Z \to \Bbb P^3$ be 
the composition of $g$ and the double covering
of $\Bbb P^3$ branched along a sextic.
Consider the $\Bbb P^1$-bundle 
$\Bbb P(\Cal O_{\Bbb P^3} \oplus \Cal O_{\Bbb P^3} (-2))$
and denote it by $\Bbb P'$. 
Let $b': \Bbb P' \to \Bbb P^3$ be the natural projection and $T'$ 
the tautological divisor.
Note that 
by $1= (-K_Y)^2 E = ({g'}^* \Cal O (1))^2 E = (\Cal O (1))^2 {g'}_*E$,
we have $P:=g'(E) \simeq E \simeq \Bbb P^2$. 
As in the treatment of Case $\alpha$, 
we have a morphism $\iota ': Y\to \Bbb P'$
over $\Bbb P^3$ associated to the surjection
${g'}^*(\Cal O_{\Bbb P^3} \oplus \Cal O_{\Bbb P^3}(-2)) \to \Cal O_Y(E)$.
By the definition of $\iota '$, we have ${\iota '} ^* (T'|_{\iota '(Y)}) =E$.
Since $\deg \ {g'} = 2$, $\deg \ \iota ' = 1 \ \text{or} \ 2$.
But we know that $\deg \ \iota ' = 1$ 
by ${\iota '}^ * (T'|_{\iota '(Y)}) = E$.
\comment
Assume that $\deg \ \iota '= 2$.
Then $\deg b'|_{\iota '(Y)} =1$. Since $\iota '(Y)$ is contained in a 
$\Bbb P^1$-bundle, we know that 
$(b'|_{\iota '(Y)})^* P$ has exactly one component which dominates $P$. 
Clearly this component is $T'|_{\iota '(Y)}$.
Furthermore , ${g'}^* P$ has exactly one
component which dominates $P$.

a contradiction.
\endcomment
\comment
Since $\o{E}$ is not contained in the branch locus of $g'$
and $E\simeq \o{E}$, ${g'}^* \o{E}$ must have two component.
But ${g'}^* \o{E} = {\iota '}^*{T'|_{\iota(Y)}} = E$, a contradiction.
\endcomment
Hence $\text{deg} b' |_{\iota ' (Y)} =2$.
So we can write $\iota ' (Y) \sim 2T' + aL$, where
$L:= {b'}^*\Cal O (1)$ and $a$ is an integer.
By $K_{\Bbb P '} \sim -2T' -6L$, we have $K_{\iota ' (Y)}\sim
(a-6)L|_{\iota ' (Y)}$.
For a line $l$ in $T' | _{\iota ' (Y)}$, we have $K_{\iota ' (Y)}.l
=-1$.
So $a=5$, which in turn shows that $K_Y = {\iota '}^* K_{\iota ' (Y)}$.
Hence $\iota ' (Y)$ is normal since it is Gorenstein and 
\comment
Let $\iota '(Y) \to Z' \to \Bbb P^3$ is the Stein factorization of
$b'|_{\iota '(Y)}$.
Since $T'|_{\iota '(Y)} \simeq \Bbb P^2$, $\iota '(Y)$ is smooth at there.
Hence any exceptional curve for $\iota '(Y) \to Z'$
is not contained in the singular locus of $\iota '(Y)$
because every fiber of $\iota '(Y) \to Z'$ intersects $T'$.
In particular $\iota '(Y)$ satisfies $R_1$ since so does $Z'$. 
$S_2$ is clearly satisfied since
$\iota '(Y)$ is a Cartier divisor of $\Bbb P'$. Hence $\iota '(Y)$ 
is normal and
\endcomment
$Y\simeq \iota '(Y)$ by the Zariski's Main Theorem.

Contracting $T'$, $\Bbb P'$ is transformed into $\Bbb P(1^4, 2)$
and $Y$ is transformed into $X$. Hence we have the assertion.
\enddefinition

\definition{Case $\gamma$}
We will prove this case does not occur.
Let $g' : Y \to R$ be the composition of $g$ and the double covering
$\Phi _{|-K_X|}$. Then $g'(E)$ is a plane in $R$ by the same reason as 
in Case $\beta$. Hence we can assume that in $\Bbb P (1^5, 2)$
(note that $Z \subset \Bbb (1^5, 2)$ in this case),
$\o{E}$ is $(x_4= x_5 = y =0)$, where
$x_i$'s ($i=1,2,3,4,5$) are the coordinates of degree $1$ and $y$
is the coordinate of degree $2$.
So
the weighted equation of degree 2 of $Z$
is the form
$ay + x_4 l_1(x) + x_5 l_2(x)$, which in turn shows that
$Z \simeq ((4) \subset \Bbb P(1^4))$, a contradiction.
\enddefinition

\comment
Let $g:Y \to Z$ be the morphism defined by a sufficient multiple
of $-K_X$ and $\o{E}:=g(E)$. 

Since $\o{E}$ is a surface of degree $1$,
$\o{E}\simeq \Bbb P^2$. 
By [ibid.],
$Z$ is $(4) \subset \Bbb P (1^5)$ or 
$Z$ is 
In the former case,
the assertion follows from the argument of the first part of the proof.
In the latter case,
the assertion follows from the argument of the case that $h=4$.

\endcomment
This complete the induction.
 
Finally we describe the graded ring of $X$.

First we note that $|-2K_X|$ is free since
$-2K_X = \Cal O_X(2)$.
So we can take a smooth curve which is the intersection of
general members of $|-K_X|$ and $|-2K_X|$.
We fix such a curve and denote it by $C$ and $L:= -K_X|_C$.
Note that $L$ is a Cartier divisor such that $K_C = 2L$.
Since $-K_X = \Cal O(1)$, 
we may assume that $C\subset \Bbb P(1^{h-1}, 2^{N-1})$.
It suffices to describe the graded ring of $C$.
It is done by [RM4, Theorem 3.4].
Let $R(C, L):= \oplus _{m\geq 0} H^0 (\Cal O_C (mL))$.
Let $X_N \subset \Bbb P(1^h)$ be the image of the restriction
of the projection $\Bbb P(1^h , 2^N) \dashrightarrow \Bbb P(1^h )$.
The rational map $X\dashrightarrow X_N$ is a composition of
blow ups of $\frac 12(1,1,1)$-singularities and crepant contractions
in case $h\geq 5$ (resp. 
a composition of
blow ups of $\frac 12(1,1,1)$-singularities, crepant contractions
and the double covering of $\Bbb P^3$ in case $h=4$).
So the restriction of the projection to $C$ is a birational map in case
$h\geq 5$ (resp. a birational map or 
a double cover of a plane curve of degree $\geq 
3$ in case $h=4$),
which in turn show that the image of $C$ by the morphism of $L$ is not
a normal rational curve in $\Bbb P(1^{h-1})$.
Hence 
$H^0 (\Cal O_C (L)) \otimes H^0 (\Cal O_C (2L)) \to H^0 (\Cal O_C (3L))$
is surjective.
(Note that $K_C = 2L$.)
So by [RM4, Theorem 3.4], 
$R(C, L)$ is generated by elements of degree $\leq 2$
and related by elements of degree $\leq 6$,
which in turn show that the same things hold for
$\oplus _{m\geq 0} H^0 (\Cal O_X (-mK_X))$.
Let $N'$ be the number of sub bases of degree $2$
which do not come from degree $1$.
Since the above embedding $X\subset \Bbb P(1^h , 2^N)$
come from (possibly) some projection
$\Bbb (1^h , 2^{N'}) \dashrightarrow \Bbb P(1^h , 2^N)$,
$X$ is an intersection of weighted hypersurfaces of degree $\leq 6$. 

Finally we determine $X$ in $3$ cases as in the statement of this theorem.
It suffices to determine $C$ as above.
If $h=4$ and $N=1$, the assertion is clear.
Assume that $h=4$ and $N=2$.
If there is a relation of degree $2$ in $R(C, L)$,
the image of the restriction to $C$
of the projection $\Bbb P(1^4 , 2) \dashrightarrow \Bbb P(1^4 )$
is a conic in $\Bbb P^2$, a contradiction.
Hence there is no relation of degree $2$ in $R(C, L)$.
Then we find easily the relation of $R(C, L)$.
 
Assume that $h=5$ and $N=1$.
Note that we know that $C\subset \Bbb P^3$
and $\deg C = 9$. Hence if there is a relation of degree $2$ in $R(C, L)$,
there is exactly one relation and one degree $2$ base which does not come
from degree $1$. But from this, we can see that 
there is $2$ relation in degree $3$, a contradiction to that 
$\deg C = 9$. Hence there is no relation of degree $2$ in $R(C, L)$.
The rest are easy calculations.
Now we complete the proof of the main theorem. 
\qed
\enddemo

\definition{Remark}
The assumption that $h^0 (-K_X) \geq 4$ is necessary for Theorem 8.2
by the existence of the following:

$$(X\simeq (((12)\subset \Bbb P (1^3,4,6)))$$
which satisfies $h^0(-K_X) =3$.
\enddefinition

\proclaim{Corollary 8.3}
Let $X$ be a $\Bbb Q$-Fano $3$-fold as in the main theorem.
Assume that $X$ has only $\frac 12 (1,1,1)$-singularities as its
non Gorenstein points.
Then $X$ is embedded into $\Bbb P (1^h, 2^N)$
and $-K_X$ is the restriction of $\Cal O(1)$,
where $h:=h^0 (-K_X)$ and 
$N$ is the number of $\frac 12 (1,1,1)$-singularities.
Furthermore $X$ is an intersection of weighted hypersurfaces of degree
$\leq 6$.

If $h=4$ and $N=1$,
then $X \simeq ((5)\subset \Bbb P (1^4,2))$.

If $h=4$ and $N=2$,
then $X \simeq ((3,4)\subset \Bbb P (1^4,2^2))$.

If $h=5$ and $N=1$,
then $(X \simeq ((3,3)\subset \Bbb P (1^5,2))$.
\endproclaim

\demo{Proof}
By Corollary 6.3, Proposition 7.0 and Corollary 7.2, 
we can see that the assumptions of Theorem 8.2 are satisfied for $X$.
Hence we are done.
\qed
\enddemo

By this Corollary 8.3, we can improve Theorem 4.0 for $X$ as 
in Corollary 8.3 and Proposition 7.0 as follows:

\proclaim{Corollary 8.4}
Let $X$ be a $\Bbb Q$-Fano $3$-fold as in the main theorem.
Then 
\roster
\item
$-2K_X$ is very ample;
\item
$|-K_X|$ is free outside $\frac 12 (1,1,1)$-singularities and
its general member has only ordinary double points as its
singularities.
\endroster 
\endproclaim

\demo{Proof}
The proof is clear from Corollary 8.3.
\qed
\enddemo

\enddocument

\enddocument